\documentclass[preprint,10pt]{elsarticle}
\usepackage{hyperref}
\usepackage{amsmath,amssymb}
\usepackage {graphicx,fancyhdr}
\usepackage{graphics,color}
\usepackage{flafter}
\usepackage{multirow}
\usepackage{mathrsfs}
\usepackage{subfigure}
\usepackage{tabularx}
\usepackage{overpic}
\usepackage{stmaryrd}
\usepackage{times}
\usepackage{epstopdf}
\usepackage{caption}
\numberwithin{equation}{section}

\newtheorem{lemma}{Lemma}[section]
\newtheorem{theorem}{Theorem}[section]

\begin{document}
	
\begin{frontmatter}
	\tnotetext[label1]{This work is partially supported by Key Program Special Fund in XJTLU (KSF-E-50, KSF-P-02) and XJTLU Research Development Funding (RDF-19-01-15). The research of Eric Chung is partially supported by the Hong Kong RGC General Research Fund (Project numbers 14304719 and 14302018) and CUHK Faculty of Science Direct Grant 2019-20.}

\title{Adaptive generalized multiscale approximation of a mixed finite element method with velocity elimination \tnoteref{label1}}	

\author[lab1]{Zhengkang He}
\author[lab2]{Eric T. Chung}
\author[lab1,lab3]{Jie Chen\corref{cor1}}
\author[lab1,lab4]{Zhangxin Chen\corref{cor2}}

\cortext[cor1]{Corresponding author. Jie Chen E-mail address: jie.chen01@xjtlu.edu.cn}
\cortext[cor2]{Corresponding author. Zhangxin Chen E-mail address: zhachen@ucalgary.ca}

\cortext[lab1]{Zhengkang He E-mail address:  hzk2abc@163.com}
\cortext[lab2]{Eric T.Chung E-mail address:  tschung@math.cuhk.edu.hk.}

\address[lab1]{School of Mathematics and Statistics, Xi'an Jiaotong University, Xi'an, 710049, China}
\address[lab2]{Department of Mathematics, The Chinese University of Hong Kong (CUHK), Hong Kong SAR}
\address[lab3]{Department of Mathematical Sciences, Xi'an Jiaotong-Liverpool University, Suzhou, 215123, China}
\address[lab4]{Department of Chemical $\&$ Petroleum Engineering, Schulich School of Engineering, University of Calgary, 2500 University Drive N.W., Calgary, Alberta T2N 1N4, Canada}
	
\begin{abstract}
  In this paper, we propose offline and online adaptive enrichment algorithms for the generalized multiscale approximation of a mixed finite element method with velocity elimination to solve the subsurface flow problem in high-contrast and heterogeneous porous media. In the offline adaptive method, we first derive an a-posteriori error indicator based on one weighted $L^2$-norm of the local residual operator, where the weighted $L^2$-norm is related to the pressure fields of the local snapshot space. Then, we enrich the multiscale space by increasing the number of offline basis functions iteratively on coarse elements where the error indicator takes large values. While in the online adaptive method, we add online basis functions on selected coarse elements based on another weighted $L^2$-norm of the local residual operator to enrich the multiscale space, here the weighted $L^2$-norm is related to the velocity fields of the local snapshot space. Online basis functions are constructed in the online stage depending on the solution of the previous iteration and some optimal estimates. We give the theoretical analysis for the convergence of these two adaptive methods, which shows that sufficient initial basis functions (belong to the offline space) leads to a faster convergence rate. A series of numerical examples are provided to highlight the performance of both these two adaptive methods and also validate the theoretical analysis. Both offline and online adaptive methods are effective that can reduce the relative error substantially. In addition, the online adaptive method generally performs better than the offline adaptive method as online basis functions contain important global information such as distant effects that cannot be captured by offline basis functions. The numerical results also show that with a suitable initial multiscale space that includes all offline basis functions corresponding to relative smaller eigenvalues of local spectral decompositions in the offline stage, the convergence rate of the online enrichment is independent of the permeability contrast.
\end{abstract}
	
	
\end{frontmatter}
\section{Introduction}
In many real-world subsurface flow applications, such as petroleum recovery, groundwater resource management, and geothermal energy production, the geological porous-media are usually governed by coefficients with high heterogeneities and complex spatial distributions. Solving these problems directly on the fine grid will result in large-scale discrete systems which are challenging to tackle. Multiscale model reduction techniques such as upscaling techniques \cite{durlofsky1991numerical,wu2002analysis,gao2015numerical} and multiscale methods \cite{chu2010new,chung2010reduced,chung2011energy,chung2014generalized,chung2013sub,efendiev2011multiscale,efendiev2009multiscale,efendiev2004multiscale,ghommem2013mode} are often employed to reduce the computational complexity.

For the multiscale methods, multiscale basis functions are constructed locally to capture the local multiscale information of the fine-grid solution. There have been many works \cite{efendiev2006accurate,efendiev2000convergence,owhadi2007metric} proposed to optimize the appropriate number of multiscale basis functions for improving the accuracy of the multiscale solution. In works \cite{efendiev2013generalized_1,efendiev2013generalized_2,efendiev2014generalized,chung2014generalized} of last several years, the authors have developed a flexible framework, generalized multiscale finite elements methods (GMsFEM), that generalizes the multiscale finite element method \cite{hou1997multiscale} by adding additional multiscale basis function that can capture extra local multiscale information to enrich the multiscale space. The computational procedure of GMsFEM is separated into two stages: the offline stage and the online stage. In the offline stage, a small dimensional offline space is established through solving a series of local problems and by conducting some suitable spectral decompositions. In the online stage, the offline multiscale space is used to compute the multiscale basis functions for the construction of the multiscale space, then based on the multiscale space, the multiscale solutions are derived by solving the problem on the coarse grid. There have also been a mixed generalized multiscale finite element method (mixed GMsFEM) \cite{chung2015mixed} and a generalized multiscale discontinuous Galerkin method (GMsDGM) \cite{chung2018adaptive} presented for high-contrast flow problems following the GMsFEM framework, where the multiscale basis functions are coupled through the standard mixed finite element method and symmetric interior penalty discontinuous Galerkin method, respectively, on the coarse grid.

In the framework of GMsFEM, adaptive enrichment of the multiscale space is of great importance as it identifies regions where the local multiscale basis functions in use are not adequate to capture all dominated local multiscale information. Therefore, the offline adaptive GMsFEM \cite{chung2014adaptive,chan2016adaptive,chung2018adaptive,chung2019adaptive} is developed subsequently. In the offline adaptive method, the multiscale space is enriched by the use of offline basis functions iteratively on coarse elements according to an a-posteriori error indicator. The error indicator is developed by an appropriate norm of the local residual operator together with the eigenvalue structure of the spectral decomposition. However, offline basis functions only contain the multiscale information locally, so after adding some offline basis functions such that the multiscale space is capable of capturing all dominated local multiscale information, some global information needs to be taken into consideration as the distant effects can be significant. Consequently, in later works \cite{chung2015residual,chan2016adaptive,chung2017online}, the online basis functions are introduced and the multiscale space is enriched adaptively with these online basis functions. Both offline and online adaptive enrichments can substantially accelerate the convergence of GMsFEM.

Recently, a generalized multiscale approximation of a mixed finite element method (MFEM) with velocity elimination has been developed in \cite{chen2020generalized} for the subsurface flow problem, which also follows the GMsFEM framework. Different from the mixed GMsFEM developed in \cite{chung2015mixed} where the multiscale basis functions are constructed for the approximation of velocity and piecewise constant on the coarse grid is used for pressure, the method in \cite{chen2020generalized} constructs multiscale basis functions to approximate the pressure and makes use of the trapezoidal quadrature rule for local velocity elimination, i.e., the velocity is solved directly on the fine grid, and in the end, only a symmetric and positive definite system related to the multiscale pressure need to be solved. Both of these two mixed GMsFEMs are locally mass conservative on the coarse-grid scale. One only needs to conduct simple post-processing in the region of interest to obtain the local mass conservation in the fine-grid scale. The local mass conservation is essential in many subsurface flow applications, especially when the flow is coupled with the transport or the streamline needs to be constructed.

In this paper, following the overall idea of the adaptive multiscale model reduction with GMsFEM \cite{chung2014adaptive,chan2016adaptive,chung2018adaptive,chung2015residual,chung2017online,chung2016adaptive}, we develop efficient offline and online adaptive enrichment algorithm, respectively, for the generalized multiscale approximation of a MFEM with velocity elimination. The offline adaptive method enriches the multiscale space with offline basis functions that are precomputed in the offline stage before the enrichment algorithm and will be used in the online stage for any given source terms and boundary conditions, while the online adaptive method enriches the multiscale space with online basis functions which need to be calculated in the actual simulation, belonging to the online stage. We introduce two different weighted $L^2$-norms of the local residual operator on the local snapshot space, one is related to the pressure fields of the local snapshot space and the other is associated with the velocity fields of the local snapshot space. For simplicity, we will call them pressure-related and velocity-related weighted $L^2$-norm, respectively. In the offline adaptive method, we employ the error indicator based on the pressure-related weighted $L^2$-norm to select coarse elements where offline basis functions need to be added, where the eigenvalue structures of spectral decompositions in the offline stage are also coupled into the error indicator. In the online adaptive method, all coarse elements are separated into non-overlapping subsets, and at each iteration of the multiscale space enrichment, the online basis functions are constructed and added on coarse elements of each these subsets based on the velocity-related weighted $L^2$-norm of the local residual operator. Moreover, we give the corresponding convergence analysis of these two adaptive methods, respectively. In our analysis, some stability and approximation properties of several projection operators from the local snapshot space to the local offline space are obtained by utilizing spectral estimates of the spectral decomposition in the offline stage, from these properties, it can be shown that the error of the multiscale solution is bounded by the proposed a-posteriori error indicator. We shows that the initial multiscale space equipped with sufficient initial basis functions leads to a faster convergence rate. In the end, we also conduct ample numerical tests to confirm the theoretical analysis and show the convergence behaviour of the proposed two adaptive enrichment algorithms. The numerical results show that both the offline and online adaptive methods are effective, reliable and can achieve a substantial error reduction of the multiscale solution. The error decay of the online adaptive method is usually more quickly than the error decay of the offline adaptive method, so the online adaptive method generally performs better than the offline adaptive method. Furthermore, when the initial multiscale space contains all offline basis functions corresponding to the small eigenvalues (which are contrast sensitive) of the local spectral decompositions on all coarse elements, the convergence rate of the online adaptive method is independent of the contrast of the permeability field.

We organize the rest of the paper as follows. In section 2, we review the generalized multiscale approximation of a MFEM with velocity elimination. In section 3 and section 4, details of offline and online adaptive enrichment algorithms are presented, respectively, the corresponding convergence analyses are also given. In section 5, a series of numerical examples are shown to illustrate the convergence behaviour and verify the convergence analysis of these two adaptive enrichment algorithms. Finally, some conclusions are presented in section 6.
\section{Preliminaries}
In this section, we briefly introduce the generalized multiscale approximation of a mixed finite element method with velocity elimination proposed in \cite{chen2020generalized} for the single-phase flow problem. First, the governing equations together with the fine-grid approximation, and then the generalized multiscale approximation are explained.
\subsection{Model problem and fine-grid approximation}
Let $\Omega$ be a bounded and simply connected porous-media domain in $\mathbb{R}^2$ with a Lipschitz continuous boundary $\partial\Omega$. We consider the single-phase flow described by Darcy's law and a mass conservation equation as
\begin{eqnarray}
{\kappa}^{-1}\mathbf{u} + \nabla{p} = \mathbf{0} &&\textrm{in}\ \ \Omega,\label{eqn_model_1}\\
\nabla\cdot\mathbf{u} = f &&\textrm{in}\ \ \Omega,\label{eqn_model_2}
\end{eqnarray}
with following boundary conditions on $\partial\Omega$
\begin{eqnarray*}
p = g_D &&\textrm{on}\ \ \partial\Omega_D,\\
\mathbf{u}\cdot\mathbf{n} = g_N &&\textrm{on}\ \ \partial\Omega_N,
\end{eqnarray*}
where $\kappa$ is the high-contrast and heterogeneous permeability, $\mathbf{n}$ is the unit outward normal vector on $\partial\Omega$, $\partial\Omega_D$ and $\partial\Omega_N$ are the Dirichlet and Neumann boundaries, respectively, with the corresponding boundary data $g_D$ and $g_N$. For simplicity, we set $g_N=0$ in this paper, although more general boundary conditions can also be treated.

Standard notations and definitions for Sobolev spaces are used for the weak formulation of the problem (\ref{eqn_model_1})-(\ref{eqn_model_2}). Define the following two spaces
\begin{equation*}
{V} = \{\mathbf{v}\in H(\textrm{div},\Omega): \ \mathbf{v}\cdot\mathbf{n}=0\ \ \textrm{on}\ \ \partial\Omega_N\} \quad\textrm{and}\quad W = L^2(\Omega).
\end{equation*}
where $ H(\textrm{div},\Omega) = \{\mathbf{v}\in({L^2(\Omega)})^2: \ \nabla\cdot\mathbf{v}\in{L^2(\Omega)}\}$. Then the weak formulation for (\ref{eqn_model_1})-(\ref{eqn_model_2}) can be written as: find $(\mathbf{u},p)\in{V\times{W}}$, such that
\begin{eqnarray}
(\kappa^{-1}\mathbf{u},\mathbf{v}) - (p,\nabla\cdot\mathbf{v}) = -(g_D,\mathbf{v}\cdot\mathbf{n})_{\partial\Omega_D}, &&\forall\mathbf{v}\in V,\label{eqn_weak_1}\\
-(\nabla\cdot\mathbf{u},q) = -(f,q), \hspace{1.35cm}&&\forall q\in W.\label{eqn_weak_2}
\end{eqnarray}

\begin{figure}
	\centering\footnotesize
	\begin{overpic}[height=6.1cm,width=13cm]{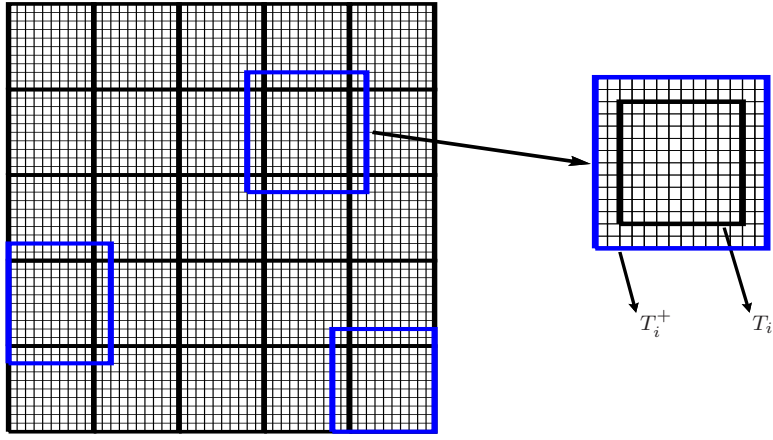}
		\put(79,12){$T_i$}
		\put(67.5,12){$T^+_i$}
	\end{overpic}
	\caption{The illustration of a multiscale mesh in the left and an oversampling coarse block $T^+_i$ associated with a coarse element $T_i$ in the right.}
	\label{Fig_mesh_Oversampling}
\end{figure}

To describe the framework of mixed GMsFEM, which solves the problem on two meshes with different scales, we construct the fine grid and coarse grid as follows. We assume the fine grid $\mathcal{T}_h$ is a uniform, regular partition of $\Omega$ composed of rectangles with mesh size $h$. The set of all fine-grid edges of $\mathcal{T}_h$ is denoted by $\mathcal{E}_h$.  We use $\mathcal{E}^0_h$ to denote the set of all interior fine-grid edges, and use $\mathcal{E}^N_h$ and $\mathcal{E}^D_h$ to denote the set of all fine-grid edges located on $\partial\Omega_N$ and $\partial\Omega_D$, respectively. For the coarse grid $\mathcal{T}_H$, each coarse element $T$ in $\mathcal{T}_H$ is defined as a connected collection of fine-grid elements belonging to $\mathcal{T}_h$, i.e., the coarse element numbered $i$ is given as $T_i = \cup^{N_i}_{k=1}t_k$, where $N_i$ is the number of fine-grid elements contained in $T_i$. In the simplest case that the coarse grid is formed as a uniform partition of the fine grid, then each coarse element $T$ turns into a rectangle, see Figure \ref{Fig_mesh_Oversampling} for an example of a multiscale mesh and a coarse element $T_i$. We use $N_T$ to denote the total number of coarse elements in $\mathcal{T}_H$.

Let $\hat{t}$ be the reference unit square $[0,1]^2$ with vertices $\hat{\mathbf{r}}_1=(0,0)^T$, $\hat{\mathbf{r}}_2=(1,0)^T$, $\hat{\mathbf{r}}_3=(0,1)^T$ and $\hat{\mathbf{r}}_4=(1,1)^T$, and let $t$ be any physical element in $\mathcal{T}_h$ with vertices $\mathbf{r}_i=(x_i,y_i)^T, i=1,...,4$, then there exists a bijective mapping $F_t:\hat{t}\to t$ as defined in (\ref{eqn_bi_mapping}). 
Denote Jacobian matrix of $F_t$ by $DF_t$, determinant of $DF_t$ by $J_t=|\textrm{det}({DF}_t)|$ and inverse mapping of $F_t$ by ${F_t^{-1}}$, respectively.
\begin{equation}\label{eqn_bi_mapping}
F_t(\hat{x},\hat{y})=\mathbf{r}_1(1-\hat{x})(1-\hat{y})+\mathbf{r}_2\hat{x}(1-\hat{y})+\mathbf{r}_3\hat{x}\hat{y}+\mathbf{r}_4(1-\hat{x})\hat{y}.
\end{equation}
Using the bijective mapping $F_t$, for any scalar function $\hat{w}$ defined in $\hat{t}$, we denote the scalar transformation of $\hat{w}$ in $t$ by $w$, defined as $w\leftrightarrow\hat{w}:w=\hat{w}\circ F_t^{-1}$.

In this paper, for simplicity, we assume that all elements in $\mathcal{T}_h$ are squares with the same size $h$. Then the Jacobian matrix $DF_t$ and the determinant $J_t$ are reduced into
\begin{equation}\label{eqn_bi_mapping_square}
DF_t = \textrm{diag}(h,h)\quad \textrm{and} \quad J_t = h^2,
\end{equation}
where $\textrm{diag}(h,h)$ denotes a $2\times 2$ diagonal matrix with elements of $[h,h]$ on the main diagonal.

On the reference unit square $\hat{t}$, the lowest order Raviart-Thomas space is defined as
$\textrm{RT}_0(\hat{t}) = \left(
\begin{matrix}
\alpha_1 + \beta_1\hat{x}\\
\alpha_2 + \beta_2\hat{y}
\end{matrix}\right)$ with ${\alpha_i}|_{i=1,2}, {\beta_i}|_{i=1,2}\in\mathbb{R}$ are arbitrary constants. On any element $t\in \mathcal{T}_h$, the local space $\textrm{RT}_0(t)$ is defined via the following vector transformation,
\begin{equation}\label{eqn_Piola_transformation}
\mathbf{v}\leftrightarrow\hat{\mathbf{v}}:\mathbf{v} = \frac{1}{J_t}DF_t\hat{\mathbf{v}}\circ F_t^{-1} = \frac{1}{h}\hat{\mathbf{v}}\circ F_t^{-1},\quad \forall\hat{\mathbf{v}}\in\textrm{RT}_0(\hat{t}),
\end{equation}
which is known as the Piola transformation, preserving the normal component of the vector on each edge, i.e., $\int_e\mathbf{v}\cdot\mathbf{n}_e=\int_{\hat{e}}\hat{\mathbf{v}}\cdot\hat{\mathbf{n}}_{\hat{e}}$, $\forall e\in\partial t$, here $\mathbf{n}_{e}$ and $\hat{\mathbf{n}}_{\hat{e}}$ are unit outward normal vectors on the edge $e\in\partial t$ and the corresponding reference edge $\hat{e}\in\hat{t}$, respectively. Therefore, by the use of $\textrm{RT}_0$ space, we have that $\mathbf{v}\cdot\mathbf{n}_e|e|=\hat{\mathbf{v}}\cdot\hat{\mathbf{n}}_{\hat{e}}|\hat{e}|=\hat{\mathbf{v}}\cdot\hat{\mathbf{n}}_{\hat{e}}$.

We make use of the following $\textrm{RT}_0$ mixed finite element spaces for the approximation of velocity and pressure to discretize the weak formulation (\ref{eqn_weak_1})-(\ref{eqn_weak_2}),
\begin{equation*}
V_h = \{\mathbf{v}\in V: \mathbf{v}|_{t}\in\textrm{RT}_0(t), \ \ \forall t\in\mathcal{T}_h \} \quad\textrm{and}\quad
W_h = \{q\in W: q|_{t}\in\mathcal{P}_0, \ \ \forall t\in\mathcal{T}_h \},
\end{equation*}
where $\mathcal{P}_0$ denotes the polynomial space of zero degree, and we obtain the corresponding discrete weak formulation as: find $(\mathbf{u}_h,p_h)\in{V_h\times{W}_h}$, such that
\begin{eqnarray}
(\kappa^{-1}\mathbf{u}_h,\mathbf{v}) - (p_h,\nabla\cdot\mathbf{v}) = -(g_D,\mathbf{v}\cdot\mathbf{n})_{\partial\Omega_D}, &&\forall\mathbf{v}\in V_h,\label{eqn_discrete_1}\\
-(\nabla\cdot\mathbf{u}_h,q) = -(f,q), \hspace{1.35cm}&&\forall q\in W_h.\label{eqn_discrete_2}
\end{eqnarray}
It is well known that the above velocity-pressure system (\ref{eqn_discrete_1})-(\ref{eqn_discrete_2}) is in a saddle-point structure, which is computational expensive. To avoid tackling the saddle-point algebraic system, we apply the trapezoidal quadrature rule that allows for local velocity elimination and results in a symmetric and positive definite algebraic system for pressure.

Next, we explain how to apply the trapezoidal quadrature rule \cite{russell1983finite,klausen2006robust,wheeler2006multipoint,hou2014numerical} to compute the integration $(\kappa^{-1}\mathbf{w},\mathbf{v})$, for any $\mathbf{w},\mathbf{v}\in V_h$. By use of the bilinear mapping (\ref{eqn_bi_mapping}) and the Piola transformation (\ref{eqn_Piola_transformation}), the integration on any physical element $t\in\mathcal{T}_h$ is mapped to the reference element $\hat{t}$, that is,
\begin{equation*}
(\kappa^{-1}\mathbf{w},\mathbf{v})_t = (\frac{1}{J_t}DF_t^T\hat{\kappa}^{-1}DF_t\hat{\mathbf{w}},\hat{\mathbf{v}})_{\hat{t}}=(\widehat{\mathcal{M}}_t\hat{\mathbf{w}},\hat{\mathbf{v}})_{\hat{t}},
\end{equation*}
where $\hat{\kappa}=\kappa\circ F_t$, $\widehat{\mathcal{M}}_t=DF_t^T\hat{\kappa}^{-1}DF_t/{J_t}$, and $\hat{\mathbf{w}}$, $\hat{\mathbf{v}}\in\textrm{RT}_0(\hat{t})$ are the inverse functions of $\mathbf{w}$, $\mathbf{v}$ through the Piola transformation (\ref{eqn_Piola_transformation}), respectively. By applying the trapezoidal quadrature rule on the reference element $\hat{t}$, we get the quadrature rule on the physical element $t\in\mathcal{T}_h$ as
\begin{equation*}
(\kappa^{-1}\mathbf{w},\mathbf{v})_{Q,t}=(\widehat{\mathcal{M}}_t\hat{\mathbf{w}},\hat{\mathbf{v}})_{\hat{Q},\hat{t}}=
\frac{|\hat{t}|}{4}\sum\limits_{i=1}^{4}{\widehat{{\mathcal{M}}}_t(\hat{\mathbf{r}}_i)\hat{\mathbf{w}}(\hat{\mathbf{r}}_i)\cdot\hat{\mathbf{v}}(\hat{\mathbf{r}}_i)}.
\end{equation*}
In the following of this paper, we suppose $\{\mathbf{v}_1, \mathbf{v}_2, \cdots, \mathbf{v}_{N_e}\}$ is a set of basis functions of $V_h$, transformed through the Piola transformation (\ref{eqn_Piola_transformation}) with the reference basis functions satisfying $\hat{\mathbf{v}}_{\hat{e}}\cdot\hat{\mathbf{n}}_{\hat{e}}=1$, where $N_e$ is the number of all fine-grid edges in $\mathcal{E}_h$, and $\mathbf{w},\mathbf{v}\in V_h$ can be expressed as $\mathbf{w} = \sum_{e\in\mathcal{E}_h}w_e\mathbf{v}_e$, $\mathbf{v} = \sum_{e\in\mathcal{E}_h}v_e\mathbf{v}_e$, respectively. Equation (\ref{eqn_bi_mapping_square}) gives that $\widehat{\mathcal{M}}_t=DF_t^T\hat{\kappa}^{-1}DF_t/{J_t} = \kappa^{-1}_t$, then we obtain
\begin{equation*}
(\kappa^{-1}\mathbf{w},\mathbf{v})_{Q,t}=
\frac{1}{4}\sum\limits_{i=1}^{4}{\kappa^{-1}_t\hat{\mathbf{w}}(\hat{\mathbf{r}}_i)\cdot\hat{\mathbf{v}}(\hat{\mathbf{r}}_i)}=
\frac{1}{2}\sum\limits_{e\in\partial t}{\kappa^{-1}_tw_ev_e},
\end{equation*}
and correspondingly, the global quadrature rule for the integration $(\kappa^{-1}\mathbf{w},\mathbf{v})$ in the whole domain $\Omega$ is defined as
\begin{equation}\label{eqn_global_quad_rule}
(\kappa^{-1}\mathbf{w},\mathbf{v})_Q = \sum\limits_{t\in\mathcal{T}_h}{(\kappa^{-1}\mathbf{w},\mathbf{v})_{Q,t}}
=\frac{1}{2}\sum\limits_{t\in\mathcal{T}_h}\sum\limits_{e\in\partial t}{\kappa^{-1}_tw_ev_e}.
\end{equation}
By the above quadrature rule, we define the related norm on the space $V_h$ as $\|\mathbf{v}\|_{\kappa^{-1}} = (\kappa^{-1}\mathbf{v},\mathbf{v})^{\frac{1}{2}}_Q$, $\forall\mathbf{v}\in V_h$, and from \cite{wheeler2006multipoint}, we know that the norm $\|\mathbf{v}\|_{\kappa^{-1}}$ is equivalent with the $L^2$-norm $\|\mathbf{v}\|$.

We obtain the corresponding discrete weak formulation using the above quadrature rule (\ref{eqn_global_quad_rule}): find $(\mathbf{u}_h,p_h)\in{V_h\times{W}_h}$, such that
\begin{eqnarray}
(\kappa^{-1}\mathbf{u}_h,\mathbf{v})_Q - (p_h,\nabla\cdot\mathbf{v}) = -(g_D,\mathbf{v}\cdot\mathbf{n})_{\partial\Omega_D}, &&\forall\mathbf{v}\in V_h,\label{eqn_mfmfe_1}\\
-(\nabla\cdot\mathbf{u}_h,q) = -(f,q), \hspace{1.35cm}&&\forall q\in W_h.\label{eqn_mfmfe_2}
\end{eqnarray}

Suppose the dimensions of $V_h$ and $W_h$ are $m_1$ and $m_2$, respectively, then the above discrete system (\ref{eqn_mfmfe_1})-(\ref{eqn_mfmfe_2}) in the mixed formulation can be written into a matrix form as: find $({U}_h,{P}_h)\in\mathbb{R}^{m_1}\times\mathbb{R}^{m_2}$, such that
\begin{equation}\label{eqn_matrix_fine}
\left(
\begin{matrix}
A_h &B_h\\B^T_h &0
\end{matrix}\right)
\left(
\begin{matrix}
U_h\\P_h
\end{matrix}\right)=
\left(
\begin{matrix}
G_h\\F_h
\end{matrix}\right),
\end{equation}
where matrixs $A_h$, $B_h$, $G_h$ and $F_h$ are associated with terms $(\kappa^{-1}\mathbf{u}_h,\mathbf{v})_Q$, $-(p_h,\nabla\cdot\mathbf{v})$, $-(g_D,\mathbf{v}\cdot\mathbf{n})_{\partial\Omega_D}$ and $-(f,q)$, respectively. By use of the $\textrm{RT}_0$ mixed finite element spaces and the trapezoidal quadrature rule $(\cdot,\cdot)_Q$, from (\ref{eqn_global_quad_rule}), we know that $A_h$ is a diagonal matrix with positive diagonal elements, so that $A_h$ is inverted easily and we can solve the system (\ref{eqn_matrix_fine}) in the following way
\begin{equation}\label{eqn_matrix_form}
-B^T_h({A_h})^{-1}B_h P_h = F_h-B^T_h({A_h})^{-1}G_h,
\end{equation}
that is, we only need to solve a symmetric and positive definite system for pressure.

\begin{figure}[h!]
	\centering\footnotesize
	\begin{overpic}[height=1.8cm,width=3.6cm]{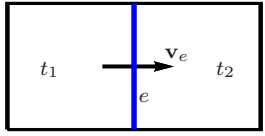}
		\put(15,22){$t_1$}
		\put(80,22){$t_2$}
		\put(51.5,12){$e$}
		\put(61,28){$\mathbf{v}_e$}
	\end{overpic}
	\caption{Two elements $t_1$ and $t_2$ in $\mathcal{T}_h$ sharing with an edge $e$ in $\mathcal{E}_h$, $\mathbf{v}_e$ is a basis functions for velocity associated with $e$}
	\label{Fig_MFEM_two_elements}
\end{figure}

Accordingly, in the following, we will define the corresponding discrete weak formulation with respect to the above linear system (\ref{eqn_matrix_form}). we use $a(\cdot,\cdot)$ to denote the bilinear form relating to the matrix $B^T_h({A_h})^{-1}B_h$ in the left-hand side of (\ref{eqn_matrix_form}). For any two elements $t_1$, $t_2\in\mathcal{T}_h$, sharing with the same edge $e$ in $\mathcal{E}_h$, as shown in Figure \ref{Fig_MFEM_two_elements}, we denote $\mathbf{v}_e$ be the basis function of $V_h$ associated with the edge $e$. Let $\mathbf{v} = \mathbf{v}_e$ in (\ref{eqn_mfmfe_1}), then the first term in the left-hand side of (\ref{eqn_mfmfe_1}) becomes
\begin{equation}\label{left_hand_side}
\begin{split}
(\kappa^{-1}\mathbf{u}_h,\mathbf{v}_e)_Q &= (\kappa^{-1}\mathbf{u}_h,\mathbf{v}_e)_{Q,t_1} + (\kappa^{-1}\mathbf{u}_h,\mathbf{v}_e)_{Q,t_2} \\
&= \frac{1}{2}{\kappa^{-1}_1u_e} + \frac{1}{2}{\kappa^{-1}_2u_e} = \bar{\kappa}^{-1}_eu_e,
\end{split}
\end{equation}
here $\bar{\kappa}_e = 2/(\kappa^{-1}_1+\kappa^{-1}_2)$ is the harmonic average of $\kappa_1$ and $\kappa_2$, with $\kappa_1$ and $\kappa_2$ are the permeability $\kappa$ on $t_1$ and $t_2$, respectively, and the second term in the right-hand side of (\ref{eqn_mfmfe_1}) turns into
\begin{equation}\label{right_hand_side}
\begin{split}
(p_h,\nabla\cdot\mathbf{v}_e) &= (p_h,\nabla\cdot\mathbf{v}_e)_{t_1} + (p_h,\nabla\cdot\mathbf{v}_e)_{t_2} \\
&= (p_h,\mathbf{v}_e\cdot\mathbf{n}_e)_e + (p_h,\mathbf{v}_e\cdot\mathbf{n}_e)_e = p_1 - p_2.
\end{split}
\end{equation}
Denote $\llbracket\cdot \rrbracket_e$ as the jump operator across edge $e$, such that $\llbracket p_h \rrbracket_e = p_1-p_2$. Combining (\ref{left_hand_side}) and (\ref{right_hand_side}), we have that
\begin{equation}\label{eqn_velocity_pressure_relationship}
\bar{\kappa}^{-1}_eu_e = \llbracket p_h \rrbracket_e, \quad\textrm{i.e.},\quad u_e = \bar{\kappa}_e\llbracket p_h \rrbracket_e,\quad\forall e\in\mathcal{E}^0_h.
\end{equation}
Note that when $e\in\mathcal{E}^D_h$, $\llbracket p_h \rrbracket_e = p_h|_t-\overline{g_D}|_e$, where $t$ is the fine-grid element such that $e\in\partial t$, and $\overline{g_D}|_e$ is the average value of function $g_D$ on edge $e$.
\\
For any fine-grid edge $e\in\mathcal{E}_h$ and element $t\in\mathcal{T}_h$, satisfying $e\in\partial t$, we define a scalar $s_{et} = \mathbf{v}_e\cdot\mathbf{n}_{et}|e|$, where $\mathbf{n}_{et}$ is the unit normal vector on $e$ pointing outward of $t$. For example, in the situation of Figure \ref{Fig_MFEM_two_elements}, $s_{et_1} = 1$ and $s_{et_2} = -1$. By the relationship in (\ref{eqn_velocity_pressure_relationship}), the term on the left-hand side of (\ref{eqn_mfmfe_2}) becomes
\begin{equation*}
\begin{split}
-(\nabla\cdot\mathbf{u}_h,q) &=-\sum\limits_{t\in\mathcal{T}_h}(\nabla\cdot\mathbf{u}_h,q)_t = -\sum\limits_{t\in\mathcal{T}_h}\sum\limits_{e\in\partial t}u_eq_ts_{et}\\
&= -\sum\limits_{t\in\mathcal{T}_h}\sum\limits_{e\in\partial t}\bar{\kappa}_e\llbracket p_h \rrbracket_eq_ts_{et} = -\sum\limits_{e\in\mathcal{E}^0_h\cup\mathcal{E}^D_h}\bar{\kappa}_e\llbracket p_h \rrbracket_e \llbracket q \rrbracket_e.
\end{split}
\end{equation*}
Thus, we define the bilinear form $a(\cdot,\cdot)$ as
\begin{equation}\label{eqn_definition_apq}
a(r,q) = \sum\limits_{e\in\mathcal{E}^0_h\cup\mathcal{E}^D_h}\bar{\kappa}_e\llbracket r \rrbracket_e \llbracket q \rrbracket_e, \quad\forall r,q\in W_h,
\end{equation}
and the original discrete weak formulation in the mixed form (\ref{eqn_mfmfe_1})-(\ref{eqn_mfmfe_2}) can be turned into the following discrete weak formulation as: find $p_h\in W_h$, such that
\begin{equation}\label{eqn_pressure_form}
a(p_h,q) = (f,q), \quad \forall q\in W_h,
\end{equation}
where the bilinear form $a(\cdot,\cdot)$ in (\ref{eqn_definition_apq}) is symmetric, continuous and coercive. Meanwhile, the matrix form of the above discrete weak formulation (\ref{eqn_pressure_form}) is the same as the linear system (\ref{eqn_matrix_form}).\\

\noindent\textbf{Remark 2.1.} By the definition of the harmonic average of permeability following (\ref{left_hand_side}), and the relationship (\ref{eqn_velocity_pressure_relationship}) between the velocity and pressure variables, we can write the bilinear form $(\kappa^{-1}\cdot,\cdot)_Q$ as
\begin{equation*}
(\kappa^{-1}\mathbf{w},\mathbf{v})_Q = \frac{1}{2}\sum\limits_{t\in\mathcal{T}_h}\sum\limits_{e\in\partial t}{\kappa^{-1}_tw_ev_e}
= \sum\limits_{e\in\mathcal{E}^0_h\cup\mathcal{E}^D_h}\bar{\kappa}^{-1}_ew_ev_e
= \sum\limits_{e\in\mathcal{E}^0_h\cup\mathcal{E}^D_h}\bar{\kappa}_e\llbracket r \rrbracket_e \llbracket q \rrbracket_e,
\end{equation*}
where $\mathbf{w}$, $\mathbf{v}\in V_h$ are the velocity fields with respect to $r$, $q\in W_h$, respectively, by the relationship (\ref{eqn_velocity_pressure_relationship}), then from the definition of $a(\cdot,\cdot)$ in (\ref{eqn_definition_apq}), we have
\begin{equation}\label{eqn_bilinear_form_equal}
(\kappa^{-1}\mathbf{w},\mathbf{v})_Q = a(r,q).
\end{equation}
On the coarse element $T_i\in\mathcal{T}_H$, $i=1,2,\cdots,N_T$, we define following three local norms for the velocity variable as
\begin{equation}\label{eqn_three_local_norms}
\|\mathbf{v}\|^2_{\kappa^{-1},i^0} = (\kappa^{-1}\mathbf{v},\mathbf{v})_{Q,{T}^0_i}, \ \
\|\mathbf{v}\|^2_{\kappa^{-1},i^+} = (\kappa^{-1}\mathbf{v},\mathbf{v})_{Q,{T}^+_i}, \ \
\|\mathbf{v}\|^2_{\kappa^{-1},i} = (\kappa^{-1}\mathbf{v},\mathbf{v})_{Q,{T}^0_i} + \sum\limits_{e\in\partial T_i}\bar{\kappa}_e\llbracket q \rrbracket^2_e,
\end{equation}
where $(\kappa^{-1}\mathbf{w},\mathbf{v})_{Q,{T}^0_i} = \sum_{e\in\mathcal{E}^0_{T_i}}\llbracket r \rrbracket_e\llbracket q \rrbracket_e$ with $\mathcal{E}^0_{T_i}$ denoting the set of all interior fine-grid edges in the partition for $T_i$, and note that in the definition of $\|\mathbf{v}\|_{\kappa^{-1},i}$, we assume zero values for $q$ outside of $T_i$, $i = 1,2,\cdots,N_T$.
\subsection{Generalized multiscale approximation}
In this subsection, we illustrate the multiscale method described in \cite{chen2020generalized} for solving the single-phase flow problem (\ref{eqn_model_1})-(\ref{eqn_model_2}) on the coarse grid, which follows the GMsFEM framework to compute the multiscale basis functions for pressure. We first derive the local snapshot space on each coarse element by solving a series of local problems with different boundary conditions, the snapshot space provided a solution space on each coarse element locally. Then we perform a spectral decomposition in each local snapshot space to obtain the dominant modes (offline basis functions) of the snapshot basis functions and get the corresponding local offline space with a smaller dimension. In addition, as the multiscale finite volume method \cite{jenny2003multi,lunati2006multiscale,wolfsteiner2006well,hajibeygi2008iterative,lunati2009operator}, the correction function is also introduced to consistently deal with the source term.
\subsubsection{Snapshot space}
We make use of oversampling techniques proposed in \cite{efendiev2014generalized} to get the more effective local snapshot space. Let $T_i \in \mathcal{T}_H$ be a coarse element in $\Omega$ and $T^+_i$ be a coarse block defined by adding some fine-grid layers around $T_i$, such that $T_i\subset T^+_i$, as shown in the right graph of Figure \ref{Fig_mesh_Oversampling}. Basis functions of the local snapshot space $W^{i,+}_{\textrm{snap}}$ are derived by numerically solving the following problems on the oversampling coarse block $T^+_i$: find $(\mathbf{\psi}_j^{i,+}, \phi_j^{i,+})$, such that
\begin{align}\label{eqn_local_problem}
\begin{split}
\kappa^{-1}\mathbf{\psi}_j^{i,+} + \nabla \phi_j^{i,+} &= 0 \qquad \mathrm{in} \ \ T^+_i,\\
\nabla\cdot\mathbf{\psi}_j^{i,+} & = 0 \qquad \mathrm{in} \ \ T^+_i.
\end{split}
\end{align}
Suppose the boundary of coarse block $T^+_i$ can be written as a union of fine-grid edges, i.e., $\partial T^+_i = \bigcup_{j = 1}^{J_i} e_j $, where $J_i$ is the total number of fine-grid edges on $\partial T^+_i$. Let $\delta^{i,+}_j $ be a piecewise constant function defined on $\partial T^+_i$ with respect to the fine-grid edges such that it has value $1$ on $e_j$ and value $0$ on the other fine-grid edges, that is
\begin{align*}
\begin{split}
\delta^{i,+}_j = \left\{
\begin{array}{ll}
1 & \text{on\ }\ e_j, \\
0 & \text{on\ other\ fine-grid\ edges\ on\ } \partial T^+_i,
\end{array}\right.
\qquad j = 1, 2, \cdots, J_i.
\end{split}
\end{align*}
The boundary conditions on the boundary of coarse block $T^+_i$ for the local problem (\ref{eqn_local_problem}) are taken as
\begin{align*}
\phi^{i,+}_j = \delta^{i,+}_j \quad \text{on} \ \partial T^+_i, \qquad j = 1, 2, \cdots, J_i.
\end{align*}
Therefore, we can obtain the local snapshot space on the coarse block $T^+_i$ as
\begin{equation*}
W^{i,+}_{\textrm{snap}} = \textrm{span}\{\phi^{i,+}_1,\phi^{i,+}_2,\cdots,\phi^{i,+}_{J_i}\},
\end{equation*}
we also define the following local space $V^{i,+}_{\textrm{snap}}$ expanded by the divergence-free velocity fields of local snapshot basis functions in $W^{i,+}_{\textrm{snap}}$, which will be used in the spectral decomposition of the next subsection to derive the offline basis functions,
\begin{equation*}
V^{i,+}_{\textrm{snap}} = \textrm{span}\{\psi^{i,+}_1,\psi^{i,+}_2,\cdots,\psi^{i,+}_{J_i}\}.
\end{equation*}
In addition, for the convergence analysis, we also introduce the local snapshot space $W^i_{\textrm{snap}}$ on the coarse element $T_i$ by restricting the space $W^{i,+}_{\textrm{snap}}$ on $T_i$, written as
\begin{equation}\label{eqn_snapshot_T_i}
W^i_{\textrm{snap}} = W^{i,+}_{\textrm{snap}}|_{T_i},
\end{equation}
and the corresponding space $V^i_{\textrm{snap}}$ which is composed of velocity fields of the local snapshot space $W^i_{\textrm{snap}}$, we notice that all functions in $W^i_{\textrm{snap}}$ take zero values outside of $T_i$.\\

\noindent\textbf{Remark 2.2.} We have employed oversampling techniques to obtain the more effective snapshot space. Note that the number of fine-grid layers enlarged by the coarse block $T_i^+$ around $T_i$ is relative smaller than the fine-grid partition for $T_i$ in each coordinate direction, for example, as shown in the right graph of Figure \ref{Fig_mesh_Oversampling}, the partition for $T_i$ is a $10\times 10$ uniform fine grid, $T_i^+$ is obtained by adding two fine-grid layers outside of $T_i$ with the partition of a $14\times 14$ uniform fine grid, so that $W^{i,+}_{\textrm{snap}}$ and $W^i_{\textrm{snap}}$ are able to have the same dimensions, besides,
we can use the randomized oversampling technique \cite{calo2016randomized} to achieve this end, where the basis functions of the local snapshot space $W^{i,+}_{\textrm{snap}}$ are solved with a smaller number of randomized boundary conditions on $\partial T^+_i$, so that we can obtain a smaller dimensional local snapshot space $W^{i,+}_{\textrm{snap}}$, which can also greatly reudce the computational costs in the offline stage. Then, for any function $q\in W^i_{\textrm{snap}}$, there exists a unique $q^+\in W^{i,+}_{\textrm{snap}}$, satisfing $q=q^+|_{T_i}$, i.e., $q^+$ can be uniquely determined by $q$ in reverse. Correspondingly, let $\mathbf{v}\in V^i_{\textrm{snap}}$ be the velocity field with respect to $q$ and $\mathbf{v}^+\in V^{i,+}_{\textrm{snap}}$ be the divergence-free velocity field with respect to $q^+$, respectively, through the relationship (\ref{eqn_velocity_pressure_relationship}), then $\|\mathbf{v}^+\|_{\kappa^{-1},i^+} \le C_{\textrm{err}}\|\mathbf{v}\|_{\kappa^{-1},i^o}$, where the constant $C_{\textrm{err}}$ depends on the permeability $\kappa$ around the fine-grid layers $T^+_i\backslash T_i$ and is independent of $\mathbf{v}$.
\subsubsection{Offline space}
We construct the local offline space for each coarse element by performing a dimension reduction in each local snapshot space. To this end, we conduct the following spectral decomposition to capture the dominant modes of each local snapshot space $W^{i,+}_{\textrm{snap}}$: find a real number $\lambda^i_k\ge 0$ and a vector $\Phi^i_k$, such that
\begin{equation}\label{eqn_spectral_decomposition}
A^{i,+}_{\textrm{off}}\Phi^i_k = \lambda^i_k M^{i,+}_{\textrm{off}}\Phi^i_k,\quad A^{i,+}_{\textrm{off}} = {R^{\scriptscriptstyle{V},i}_{\textrm{off}}}'A^{i,+}R^{\scriptscriptstyle{V},i}_{\textrm{off}},\quad M^{i,+}_{\textrm{off}} = {R^{\scriptscriptstyle{W},i}_{\textrm{off}}}'M^{i,+}R^{\scriptscriptstyle{W},i}_{\textrm{off}},
\end{equation}
where, $R^{\scriptscriptstyle{V},i}_{\textrm{off}}$ and $R^{\scriptscriptstyle{W},i}_{\textrm{off}}$ denote the coefficient matrices of snapshot basis functions in the expansion of fine-grid basis functions, written as
\begin{equation*}
R^{\scriptscriptstyle{V},i}_{\textrm{off}} = [\psi^{i,+}_1,\psi^{i,+}_2,\cdots,\psi^{i,+}_{J_i}] \quad\textrm{and}\quad R^{\scriptscriptstyle{W},i}_{\textrm{off}} = [\phi^{i,+}_1,\phi^{i,+}_2,\cdots,\phi^{i,+}_{J_i}],
\end{equation*}
$A^{i,+}$ and $M^{i,+}$ are fine-grid matrices associated with the following two bilinear forms
\begin{equation*}
A^{i,+} = [a_{rl}] = (\kappa^{-1}\psi_r,\psi_l)_{Q,T^+_i} \quad\textrm{and}\quad M^{i,+} = [s_{rl}] = (\bar{\kappa}\phi_r,\phi_l)_{T^+_i},
\end{equation*}
where $\bar{\kappa}$ is a piecewise constant function on $\mathcal{T}_h$ defined by
\begin{equation}\label{eqn_kappa_bar}
\bar{\kappa}_t = \sum\limits_{e\in\partial t}\bar{\kappa}_e, \quad \forall t\in\mathcal{T}_h.
\end{equation}
We arrange the eigenvalues of the spectral decomposition (\ref{eqn_spectral_decomposition}) in increasing order, $\lambda_1^{i}\le \lambda_2^{i}\le \cdots \le \lambda_{J_i}^{i}$,
and choose the first $l_{i}$ eigenvalues $\lambda_k^{i}$ and the corresponding eigenvectors $\Phi^i_k=(\Phi_{kj})^{J_i}_{j=1}$ to form the local offline space with respect to $T_i$, where $(\Phi^i_{kj})$ is the $j$-th component of the vector $\Phi^i_k$ for $k=1,\cdots,l_i$. We get the following eigenfunctions $\phi^{i,\textrm{off},+}_k$ as
\begin{equation*}
\phi^{i,\textrm{off},+}_k = \sum^{J_i}_{j=1}\Phi^i_{kj}\phi^{i,+}_j, \quad k=1,2,\dots,l_{i},
\end{equation*}
we assume that these eigenfunctions are normalized such that $\int_{T^+_i}\bar{\kappa}(\phi^{i,\textrm{off},+}_k)^2=1$, $k=1,2,\dots,l_{i}$. Then, restricting these eigenfunctions on $T_i$, we obtain the following local offline basis functions
\begin{equation*}
\phi^{i,\textrm{off}}_k = \phi^{i,\textrm{off},+}_k|_{T_i}, \quad k=1,2,\dots,l_i,
\end{equation*}
and the local offline space with respect to $T_i$ is defined as
\begin{equation*}
W^{i}_{\textrm{off}} = \textrm{span}\{\phi^{i,\textrm{off}}_1,\phi^{i,\textrm{off}}_2,\cdots,\phi^{i,\textrm{off}}_{l_{i}}\}.
\end{equation*}
In the same way, for the convenience of convergence analysis, we introduce the local space $V^{i}_{\textrm{off}}$ which is composed of velocity fields of the local offline space $W^i_{\textrm{off}}$, we notice that all functions in $W^i_{\textrm{off}}$ take zero values outside of $T_i$.
And we also define the following space
\begin{equation}
W^{i,+}_{\textrm{off}} = \textrm{span}\{\phi^{i,\textrm{off},+}_1,\phi^{i,\textrm{off},+}_2,\cdots,\phi^{i,\textrm{off},+}_{l_{i}}\}.
\end{equation}

Combining all these local offline spaces $W^{i}_{\textrm{off}}$, $i=1,2,\cdots,N_T$, together, we get the global offline space $W_{\textrm{off}}$ for pressure, and by use of the single-index notation, it can be written as $W_{\textrm{off}} = \textrm{span}\{\psi^{\textrm{off}}_{k} : 1\le k\le M_{\textrm{off}}\}$,
where $M_{\textrm{off}} = \sum^{N_T}_{i=1}l_i$ is the total dimensions of the global offline space. Having gotten the offline space $W_{\textrm{off}}$ for pressure, we define the following multiscale spaces for velocity and pressure, respectively, as
\begin{equation}\label{eqn_multiscale_space}
V_{\textrm{ms}}=V_h \quad\textrm{and} \quad W_{\textrm{ms}}=W_{\textrm{off}},
\end{equation}
we can see that the fine-grid space is employed for the approximation of velocity.
\subsubsection{Correction function}
The correction function is introduced for the consistent handling of the related source term $f$ in (\ref{eqn_model_2}). We solve the following local problem on the coarse block $T^+_i$,
\begin{align}\label{eqn_local_problem_correction}
\begin{split}
\kappa^{-1}\mathbf{\tilde{\psi}}^{i,+} + \nabla\tilde{\phi}^{i,+} &= 0 \qquad \mathrm{in} \ \ T^+_i, \\
\nabla\cdot\mathbf{\tilde{\psi}}^{i,+} & = f \qquad \mathrm{in} \ \ T^+_i.
\end{split}
\end{align}
with homogeneous Neumann boundary conditions, the function $\tilde{\phi}^{i,+}$ is uniquely determined by specifying the condition that $\int_{T^+_i}\tilde{\phi}^{i,+}=0$. Then, restricting $\tilde{\phi}^{i,+}$ on $T_i$, we get the local correction function, represented by $\tilde{\phi}^{i} = \tilde{\phi}^{i,+}|_{T_i}$. And the global correction function is obtained by adding all these local correction functions together, denoted by $\tilde{p}_{\textrm{ms}} = \sum^{N_T}_{i=1}{\tilde{\phi}}^{i}$.
\subsubsection{Multiscale solution}
The multiscale pressure will be taken the form as $p_{\textrm{ms}}=\bar{p}_{\textrm{ms}} + \tilde{p}_{\textrm{ms}}$ with $\bar{p}_{\textrm{ms}}\in W_{\textrm{ms}}$ and $\tilde{p}_{\textrm{ms}}$ denoting the global correction function.
We get the corresponding mixed GMsFEM system: find $(\mathbf{u}_{\textrm{ms}},\bar{p}_{\textrm{ms}})\in V_{\textrm{ms}}\times W_{\textrm{ms}}$, such that
\begin{eqnarray}
(\kappa^{-1}\mathbf{u}_{\textrm{ms}},\mathbf{v})_Q - (\bar{p}_{\textrm{ms}},\nabla\cdot\mathbf{v})
= -(g_D,\mathbf{v}\cdot\mathbf{n})_{\partial\Omega_D} + (\tilde{p}_{\textrm{ms}},\nabla\cdot\mathbf{v}), &&\forall\mathbf{v}\in V_{\textrm{ms}},\label{eqn_mixed_gms_1}\\
-(\nabla\cdot\mathbf{u}_{\textrm{ms}},q) = -(f,q), \hspace{3.46cm}&&\forall q\in W_{\textrm{ms}}.\label{eqn_mixed_gms_2}
\end{eqnarray}
The matrix form of the above multiscale discrete weak formulation can be written as: find $({U}_{\textrm{ms}},\bar{P}_{\textrm{ms}})\in\mathbb{R}^{m_1}\times\mathbb{R}^{M_{\textrm{off}}}$, such that
\begin{equation}\label{eqn_matrix_gms}
\left(
\begin{matrix}
A_h &B_hR_{\textrm{off}}\\
R_{\textrm{off}}^TB^T_h &0
\end{matrix}\right)
\left(
\begin{matrix}
U_{\textrm{ms}}\\
\bar{P}_{\textrm{ms}}
\end{matrix}\right)=
\left(
\begin{matrix}
G_{\textrm{ms}}\\
R_{\textrm{off}}^TF_h
\end{matrix}\right).
\end{equation}
where $R_{\textrm{off}}$ denotes the coefficients matrix of offline basis functions in the expansion of fine-grid basis functions, the vector $G_{\textrm{ms}}$ is associated with the term  $-(g_D,\mathbf{v}\cdot\mathbf{n})_{\partial\Omega_D} + (\tilde{p}_{\textrm{ms}},\nabla\cdot\mathbf{v})$. In ({\ref{eqn_matrix_fine}), we have know that $A_h$ is a diagonal matrix with positive diagonal elements, and can be inverted easily, so we solve the system (\ref{eqn_matrix_gms}) in the following way
\begin{equation}\label{matrixformGMs}
-R_{\textrm{off}}^TB^T_h({A_h})^{-1}B_hR_{\textrm{off}}\bar{P}_{\textrm{ms}} = R_{\textrm{off}}^TF_h-R_{\textrm{off}}^TB^T_h({A_h})^{-1}G_{\textrm{ms}},
\end{equation}
and similarly, we can write the multiscale discrete weak formulation (\ref{eqn_mixed_gms_1})-(\ref{eqn_mixed_gms_2}) of mixed form into the following multiscale discrete weak formulation that only related to the multiscale pressure as: find $p_{\textrm{ms}}\in \tilde{p}_{\textrm{ms}} + W_{\textrm{ms}}$, such that
\begin{equation}\label{eqn_pressure_form_multiscale}
a(p_{\textrm{ms}},q) = (f,q), \quad \forall q\in W_{\textrm{ms}}.
\end{equation}
where the bilinear form $a(\cdot,\cdot)$ in (\ref{eqn_definition_apq}) is symmetric, continuous and coercive. And the matrix form of the above multiscale discrete weak formulation (\ref{eqn_pressure_form_multiscale}) is the same as the linear system (\ref{matrixformGMs}).
\section{Offline adaptive method}
In this section, we introduce an a-posteriori error indicator on each coarse element. Based on this error indicator, we develop an offline adaptive enrichment algorithm to increase the number of offline basis functions iteratively on coarse elements with large residuals for solving the multiscale problem (\ref{eqn_model_1})-(\ref{eqn_model_2}). The error indicator is constructed depending on one weighted $L^2$-norm of the local residual operator on the local snapshot space, combined with the eigenvalue structure of the local spectral decomposition in the offline stage, where the weighted $L^2$-norm is related to the pressure fields of the local snapshot space. For simplicity, we call it pressure-related weighted $L^2$-norm.

For the coarse element $T_i$, $i=1,2,\cdots,N_T$, we define the local residual operator $R_i$ as a linear functional on the local snapshot space $W^{i,+}_{\textrm{snap}}$. For any $q^+\in W^{i,+}_{\textrm{snap}}$, let $q$ be the restriction of $q^+$ on the coarse element $T_i$, satisfying $q=q^+$ in $T_i$, and $q=0$ outside of $T_i$, from the discussion in Remark 2.2, $q^+$ can also be uniquely determined by $q$ in reverse, we define the local residual operator $R_i$ as follows
\begin{equation}\label{eqn_residual_operator}
R_iq^+ = \int_{T_i}fq - a(p_{\textrm{ms}},q) = \int_{T_i}fq-\int_{T_i}\nabla\cdot\mathbf{u}_{\textrm{ms}}q.
\end{equation}
We introduce the following weighted $L^2$-norm, i.e., the so-called pressure-related weighted $L^2$-norm, on the local snapshot space $W^{i,+}_{\textrm{snap}}$, defined as
\begin{equation}\label{eqn_snapshot_norm_1}
\|q^+\|_{\bar{\kappa},i^+} = (\bar{\kappa}q^+,q^+)^{\frac{1}{2}}_{T^+_i},
\end{equation}
where $\bar{\kappa}$ is a piecewise constant function on $\mathcal{T}_h$ defined in (\ref{eqn_kappa_bar}). Correspondingly, we define the following norm of the local residual operator $R_i$ associated with the above norm $\|\cdot\|_{\bar{\kappa},i^+}$ on $W^{i,+}_{\textrm{snap}}$, by
\begin{equation}\label{eqn_residual_norm}
\|R_i\| = \sup\limits_{q^+\in W^{i,+}_{\textrm{snap}}}\frac{|R_iq^+|}{\|q^+\|_{\bar{\kappa},i^+}}.
\end{equation}
The above pressure-related weighted $L^2$-norm of $R_i$ gives estimate on fine-grid residual errors with respect to the local snapshot space $W^{i,+}_{\textrm{snap}}$. In this section, we will take $\|R_i\|^2(\lambda^{i}_{l_i+1})^{-1}$ as our error indicator, where $\lambda^{i}_1, \lambda^{i}_2, \cdots, \lambda^{i}_{J_i}$ are eigenvalues of the local spectral decomposition problem in (\ref{eqn_spectral_decomposition}). 
\begin{lemma}\label{lemma_1}
Let $\mathbf{u}_h$ be the fine-grid solution and $\mathbf{u}_{\textrm{ms}}$ be the multiscale solution on the coarse grid, then
\begin{equation}
\|\mathbf{u}_h - \mathbf{u}_{\textrm{ms}}\|^2_{\kappa^{-1}} \le C_{\textrm{err}}\sum\limits^{N_T}_{i=1}\|R_i\|^2(\lambda^{i}_{l_i+1})^{-1},\label{eqn_lemma_1_w}
\end{equation}
where $C_{\textrm{err}}$ is the constant in Remark 2.2.
\end{lemma}
Based on the above estimate in Lemma \ref{lemma_1} for the proposed error indicator, we can see that the error of the multiscale function can be bounded by the sum of all local residuals. Therefore, we present the offline adaptive enrichment algorithm in the following. We employ $m\ge 1$ to represent the enrichment level, and in the enrichment level $m$, we use $W^m_{\textrm{ms}}$ to denote the corresponding multiscale space for pressure, and use $l^m_i$ to denote the number of offline basis functions on the coarse element $T_i$, $i=1,2,\cdots,N_T$.\\

\noindent{\bf Offline adaptive enrichment algorithm:} Assume that the initial multiscale space $W^1_{\textrm{ms}}$ is given, with $l^1_i$ offline basis functions on the coarse element $T_i$, $i=1,2,\cdots,N_T$. Choose a fixed real numbers $\theta_{\textrm{off}}$, such that $0<\theta_{\textrm{off}}<1$. For $m = 1,2,\cdots$, we perform the following steps,
\begin{itemize}
	\item[]
	\begin{itemize}
    \item[Step 1:] Solve the multiscale problem in the current level. That is, find solutions $(\mathbf{u}^m_{\textrm{ms}},p^m_{\textrm{ms}})\in{V_{\textrm{ms}}\times\big(\tilde{p}_{\textrm{ms}}+{W}^m_{\textrm{ms}}}\big)$, such that
    \begin{eqnarray*}
          (\kappa^{-1}\mathbf{u}^m_{\textrm{ms}},\mathbf{v})_Q - (p^m_{\textrm{ms}},\nabla\cdot\mathbf{v}) = -(g_D,\mathbf{v}\cdot\mathbf{n})_{\partial\Omega_D}, &&\forall\mathbf{v}\in V_{\textrm{ms}},\\
          -(\nabla\cdot\mathbf{u}^m_{\textrm{ms}},q) = -(f,q), \hspace{1.35cm} &&\forall q\in W^m_{\textrm{ms}}.
    \end{eqnarray*}
	
    \item[Step 2:] Calculate the local error indicator on each coarse element. By use of the pressure-related weighted $L^2$-norm of the local residual operator in (\ref{eqn_residual_norm}), on the coarse element $T_i$, $i=1,2,\cdots,N_T$, we compute the local error indicator $\eta_i$ as
    \begin{equation}\label{eqn_eta_offline}
          \eta^2_i = \|R_i\|^2(\lambda^{i}_{l_i+1})^{-1},
    \end{equation}
    where $R_i$ is the local residual operator defined in (\ref{eqn_residual_operator}). After we have computed the local error indicator $\eta_i$ on all coarse elements, we rearrange them in decreasing order, $\eta_1\ge\eta_2\ge\cdots\ge\eta_{\scriptscriptstyle{N_T}}$.
	
    \item[Step 3:] Choose coarse elements where the enrichment with offline basis functions is needed. We choose the smallest integer $N_{\textrm{add}}$ such that the cumulative residuals on the selected coarse elements is $\theta_{\textrm{off}}$ fraction of the sum of all residuals, that is
    \begin{equation}\label{eqn_criterion_offline}
          \theta_{\textrm{off}}\sum\limits^{N_T}_{i=1}\eta^2_i \le \sum\limits^{N_{\textrm{add}}}_{i=1}\eta^2_i.
    \end{equation}
	The number of coarse elements that offline basis functions need to be added is determined by the parameter $\theta_{\textrm{off}}$. We will add offline basis functions on coarse elements where the corresponding error indicator takes values $\eta_1, \eta_2, \cdots, \eta_{\scriptscriptstyle{N_{\textrm{add}}}}$, respectively, to enrich the multiscale space.
	
    \item[Step 4:] Enrich the multiscale space. For the coarse element $T_i$ selected by the above criterion (\ref{eqn_criterion_offline}), we add offline basis functions in the following way. Let $s$ be the smallest positive integer such that $\lambda^{i}_{l^m_i+1+s}$ is large enough compared with $\lambda^{i}_{l^m_i+1}$, then we take $l^{m+1}_i = l^m_i+s$ and add offline basis functions $\phi^{i,\textrm{off}}_{k}$, $k = l^m_i + 1, l^m_i+2, \cdots, l^{m+1}_i$, to enrich the multiscale space.
    \end{itemize}
\end{itemize}
After step 4, we repeat the above procedure form step $1$ again until the global error indicator $\sum^{N_T}_{i=1}\eta^2_i$ is small enough or the total dimension of the multiscale space $W_{\textrm{ms}}$ is large enough.

Next, we define some projection operators that will be used in the convergence analysis of the proposed offline adaptive method.

{\bf Projection operator $P^+_i$}: For the coarse block $T^+_i$, $i = 1, 2, \cdots, N_T$, we define the projection operator $P^+_i: W^{i,+}_{\textrm{snap}} \rightarrow  W^{i,+}_{\textrm{off}}$ from the local snapshot space to the corresponding local offline space in the oversampling region as
\begin{equation*}
\int_{T^+_i}\bar{\kappa}(P^+_iq^+)r = \int_{T^+_i}\bar{\kappa}q^+r, \quad \forall r\in W^{i,+}_{\textrm{off}}.
\end{equation*}
For any $q^+\in W^{i,+}_{\textrm{snap}}$, suppose we can express it as $q^+=\sum^{J_i}_{l=1}c_l\phi^{i,\textrm{off},+}_l$. Then, by the fact that eigenfunctions of the spectral decomposition (\ref{eqn_spectral_decomposition}) are orthogonal, we have $P^+_iq^+=\sum^{l_i}_{l=1}c_l\phi^{i,\textrm{off},+}_l$.

{\bf Projection operator $P_i$}: For the coarse element $T_i$, $i = 1, 2, \cdots, N_T$, making use of the projection operator $P^+_i$, we define the projection operator $P_i: W^{i}_{\textrm{snap}} \rightarrow  W^{i}_{\textrm{off}}$ from the local snapshot space to the corresponding local offline space of $T_i$ in the following way. For any $q\in W^{i}_{\textrm{snap}}$, find $q^+\in W^{i,+}_{\textrm{snap}}$, satisfying $q^+|_{T_i} = q$, by Remark 2.2, $q^+$ is uniquely determined by $q$, then $P_iq$ is defined as the restriction of $P^+_iq^+$ on the coarse element $T_i$, that is 
\begin{equation*}
P_iq = (P^+_iq^+)|_{T_i},
\end{equation*}
we assume that $P_iq$ also takes zero values outside of $T_i$.

{\bf Projection operator $\boldsymbol{P}_i$}: For the coarse element $T_i$, $i = 1, 2, \cdots, N_T$, we define the projection operator $\boldsymbol{P}_i: V^{i}_{\textrm{snap}} \rightarrow V^{i}_{\textrm{off}}$, with $V^{i}_{\textrm{snap}}$ and $V^{i}_{\textrm{off}}$ representing the related velocity fields of the local snapshot space $W^{i}_{\textrm{snap}}$ and the local offline space $W^{i}_{\textrm{off}}$, respectively. For any $\mathbf{v}\in V^{i}_{\textrm{snap}}$, find $q\in W^{i}_{\textrm{snap}}$, such that $\mathbf{v}$ is the velocity field with respect to $q$ on the coarse element $T_i$, then $\boldsymbol{P}_i\mathbf{v}\in V^{i}_{\textrm{off}}$ is defined as the velocity field with respect to $P_iq$ on the coarse element $T_i$, note that $P_iq$ takes zero values outside of $T_i$.

From the above definitions of projection operators $P^+_i$, $P_i$ and $\boldsymbol{P}_i$, we can see that operators $P_i$ and $\boldsymbol{P}_i$ are determined by the operator $P^+_i$. Indeed, they are the same operator with different domains and ranges. In regard to these projection operators, by the spectral decomposition (\ref{eqn_spectral_decomposition}), we have the following estimates
\begin{equation}\label{eqn_projection_norm_inner}
(\kappa^{-1}(\mathbf{v}-\boldsymbol{P}_i{\mathbf{v}}),\mathbf{v}-\boldsymbol{P}_i{\mathbf{v}})_{Q,{T}^0_i}
\le (\kappa^{-1}(\mathbf{v}-\boldsymbol{P}_i{\mathbf{v}}),\mathbf{v}-\boldsymbol{P}_i{\mathbf{v}})_{Q,{T}^+_i}
= \sum\limits^{J_i}_{l=l_i+1}\lambda^i_lc^2_l,
\end{equation}
and
\begin{equation}\label{eqn_projection_norm_boundary}
\sum\limits_{e\in\partial T_i}\bar{\kappa}_e(q-P_iq)^2_{t_e}
\le \sum\limits_{t\in T^+_i}\bar{\kappa}_t(q^+-P^+_iq^+)^2_t
= \frac{1}{h^2}\int_{T^+_i}\bar{\kappa}(q^+-P^+_iq^+)^2
= \frac{1}{h^2}\sum\limits^{J_i}_{l=l_i+1}c^2_l.
\end{equation}
where the subscript $t_e$ denotes the fine-grid element $t\in T_i$ such that $e\in\partial t$. Thereby, making use of (\ref{eqn_projection_norm_inner}), (\ref{eqn_projection_norm_boundary}) and the fact that the eigenvalues $\lambda^i_l$ of the spectral decomposition are increasingly ordered, we derive the following bound for $\|\boldsymbol{P}_i\mathbf{v}-\mathbf{v}\|^2_{\kappa^{-1},i}$, written as
\begin{equation}\label{eqn_projection_bound_2}
\begin{split}
\|\boldsymbol{P}_i\mathbf{v}-\mathbf{v}\|^2_{\kappa^{-1},i}
&= (\kappa^{-1}(\mathbf{v}-\boldsymbol{P}_i{\mathbf{v}}),\mathbf{v}-\boldsymbol{P}_i{\mathbf{v}})_{Q,{T}^o_i} + \sum\limits_{e\in\partial T_i}\bar{\kappa}_e\llbracket q-P_iq\rrbracket^2_e \\
&= (\kappa^{-1}(\mathbf{v}-\boldsymbol{P}_i{\mathbf{v}}),\mathbf{v}-\boldsymbol{P}_i{\mathbf{v}})_{Q,{T}^o_i} + \sum\limits_{e\in\partial T_i}\bar{\kappa}_e(q-P_iq)^2_{t_e} \\
&\le ({\lambda^i_l} + \frac{1}{h^2} ) \sum\limits^{J_i}_{l=l_i+1}c^2_l \le (\Lambda_i + \frac{1}{h^2} ) \sum\limits^{J_i}_{l=1}c^2_l = (\Lambda_i + \frac{1}{h^2} )\|q^+\|^2_{\bar{\kappa},i^+},
\end{split}
\end{equation}
where $\Lambda_i = \max\limits_{1\le l \le J_i}\lambda^i_l = \lambda^i_{J_i}$. And likewise,
\begin{equation*}
\begin{split}
\|\boldsymbol{P}_i\mathbf{v}\|^2_{\kappa^{-1},i}
&= (\kappa^{-1}\boldsymbol{P}_i{\mathbf{v}},\boldsymbol{P}_i{\mathbf{v}})_{Q,{T}^o_i} + \sum\limits_{e\in\partial T_i}\bar{\kappa}_e\llbracket P_iq\rrbracket^2_e \\
&= (\kappa^{-1}\boldsymbol{P}_i{\mathbf{v}},\boldsymbol{P}_i{\mathbf{v}})_{Q,{T}^o_i} + \sum\limits_{e\in\partial T_i}\bar{\kappa}_e(P_iq)^2_{t_e}
\le \sum\limits^{l_i}_{l=1}\lambda^i_lc^2_l + \frac{1}{h^2}\sum\limits^{l_i}_{l=1}c^2_l \\
&\le ({\lambda^i_{l_i}} + \frac{1}{h^2} ) \sum\limits^{l_i}_{l=1}c^2_l \le ({\lambda^i_{l_i}} + \frac{1}{h^2} ) \sum\limits^{J_i}_{l=1}c^2_l \le ({\lambda^i_{l_i}} + \frac{1}{h^2} )\|q^+\|^2_{\bar{\kappa},i^+}.
\end{split}
\end{equation*}
Note that, in the above derivations, $q$ and $P_iq$ take zero values outside of $T_i$. Thus, the projection operator $\boldsymbol{P}_i$ satisfies the following stability property,
\begin{equation}\label{eqn_projection_bound_1}
\|\boldsymbol{P}_i\mathbf{v}\|^2_{\kappa^{-1},i} \le ({\lambda^i_{l_i}} + \frac{1}{h^2})\|q^+\|^2_{\bar{\kappa},i^+}, \quad i=1,2,\cdots,N_T,
\end{equation}
with $q^+\in W^{i,+}_{\textrm{snap}}$, such that $\mathbf{v}$ is the velocity field with respect to $q = q^+|_{T_i}$.

Furthermore, we establish the approximation property for the projection operator $P^+_i$. In fact, utilizing the definition of the operator $P^+_i$, for any $q^+\in W^{i,+}_{\textrm{snap}}$, we have
\begin{equation} \label{eqn_off_projection_approximation}
\begin{split}
\|q^+-P^+_iq^+\|^2_{\bar{\kappa},i^+} &= \sum\limits_{l\ge l_i+1}c^2_l \le (\lambda^i_{l_i+1})^{-1}\sum\limits_{l\ge l_i+1} {\lambda^i_l}c^2_l \\
&\le (\lambda^i_{l_i+1})^{-1}\sum\limits^{J_i}_{l=1} {\lambda^i_l}c^2_l \le (\lambda^i_{l_i+1})^{-1}\|\mathbf{v}^+\|^2_{\kappa^{-1},i^+},
\end{split}
\end{equation}
with $\mathbf{v}^+\in V^{i,+}_{\textrm{snap}}$ representing the divergence-free velocity field with respect to $q^+$, then the following approximation property for $P^+_i$ holds
\begin{equation}\label{eqn_projection_approximation_pressure}
\|q^+-P^+_iq^+\|_{\bar{\kappa},i^+} \le (\lambda^i_{l_i+1})^{-\frac{1}{2}}\|\mathbf{v}^+\|_{\kappa^{-1},i^+}, \quad i=1,2,\cdots,N_T.
\end{equation}
For the theoretical analysis presented below, we also define the global projection operator $\Pi$ as  
\begin{equation}\label{eqn_interpolator}
\Pi: W_{\textrm{snap}}\rightarrow W_{\textrm{off}} \quad\textrm{with}\quad \Pi q = \sum^{N_T}_{i=1}P_iq, \quad \forall q\in W_{\textrm{snap}}.
\end{equation}
where $W_{\textrm{snap}}=\bigcup^{N_T}_{i=1}W^i_{\textrm{snap}}$ is the global snapshot space, the combination of all local snapshot spaces.

Having gotten the stability and approximation properties for the above projection operators, we give the proof of Lemma \ref{lemma_1} in the following, and after that, we will also conduct the convergence analysis for the proposed offline adaptive method.\\

{\bf Proof :} Let $q$ be an arbitrary function in $W_{\textrm{snap}}$ and $\mathbf{v}$ be the velocity filed with respect to $q$ by the relationship (\ref{eqn_velocity_pressure_relationship}). Combining (\ref{eqn_pressure_form}), (\ref{eqn_bilinear_form_equal}) and (\ref{eqn_pressure_form_multiscale}), we have
\begin{equation}\label{eqn_lemma_1}
\begin{split}
(\kappa^{-1}(\mathbf{u}_h-\mathbf{u}_{\textrm{ms}}),\mathbf{v})_Q &= a(p_h-p_{\textrm{ms}},q)= (f,q)-a(p_{\textrm{ms}},q) \\
&= (f,q-\Pi q) + (f,\Pi q) - a(p_{\textrm{ms}},\Pi q) - a(p_{\textrm{ms}},q-\Pi q).
\end{split}
\end{equation}
Since $\Pi q \in W_{\textrm{off}} = W_{\textrm{ms}}$, then from (\ref{eqn_pressure_form_multiscale}), we get
\begin{equation*}
(f,\Pi q) - a(p_{\textrm{ms}},\Pi q) = 0,
\end{equation*}
therefore, (\ref{eqn_lemma_1}) becomes
\begin{equation*}
(\kappa^{-1}(\mathbf{u}_h-\mathbf{u}_{\textrm{ms}}),\mathbf{v})_Q = (f,q-\Pi q) - a(p_{\textrm{ms}},q-\Pi q).
\end{equation*}
We can write the function $q$ in the form of a summation, i.e., $q = \sum^{N_{T}}_{i=1}q^i$ with $q^i\in W^i_{\textrm{snap}}$, $i=1,2,\cdots,N_T$. Denote $\mathbf{v}^i\in V^i_{\textrm{snap}}$ as the velocity field with respect to $q^i$ defined in $T_i$. By the definition of the interpolator $\Pi$ in (\ref{eqn_interpolator}) and the residual operator $R_i$ in (\ref{eqn_residual_operator}), we obtain
\begin{equation*}
\begin{split}
(\kappa^{-1}(\mathbf{u}_h-\mathbf{u}_{\textrm{ms}}),\mathbf{v})_Q &= (f,q-\Pi q) - a(p_{\textrm{ms}},q-\Pi q) \\
&= \sum\limits^{N_T}_{i=1}(f,q^i-P_iq^i) - a(p_{\textrm{ms}},q^i-P_iq^i) \\
&= \sum\limits^{N_T}_{i=1}(f,q^i-P_iq^i) - (\nabla\cdot\mathbf{u}_{\textrm{ms}},q^i-P_iq^i) = \sum\limits^{N_T}_{i=1}R_i(q^{i,+}-P^+_iq^{i,+}),
\end{split}
\end{equation*}
where $q^{i,+}\in W^{i,+}_{\textrm{snap}}$, satisfying $q^{i,+}|_{T_i} = q^i$, $i = 1, 2, \cdots, N_T$, from Remark 2.2, $q^{i,+}$ is uniquely determined by $q^i$.\\
Thus, making use of the approximation property of $P^+_i$ in (\ref{eqn_projection_approximation_pressure}), we have
\begin{equation}\label{eqn_lemma_1_estimate_1}
|(\kappa^{-1}(\mathbf{u}_h-\mathbf{u}_{\textrm{ms}}),\mathbf{v})_Q| \le \sum\limits^{N_T}_{i=1}\|R_i\|\|q^{i,+}-P^+_iq^{i,+}\|_{\bar{\kappa},i^+} \le \sum\limits^{N_T}_{i=1}\|R_i\|(\lambda^{i}_{l_i+1})^{-\frac{1}{2}}\|\mathbf{v}^{i,+}\|_{\kappa^{-1},i^+},
\end{equation}
From the discussion in Remark 2.2, we know that $\mathbf{v}^+_i$ is uniquely determined by $\mathbf{v}_i$ and there exits a constant $C_{\textrm{err}}$ independent of $\mathbf{v}_i$, such that
\begin{equation*}
\|\mathbf{v}^{i,+}\|_{\kappa^{-1},i^+} \le C_{\textrm{err}}\|\mathbf{v}^i\|_{\kappa^{-1},i^0}.
\end{equation*}
Hence, selecting appropriate $q$, such that the corresponding velocity field $\mathbf{v}$ equals $\mathbf{u}_h - \mathbf{u}_{\textrm{ms}}$ in the above equality (\ref{eqn_lemma_1_estimate_1}), and making use of the fact that $\sum^{N_T}_{i=1}\|\mathbf{v}^i\|^2_{\kappa^{-1},i^0} \le \|\mathbf{v}\|_{{\kappa}^{-1}}$, we obtain the equality (\ref{eqn_lemma_1_w}) and then complete the proof of Lemma \ref{lemma_1}. $\square$

Next, we give the convergence analysis for the offline adaptive enrichment algorithm with the proposed error indicator $\|R_i\|^2(\lambda^{i}_{l_i+1})^{-1}$. In the $m$-th enrichment level, for the coarse element $T_i$, we use $R^m_i$ to denote the residual operator $R_i$ with respect to the multiscale solution $\mathbf{u}^m_{\textrm{ms}}$, employ $P^{m,+}_i$, $P^m_i$, $\boldsymbol{P}^m_i$ to denote the projection operator $P^+_i$, $P_i$, $\boldsymbol{P}_i$, respectively, and define
\begin{equation}\label{eqn_S_m_i}
\begin{split}
S^m_i &= (\lambda^i_{l^m_i+1})^{-\frac{1}{2}}\sup\limits_{q^+\in W^{i,+}_{\textrm{snap}}}\frac{|R^m_i(q^+-P^{m,+}_iq^+)|}{\|q^+\|_{\bar{\kappa},i^+}}.
\end{split}
\end{equation}
For any $q\in W^i_{\textrm{snap}}$, we have $P^m_iq\in W^i_{\textrm{off}}\subset W^m_{\textrm{ms}}$, then by (\ref{eqn_pressure_form_multiscale}),
\begin{equation*}
\int_{T_i}fP^m_iq - a(p^m_{\textrm{ms}},P^m_iq) = 0,
\end{equation*}
and accordingly, $S^m_i$ can be written as
\begin{equation}\label{eqn_S_m_i_calculation}
\begin{split}
S^m_i &= (\lambda^i_{l^m_i+1})^{-\frac{1}{2}}\sup\limits_{q\in W^{i,+}_{\textrm{snap}}}\frac{|R^m_i(q^+-P^{m,+}_iq^+)|}{\|q^+\|_{\bar{\kappa},i^+}} \\
&= (\lambda^i_{l^m_i+1})^{-\frac{1}{2}}\sup\limits_{q^+\in W^{i,+}_{\textrm{snap}}}\frac{|\int_{T_i}f(q-P^m_iq)-a(p^m_{\textrm{ms}},q-P^m_iq)|}{\|q^+\|_{\bar{\kappa},i^+}} \\
&= (\lambda^i_{l^m_i+1})^{-\frac{1}{2}}\sup\limits_{q^+\in W^{i,+}_{\textrm{snap}}}\frac{|\int_{T_i}fq - a(p^m_{\textrm{ms}},q)|}{\|q^+\|_{\bar{\kappa},i^+}} = (\lambda^i_{l^m_i+1})^{-\frac{1}{2}}\sup\limits_{q^+\in W^{i,+}_{\textrm{snap}}}\frac{|R^m_iq^+|}{\|q^+\|_{\bar{\kappa},i^+}} = (\lambda^i_{l^m_i+1})^{-\frac{1}{2}}\|R^m_i\|.
\end{split}
\end{equation}
Note that in above derivation $q=q^+|_{T_i}\in W^i_{\textrm{snap}}$.

\begin{lemma}\label{lemma_2}
For any $\alpha>0$, we have
\begin{equation}
(S^{m+1}_i)^2 \le (1+\alpha)\frac{\lambda^i_{l^m_i+1}}{\lambda^i_{l^{m+1}_i+1}}(S^m_i)^2 + (1+\alpha^{-1})D^{m+1}_i\|\mathbf{u}^{m+1}_{\textrm{ms}}-\mathbf{u}^m_{\textrm{ms}}\|^2_{\kappa^{-1},i},
\end{equation}
where the constant $D^{m+1}_i$ depends on the enrichment level, defined by
\begin{equation}\label{eqn_lemma_2_D_i}
D^{m+1}_i = \big(\frac{\Lambda_i}{\lambda^i_{l^{m+1}_i+1}} + \frac{1}{h^2\lambda^i_{l^{m+1}_i+1}}\big)
\end{equation}
with $\Lambda_i = \max\limits_{0\le l \le J_i}\lambda^i_l = \lambda^i_{J_i}$.
\end{lemma}

{\bf Proof :} For any $q^+\in W^{i,+}_{\textrm{snap}}$, let $q=q^+|_{T_i}\in W^i_{\textrm{snap}}$, then by the definition of the residual operator $R^m_i$, $m\ge 1$, we can deduce that
\begin{equation}\label{eqn_direct_calculation}
\begin{split}
R^{m+1}_i(q^+-P^{m+1,+}_iq+) &= \int_{T_i}f(q-P^{m+1}_iq) - a(p^{m+1}_{\textrm{ms}},q-P^{m+1}_iq) \\
&= \int_{T_i}fq-a(p^{m+1}_{\textrm{ms}},q) = \int_{T_i}fq - a(p^{m}_{\textrm{ms}},q) + a(p^{m}_{\textrm{ms}}-p^{m+1}_{\textrm{ms}},q) \\
&= \int_{T_i}f(q-P^m_iq) - a(p^{m}_{\textrm{ms}},q-P^m_iq) + a(p^{m}_{\textrm{ms}}-p^{m+1}_{\textrm{ms}},q) \\
&= R^m_i(q^+-P^{m,+}_iq+) + a(p^{m}_{\textrm{ms}}-p^{m+1}_{\textrm{ms}},q),
\end{split}
\end{equation}
According to the definition of $S^m_i$ in (\ref{eqn_S_m_i}),
we multiply (\ref{eqn_direct_calculation}) by $(\lambda^i_{l^{m+1}_i+1})^{-\frac{1}{2}}\|q^+\|^{-1}_{\bar{\kappa},i}$ and take superme with respect to $q^+$, by doing this, we arrive at
\begin{equation*}
S^{m+1}_i\le\Big(\frac{\lambda^i_{l^m_i+1}}{\lambda^i_{l^{m+1}_i+1}}\Big)^{\frac{1}{2}}S^m_i + I,
\end{equation*}
where
\begin{equation*}
I = (\lambda^i_{l^{m+1}_i+1})^{-\frac{1}{2}}\sup\limits_{q^+\in W^{i,+}_{\textrm{snap}}}\frac{a(p^m_{\textrm{ms}}-p^{m+1}_{\textrm{ms}},q)}{\|q^+\|_{\bar{\kappa},i}},
\end{equation*}
note that $q=q^+|_{T_i}$, belonging to $W^{i}_{\textrm{snap}}$. For the estimation of the term $I$, we use the fact that $P^m_iq\in W^m_{\textrm{ms}}\subset W^{m+1}_{\textrm{ms}}$ to obtain
\begin{equation*}
a(p^m_{\textrm{ms}},P^m_iq) = (f,P_iq) = a(p^{m+1}_{\textrm{ms}},P^m_iq),
\end{equation*}
which implies
\begin{equation*}
\begin{split}
a(p^m_{\textrm{ms}}-p^{m+1}_{\textrm{ms}},q) &= a(p^m_{\textrm{ms}}-p^{m+1}_{\textrm{ms}},q-P^m_iq) \\
&= \big(\kappa^{-1}(\mathbf{u}^m_{\textrm{ms}}-\mathbf{u}^{m+1}_{\textrm{ms}}),\mathbf{v}-\boldsymbol{P}^m_i\mathbf{v}\big)_{Q,T_i} \\
&\le \|\mathbf{u}^m_{\textrm{ms}}-\mathbf{u}^{m+1}_{\textrm{ms}}\|_{\kappa^{-1},i}\|\mathbf{v}-\boldsymbol{P}^m_i\mathbf{v}\|_{\kappa^{-1},i},
\end{split}
\end{equation*}
where $\mathbf{v}\in V^{i}_{\textrm{snap}}$ is the velocity field respect to $q\in W^{i}_{\textrm{snap}}$, and employ the stability property in (\ref{eqn_projection_bound_2}),
\begin{equation*}
\|\mathbf{v}-\boldsymbol{P}_i\mathbf{v}\|_{\kappa^{-1},i} \le (\Lambda_i + \frac{1}{h^2}) \|q^+ \|^2_{\bar{\kappa},i^+},
\end{equation*}
to get that
\begin{equation*}
I\le(\lambda^i_{l^{m+1}_i+1})^{-\frac{1}{2}}(\Lambda_i + \frac{1}{h^2})\|\mathbf{u}^m_{\textrm{ms}}-\mathbf{u}^{m+1}_{\textrm{ms}}\|_{\kappa^{-1},i}.
\end{equation*}
Thus, by the definition of $D^{m+1}_i$ in (\ref{eqn_lemma_2_D_i}), Lemma \ref{lemma_2} is proved. $\square$

\begin{theorem}\label{theorem_3_1}
There are a sequence $\{L_m\}^M_{m=1}$ and positive constants $\epsilon$, $\rho$, $\tau$ independent of the enrichment level $m$ such that the following contracting property holds
\end{theorem}
\begin{equation}
\|\mathbf{u}_h-\mathbf{u}^{m+1}_{\textrm{ms}}\|^2_{\kappa^{-1}}+\frac{\tau}{1+\tau L_{m+1}}\sum\limits^{N_T}_{i=1}(S^{m+1}_i)^2 \le
\varepsilon_{\textrm{off}}\Big(\|\mathbf{u}_h-\mathbf{u}^m_{\textrm{ms}}\|^2_{\kappa^{-1}}+\frac{\tau}{1+\tau L_m}\sum\limits^{N_T}_{i=1}(S^m_i)^2\Big),
\end{equation}	
{\it where} $\varepsilon_{\textrm{off}}$ {\it is the convergence rate satisfying} $ 0<\varepsilon_{\textrm{off}}<1$, {\it denoted by}
\begin{equation}
\varepsilon_{\textrm{off}} = 1-\epsilon\theta_{\textrm{off}}\frac{\tau}{2C_{err}(1+\tau L_1)}.
\end{equation}

{\bf Proof :} Let $0<\theta_{\textrm{off}}<1$, according to the criterion (\ref{eqn_criterion_offline}), we choose an index set $I_{\textrm{add}}$ satisfying
\begin{equation}\label{eqn_criterion_offline_I}
\sum^{N_T}_{i=1}\eta^2_i\le\frac{1}{\theta_{\textrm{off}}}\sum_{i\in I_{\textrm{add}}}\eta^2_i,
\end{equation}
and we will increase the number of offline basis functions on the coarse element $T_i$, $i\in I_{\textrm{add}}$, to enrich the multiscale space. By the definition of $\eta_i$ in (\ref{eqn_eta_offline}), results in Lemma \ref{lemma_1} and the calculation in (\ref{eqn_S_m_i_calculation}), 
we derive that
\begin{equation}\label{eqn_criterion_offline_2}
\|\mathbf{u}_h - \mathbf{u}^m_{\textrm{ms}}\|_{\kappa^{-1}} \le C_{err}\sum\limits^{N_T}_{i=1}\eta^2_i \le \frac{C_{err}}{\theta_{\textrm{off}}}\sum\limits_{i\in I_{\textrm{add}}}\eta^2_i = \frac{C_{err}}{\theta_{\textrm{off}}}\sum\limits_{i\in I_{\textrm{add}}} (S^m_i)^2.
\end{equation}
On the other hand,
\begin{equation*}
\sum\limits^{N_T}_{i=1}(S^{m+1}_i)^2 = \sum\limits_{i\in I_{\textrm{add}}}(S^{m+1}_i)^2 + \sum\limits_{i\notin I_{\textrm{add}}}(S^{m+1}_i)^2.
\end{equation*}
Thanks to Lemma \ref{lemma_2}, if $i\in I_{\textrm{add}}$, we have
\begin{equation}\label{eqn_theorem_off_situation_1}
(S^{m+1}_i)^2 \le (1+\alpha)\frac{\lambda^i_{l^m_i+1}}{\lambda^i_{l^{m+1}_i+1}}(S^m_i)^2 + (1+\alpha^{-1})D^{m+1}_i\|\mathbf{u}^{m+1}_{\textrm{ms}}-\mathbf{u}^m_{\textrm{ms}}\|^2_{\kappa^{-1},i},
\end{equation}
otherwise, if $i\notin I_{\textrm{add}}$, then there is no new offline basis functions being added, i.e., $\lambda^i_{l^m_i+1} = \lambda^i_{l^{m+1}_i+1}$, which implies
\begin{equation}\label{eqn_theorem_off_situation_2}
(S^{m+1}_i)^2 \le (1+\alpha)(S^m_i)^2 + (1+\alpha^{-1})D^{m+1}_i\|\mathbf{u}^{m+1}_{\textrm{ms}}-\mathbf{u}^m_{\textrm{ms}}\|^2_{\kappa^{-1},i}.
\end{equation}
Adding the above two situations (\ref{eqn_theorem_off_situation_1}) and (\ref{eqn_theorem_off_situation_2}) together, we obtain
\begin{equation}\label{eqn_I_summation}
\begin{split}
\sum\limits^{N_T}_{i=1}(S^{m+1}_i)^2 &\le \sum\limits_{i\in I_{\textrm{add}}}\Big((1+\alpha)\frac{\lambda^i_{l^m_i+1}}{\lambda^i_{l^{m+1}_i+1}}(S^m_i)^2 + (1+\alpha^{-1})D^{m+1}_i\|\mathbf{u}^{m+1}_{\textrm{ms}}-\mathbf{u}^m_{\textrm{ms}}\|^2_{\kappa^{-1},i}\Big)\\
&+\sum\limits_{i\notin I_{\textrm{add}}}\Big((1+\alpha)(S^m_i)^2 + (1+\alpha^{-1})D^{m+1}_i\|\mathbf{u}^{m+1}_{\textrm{ms}}-\mathbf{u}^m_{\textrm{ms}}\|^2_{\kappa^{-1},i}\Big).
\end{split}
\end{equation}
Suppose there exists a positive constant $\delta$ independent of the enrichment level $m$, such that each enrichment of the multiscale space satisfies
\begin{equation*}
\max\limits_{i\in I_{\textrm{add}}}\frac{\lambda^i_{l^m_i+1}}{\lambda^i_{l^{m+1}_i+1}}\le\delta<1, \quad m = 1, 2, \cdots, M,
\end{equation*}
then we have the following estimate for the above inequality (\ref{eqn_I_summation}) as
\begin{equation}
\begin{split}
\sum\limits^{N_T}_{i=1}(S^{m+1}_i)^2 &\le \sum\limits_{i\in I_{\textrm{add}}}(1+\alpha)\delta(S^m_i)^2\\
&+ \sum\limits_{i\notin I_{\textrm{add}}}(1+\alpha)(S^m_i)^2  + \sum\limits^{N_T}_{i=1}(1+\alpha^{-1})D^{m+1}_i\|\mathbf{u}^{m+1}_{\textrm{ms}}-\mathbf{u}^m_{\textrm{ms}}\|^2_{\kappa^{-1},i}.
\end{split}
\end{equation}
Since $\delta = 1-(1-\delta)$, it can also be written as
\begin{equation}
\sum\limits^{N_T}_{i=1}(S^{m+1}_i)^2 \le \sum\limits^{N_T}_{i=1}(1+\alpha)(S^m_i)^2 - (1+\alpha)(1-\delta)\sum\limits_{i\in I_{\textrm{add}}}(S^m_i)^2 + L_{m+1}\|\mathbf{u}^{m+1}_{\textrm{ms}}-\mathbf{u}^m_{\textrm{ms}}\|^2_{\kappa^{-1}},
\end{equation}
with
\begin{equation}\label{eqn_L_m}
L_{m+1} = N_E(1+\alpha^{-1})\max\limits_{1\le i\le N_T}D^{m+1}_i,
\end{equation}
where $N_E$ is the maximum number of coarse edges of the coarse element in $\mathcal{T}_H$.
 Owing to (\ref{eqn_eta_offline}), (\ref{eqn_S_m_i_calculation}) and (\ref{eqn_criterion_offline_I}), we have
\begin{equation}\label{eqn_S_rho}
\sum\limits^{N_T}_{i=1}(S^{m+1}_i)^2 \le \sum\limits^{N_T}_{i=1}(1+\alpha)(S^m_i)^2 - (1+\alpha)(1-\delta)\theta_{\textrm{off}}\sum\limits^{N_T}_{i=1}(S^m_i)^2 + L_{m+1}\|\mathbf{u}^{m+1}_{\textrm{ms}}-\mathbf{u}^m_{\textrm{ms}}\|^2_{\kappa^{-1}},
\end{equation}
Let $\rho = (1+\alpha)(1-(1-\delta)\theta_{\textrm{off}})$, then the above inequality (\ref{eqn_S_rho}) turns into
\begin{equation}\label{eqn_tau}
\sum\limits^{N_T}_{i=1}(S^{m+1}_i)^2 \le \rho\sum\limits^{N_T}_{i=1}(S^m_i)^2 + L_{m+1}\|\mathbf{u}^{m+1}_{\textrm{ms}}-\mathbf{u}^m_{\textrm{ms}}\|^2_{\kappa^{-1}},
\end{equation}
note that we have chosen suitable constants $\alpha$ and $\epsilon$ such that $1-\rho>\epsilon>0$. By the definition of $L_m$ in (\ref{eqn_L_m}), we know that $\{L_m\}^M_{m=1}$ is a sequence in decreasing order. Let $\tau$ be a constant satisfying
\begin{equation*}
0<\tau\le\frac{1-\rho-\epsilon}{\rho L_1},
\end{equation*}
then, we get
\begin{equation*}
\frac{1+\tau L_m}{1+\tau L_{m+1}} \le 1+\tau L_m \le 1+\tau L_1,
\end{equation*}
and thereby,
\begin{equation}\label{eqn_beta_estimation}
1-\rho\frac{1+\tau L_m}{1+\tau L_{m+1}}\ge1-\rho(1+\tau L_1)\ge\epsilon>0.
\end{equation}
Since $p^m_{\textrm{ms}}-p^{m+1}_{\textrm{ms}}\in W^{m+1}_{\textrm{ms}}$, by a direct calculation, we have
\begin{equation*}
\begin{split}
a(p^m_{\textrm{ms}}-p_h,p^m_{\textrm{ms}}-p_h) &= a(p^m_{\textrm{ms}}-p^{m+1}_{\textrm{ms}}+p^{m+1}_{\textrm{ms}}-p_h,p^m_{\textrm{ms}}-p^{m+1}_{\textrm{ms}}+p^{m+1}_{\textrm{ms}}-p_h) \\
&= a(p^m_{\textrm{ms}}-p^{m+1}_{\textrm{ms}},p^m_{\textrm{ms}}-p^{m+1}_{\textrm{ms}}) + 2a(p^m_{\textrm{ms}}-p^{m+1}_{\textrm{ms}},p^{m+1}_{\textrm{ms}}-p_h) + a(p^{m+1}_{\textrm{ms}}-p_h,p^{m+1}_{\textrm{ms}}-p_h) \\
&= a(p^m_{\textrm{ms}}-p^{m+1}_{\textrm{ms}},p^m_{\textrm{ms}}-p^{m+1}_{\textrm{ms}}) + a(p^{m+1}_{\textrm{ms}}-p_h,p^{m+1}_{\textrm{ms}}-p_h),
\end{split}
\end{equation*}
hence,
\begin{equation}\label{eqn_orthogonality}
\|\mathbf{u}^m_{\textrm{ms}}-\mathbf{u}_h\|^2_{\kappa^{-1}} = \|\mathbf{u}^{m+1}_{\textrm{ms}}-\mathbf{u}^m_{\textrm{ms}}\|^2_{\kappa^{-1}} + \|\mathbf{u}^{m+1}_{\textrm{ms}}-\mathbf{u}_h\|^2_{\kappa^{-1}}.
\end{equation}
Multiplying (\ref{eqn_tau}) by $\tau$, and then added by the inequality $\|\mathbf{u}^{m+1}_{\textrm{ms}}-\mathbf{u}_h\|^2_{\kappa^{-1}}\le\|\mathbf{u}^m_{\textrm{ms}}-\mathbf{u}_h\|^2_{\kappa^{-1}}$, we have
\begin{equation*}
\|\mathbf{u}^{m+1}_{\textrm{ms}}-\mathbf{u}_h\|^2_{\kappa^{-1}} + \tau\sum\limits^{N_T}_{i=1}(S^{m+1}_i)^2 \le \|\mathbf{u}^m_{\textrm{ms}}-\mathbf{u}_h\|^2_{\kappa^{-1}} + \tau\rho\sum\limits^{N_T}_{i=1}(S^m_i)^2 + \tau L_{m+1}\|\mathbf{u}^{m+1}_{\textrm{ms}}-\mathbf{u}^m_{\textrm{ms}}\|^2_{\kappa^{-1}},
\end{equation*}
making use of (\ref{eqn_orthogonality}), the above inequality becomes
\begin{equation*}
\begin{split}
&\|\mathbf{u}^{m+1}_{\textrm{ms}}-\mathbf{u}_h\|^2_{\kappa^{-1}} + \tau\sum\limits^{N_T}_{i=1}(S^{m+1}_i)^2 \\
&\le\|\mathbf{u}^m_{\textrm{ms}}-\mathbf{u}_h\|^2_{\kappa^{-1}} + \tau\rho\sum\limits^{N_T}_{i=1}(S^m_i)^2 + \tau L_{m+1}(\|\mathbf{u}^{m}_{\textrm{ms}}-\mathbf{u}_h\|^2_{\kappa^{-1}}-\|\mathbf{u}^{m+1}_{\textrm{ms}}-\mathbf{u}_h\|^2_{\kappa^{-1}}),
\end{split}
\end{equation*}
which indicates that
\begin{equation}\label{eqn_theorem_1_before_beta}
\|\mathbf{u}^{m+1}_{\textrm{ms}}-\mathbf{u}_h\|^2_{\kappa^{-1}} + \frac{\tau}{1+\tau L_{m+1}}\sum\limits^{N_T}_{i=1}(S^{m+1}_i)^2 \le
\|\mathbf{u}^{m}_{\textrm{ms}}-\mathbf{u}_h\|^2_{\kappa^{-1}} + \frac{\tau\rho}{1+\tau L_{m+1}}\sum\limits^{N_T}_{i=1}(S^{m}_i)^2.
\end{equation}
Selecting
\begin{equation*}
\beta = {\theta_{\textrm{off}}(\theta_{\textrm{off}} + C_{err}\tau^{-1}(1+\tau L_m))}^{-1}{(1-\rho\frac{1+\tau L_m}{1+\tau L_{m+1}})},
\end{equation*}
due to (\ref{eqn_beta_estimation}), we have $0<\beta<1$. And accordingly, by use of (\ref{eqn_criterion_offline_2}), the inequalities (\ref{eqn_theorem_1_before_beta}) becomes
\begin{equation*}
\begin{split}
\|\mathbf{u}^{m+1}_{\textrm{ms}}-\mathbf{u}_h\|^2_{\kappa^{-1}} &+ \frac{\tau}{1+\tau L_{m+1}}\sum\limits^{N_T}_{i=1}(S^{m+1}_i)^2 \\
&\le (1-\beta)\|\mathbf{u}^{m}_{\textrm{ms}}-\mathbf{u}_h\|^2_{\kappa^{-1}} + \big(\frac{\beta C_{err}}{\theta_{\textrm{off}}} + \frac{\tau\rho}{1+\tau L_{m+1}}\big)\sum\limits^{N_T}_{i=1}(S^{m}_i)^2 \\
&\le (1-\beta)\|\mathbf{u}^{m}_{\textrm{ms}}-\mathbf{u}_h\|^2_{\kappa^{-1}} + \frac{\tau(1-\beta)}{1+\tau L_m}\sum\limits^{N_T}_{i=1}(S^{m}_i)^2.
\end{split}
\end{equation*}
Obviously, $C_{\textrm{err}}$, $L_1\ge 1$, then we obtain
\begin{equation*}
\theta_{\textrm{off}}<1<\frac{C_{\textrm{err}}(1+\tau L_1)}{\tau},
\end{equation*}
and
\begin{equation*}
\beta\ge\epsilon\theta_{\textrm{off}}\big(C_{\textrm{err}}\tau^{-1}(1+\tau L_1) + C_{\textrm{err}}\tau^{-1}(1+\tau L_m)\big)^{-1} \ge \epsilon\theta_{\textrm{off}}\big(2C_{\textrm{err}}\tau^{-1}(1+\tau L_1)\big)^{-1},
\end{equation*}
which gives the required convergence rate as
\begin{equation*}
\varepsilon_{\textrm{off}} = 1-\beta = 1-\epsilon\theta_{\textrm{off}}\big(2C_{\textrm{err}}\tau^{-1}(1+\tau L_1)\big)^{-1}.
\end{equation*}
The proof of Theorem \ref{theorem_3_1} is now completed. $\square$
\section{Online adaptive method}
In this section, we give the online adaptive enrichment algorithm which requires the construction of online basis functions in selected regions based on residual errors and some optimally estimates. Different from offline basis functions that are precomputed in the offline stage before the enrichment algorithm, the online basis functions need to be calculated in the online stage, i.e., the actual simulation. Since the online basis function contains important global information such as distant effects that the offline basis function cannot capture, we can generally get a much faster convergence rate than the offline adaptive enrichment algorithm.

In the online adaptive method, the oversampling techniques are not necessary. So we derive the local snapshot space $W^i_{\textrm{snap}}$ by directly solving snapshot basis functions $(\mathbf{\psi}_j^i, \phi_j^i)$, $j=1,2,\cdots,J_i$, on the coarse element $T_i$, through the local problem (\ref{eqn_local_problem}), and then perform the spectral decomposition (\ref{eqn_spectral_decomposition}) in the local snapshot space $W^i_{\textrm{snap}}$ to get the the local offline space $W^i_{\textrm{off}}$, i.e., find a real number $\lambda^i_k\ge 0$ and a vector $\Phi^i_k$, such that
\begin{equation}\label{eqn_spectral_decomposition_no_oversampling}
A^i_{\textrm{off}}\Phi^i_k = \lambda^i_k M^i_{\textrm{off}}\Phi^i_k,\quad A^i_{\textrm{off}} = {R^{\scriptscriptstyle{V},i}_{\textrm{off}}}'A^iR^{\scriptscriptstyle{V},i}_{\textrm{off}},\quad M^i_{\textrm{off}} = {R^{\scriptscriptstyle{W},i}_{\textrm{off}}}'M^iR^{\scriptscriptstyle{W},i}_{\textrm{off}},
\end{equation}
where, $R^{\scriptscriptstyle{V},i}_{\textrm{off}}$ and $R^{\scriptscriptstyle{W},i}_{\textrm{off}}$ denote the coefficient matrices of snapshot basis functions in the expansion of fine-grid basis functions, written as
\begin{equation*}
R^{\scriptscriptstyle{V},i}_{\textrm{off}} = [\psi^i_1,\psi^i_2,\cdots,\psi^i_{J_i}] \quad\textrm{and}\quad R^{\scriptscriptstyle{W},i}_{\textrm{off}} = [\phi^i_1,\phi^i_2,\cdots,\phi^i_{J_i}],
\end{equation*}
$A^i$ and $M^i$ are fine-grid matrices associated with the following two bilinear forms
\begin{equation*}
A^i = [a_{rl}] = (\kappa^{-1}\psi_r,\psi_l)_{Q,T^0_i} \quad\textrm{and}\quad M^i = [s_{rl}] = (\bar{\kappa}\phi_r,\phi_l)_{T_i}.
\end{equation*}
For any $q\in W^i_{\textrm{snap}}$, let $\mathbf{v}\in V^i_{\textrm{snap}}$ be the velocity field with respect to $q$, note that $q$ takes zero values outside of $T_i$, $\mathbf{v}$ and $q$ are in one-to-one correspondence by the relationship (\ref{eqn_velocity_pressure_relationship}), so we can derive another norm on $W^i_{\textrm{snap}}$ by use of the velocity fields $V^i_{\textrm{snap}}$ of $W^i_{\textrm{snap}}$, written as
\begin{equation}\label{eqn_snapshot_norm_2_V}
\|q\|^2_{\scriptscriptstyle{V_i}} = \|\mathbf{v}\|^2_{\kappa^{-1},i} = (\kappa^{-1}\mathbf{v},\mathbf{v})_{Q,{T}^0_i} + \sum\limits_{e\in\partial T_i}\bar{\kappa}_e\llbracket q \rrbracket^2_e,
\end{equation}
which is referred to as the velocity-related weighted $L^2$-norm on $W^i_{\textrm{snap}}$, where $\|\cdot\|_{\kappa^{-1},i}$ is the local norm for velocity defined in (\ref{eqn_three_local_norms}). Correspondingly, we define the following norm of the local residual operator $R_i$ associated with the above norm (\ref{eqn_snapshot_norm_2_V}) on $W^i_{\textrm{snap}}$, by
\begin{equation}\label{eqn_residual_norm_V}
\|R_i\|_{\scriptscriptstyle{V}} = \sup\limits_{q^+\in W^{i}_{\textrm{snap}}}\frac{|R_iq|}{\|q\|_{\scriptscriptstyle{V_i}}}.
\end{equation}
where the local residual operator $R_i$ is defined as
\begin{equation}\label{eqn_residual_operator_no_oversampling}
R_iq = \int_{T_i}fq - a(p_{\textrm{ms}},q) = \int_{T_i}fq-\int_{T_i}\nabla\cdot\mathbf{u}_{\textrm{ms}}q, \quad \forall q\in W^i_{\textrm{snap}}.
\end{equation}

We employ similar notations as the offline adaptive enrichment algorithm in the previous section. We use the index $m\ge 1$ to represent the enrichment level of the online adaptive method, and in the $m$-th enrichment level, we use $W^m_{\textrm{ms}}$ and $(\mathbf{u}^m_{\textrm{ms}},p^m_{\textrm{ms}})$, respectively, to denote the corresponding multiscale space for pressure and multiscale solutions. As the offline adaptive method, the initial multiscale space $W^1_{\textrm{ms}}$ is composed of the first $l_i$ offline basis functions on the coarse element $T_i$, $i=1,2,\cdots,N_T$. The online adaptive enrichment algorithm is illustrated as follows.\\

\noindent{\bf Online adaptive enrichment algorithm:} Let $m=1$. We begin with the selection of a number of offline basis functions $l_i$ for the coarse element $T_i$, $i=1,2,\cdots,N_T$, respectively, to form the initial multiscale space $W^1_{\textrm{ms}}$. Then, we choose a fixed real number $\theta_{\textrm{on}}$, such that $0<\theta_{\textrm{on}}<1$, and do the following steps.
\begin{itemize}
	\item[]
	\begin{itemize}
		\item[Step 1:] Find multiscale solutions in the current level $m$. That is, seek  $(\mathbf{u}^m_{\textrm{ms}},p^m_{\textrm{ms}})\in{V_{\textrm{ms}}\times\big(\tilde{p}_{\textrm{ms}}+{W}^m_{\textrm{ms}}}\big)$, such that
		\begin{eqnarray*}
		(\kappa^{-1}\mathbf{u}^m_{\textrm{ms}},\mathbf{v})_Q - (p^m_{\textrm{ms}},\nabla\cdot\mathbf{v}) = -(g_D,\mathbf{v}\cdot\mathbf{n})_{\partial\Omega_D}, &&\forall\mathbf{v}\in V_{\textrm{ms}},\\
		-(\nabla\cdot\mathbf{u}^m_{\textrm{ms}},q) = -(f,q), \hspace{1.34cm}&&\forall q\in W^m_{\textrm{ms}}.
		\end{eqnarray*}
		
		\item[Step 2:] Calculate the residual error estimator. For coarse element $T_i$, $i=1,2,\cdots,N_T$, we compute the local error estimator as
		\begin{equation}\label{eqn_online_adaptive_estimator}
		\eta^2_i = \|R_i\|^2_{\scriptscriptstyle{V}},
		\end{equation}
		where the local residual $R_i$ and the corresponding norm $\|\cdot\|_{\scriptscriptstyle{V}}$ are defined in (\ref{eqn_residual_operator}) and (\ref{eqn_residual_norm}), respectively. After we have computed the error estimator on all coarse elements, we rearrange them in decreasing order as $\eta_1\ge\eta_2\ge\cdots\ge\eta_{\scriptscriptstyle{N_T}}$.
		
		\item[Step 3:] Select coarse elements where online multiscale basis functions need to be added. We choose the smallest interger $N_{\textrm{add}}$ such that
		\begin{equation}\label{eqn_criterion_online}
		\theta_{\textrm{on}}\sum\limits^{N_T}_{i=1}\eta^2_i \le \sum\limits^{N_{\textrm{add}}}_{i=1}\eta^2_i,
		\end{equation}
		Obviously, the number of coarse elements that the corresponding online basis functions need to be added is determined by the parameter $\theta_{\textrm{on}}$.
		
		\item[Step 4:] Construct the online basis function and enrich the multiscale space. We construct and add online basis functions on coarse elements where the corresponding error estimator takes value $\eta_1, \eta_2, \cdots, \eta_{\textrm{add}}$, respectively, to enrich the multiscale space. The online basis functions are computed in the following way. Suppose the $m$-th level multiscale space $W^m_{\textrm{ms}}$ and the corresponding multiscale solution $(\mathbf{u}^m_{\textrm{ms}},p^m_{\textrm{ms}})$ are already known, and we need to construct online basis function $\phi$ on the coarse element $T_i$ to enrich the multiscale space, so that $W^{m+1}_{\textrm{ms}} = W^m_{\textrm{ms}} + \textrm{span}\{\phi\}$. Let $T^+_i$ be the coarse block inclusive of $T_i$, $T_i\subset T^+_i$, defined by adding one fine-grid layers around $T_i$. Then, we slove the following local problem: find $(\psi^+,\phi^+)\in V_h(T^+_i)\times W_h(T^+_i)$, such that
		\begin{eqnarray}
		(\kappa^{-1}\psi^+,\mathbf{v})_Q - (\phi^+,\nabla\cdot\mathbf{v})
		= 0, \hspace{2.15cm} &&\forall \mathbf{v}\in V_h(T^+_i),\label{eqn_online_basis_1}\\
		(\nabla\cdot \psi^+,q) = (f-\nabla\cdot\mathbf{u}^m_{\textrm{ms}},q), &&\forall q\in W_h(T^+_i),\label{eqn_online_basis_2}
		\end{eqnarray}
		with homogeneous Neumann boundary conditions $\psi^+\cdot n = 0$ on $\partial T^+_i$, and $\phi^+$ is uniquely solved under the condition that $\phi^+=0$ on the boundary elements of $T^+_i$. Restricting $\phi^+$ on $T_i$, we get the online basis function on the coarse element $T_i$, written as $\phi = \phi^+|_{T_i}$. Since $\phi^+=0$ outside of $T_i$, we have that $\phi^+=\phi$. The discrete system (\ref{eqn_online_basis_1})-(\ref{eqn_online_basis_2}) is equivalent to
		\begin{equation}\label{eqn_online_basis}
		a(\phi,q) = \int_{T_i}fq - a(p^m_{\textrm{ms}},q), \quad\forall q\in W_h(T_i).
		\end{equation}
	\end{itemize}
\end{itemize}
After step 4, we repeat the above procedure from step $1$ until the global error estimator $\sum^{N_T}_{i=1}\eta^2_i$ is small enough.\\

Since we have not used the oversampling techniques, for the coarse block $T_i$, $i = 1,2,\cdots,N_T$, the projection operator $P_i: W^i_{\textrm{snap}} \rightarrow  W^i_{\textrm{off}}$ from the local snapshot space to the corresponding local offline space is redefined as
\begin{equation*}
\int_{T_i}\bar{\kappa}(P_iq)r = \int_{T_i}\bar{\kappa}qr, \quad \forall q\in W^i_{\textrm{snap}},\ \forall r\in W^i_{\textrm{off}}.
\end{equation*}
For any $q\in W^i_{\textrm{snap}}$, suppose we can express it as $q=\sum^{J_i}_{l=1}c_l\phi^{i,\textrm{off}}_l$. Then, by the fact that eigenfunctions of the spectral decomposition (\ref{eqn_spectral_decomposition_no_oversampling}) are orthogonal, we have $P_iq=\sum^{l_i}_{l=1}c_l\phi^{i,\textrm{off}}_l$. The projection operator $\boldsymbol{P}_i: V^{i}_{\textrm{snap}} \rightarrow V^{i}_{\textrm{off}}$, $i = 1,2,\cdots,N_T$, are defiend in the same way as the previous section.

Similar to the previous section, using the spectral decomposition (\ref{eqn_spectral_decomposition_no_oversampling}) and by a direct calculation, we have
\begin{equation}\label{eqn_off_projection_approximation_no_oversampling}
\|q-P_iq\|^2_{\bar{\kappa},i} \le (\lambda^i_{l_i+1})^{-1}\|\mathbf{v}\|^2_{\kappa^{-1},i^0},
\end{equation}
where $\|q\|_{\bar{\kappa},i} = (\bar{\kappa}q,q)^{\frac{1}{2}}_{T_i}$, and $\mathbf{v}\in V^{i}_{\textrm{snap}}$ is the velocity with respect to $q$. From the definition of the norm $\|\cdot\|_{\kappa^{-1},i}$ in (\ref{eqn_three_local_norms}), we have
\begin{equation*}
\begin{split}
\|\mathbf{v}\|^2_{\kappa^{-1},i} &= (\kappa^{-1}\mathbf{v},\mathbf{v})_{Q,{T}^0_i} + \sum\limits_{e\in\partial T_i}\bar{\kappa}_e\llbracket q \rrbracket^2_e\\
&= \sum_{e\in\mathcal{E}^0_{T_i}}\bar{\kappa}_e\llbracket q \rrbracket^2_e + \sum\limits_{e\in\partial T_i}\bar{\kappa}_e(q)^2_{t_e}
\le 2\sum\limits_{t\in T_i}\bar{\kappa}_t q^2_t = 2(\bar{\kappa}q,q)_{T_i} = 2\|q\|^2_{\bar{\kappa},i},
\end{split}
\end{equation*}
where the subscript $t_e$ denotes the fine-grid element $t\in T_i$ such that $e\in\partial t$, $\bar{\kappa}$ is a piecewise constant function on $\mathcal{T}_h$ defined in (\ref{eqn_kappa_bar}), $\mathcal{E}^0_{T_i}$ denotes the set of all interior fine-grid edges in the partition for $T_i$. In the same way, $\|\mathbf{v}-\boldsymbol{P}_i\mathbf{v}\|_{\kappa^{-1},i} \le \|q-P_iq\|_{\bar{\kappa},i}$, and by (\ref{eqn_off_projection_approximation_no_oversampling}), we also get the approximation property for $\boldsymbol{P}_i$ as
\begin{equation}\label{eqn_projection_approximation_velocity_no_oversampling}
\|\mathbf{v}-\boldsymbol{P}_i\mathbf{v}\|_{\kappa^{-1},i} \le (\lambda^i_{l_i+1})^{-\frac{1}{2}}\|\mathbf{v}\|_{\kappa^{-1},i}, \quad i=1,2,\cdots,N_T.
\end{equation}
Then we have the following lemma:
\begin{lemma}\label{lemma_3}
	Let $\mathbf{u}_h$ be the fine-grid solution and $\mathbf{u}_{\textrm{ms}}$ be the multiscale solution on the coarse grid using the initial multiscale space, then
	\begin{equation}
	\|\mathbf{u}_h - \mathbf{u}_{\textrm{ms}}\|^2_{\kappa^{-1}} \le C_{\textrm{err}}\sum\limits^{N_T}_{i=1}\|R_i\|^2_{\scriptscriptstyle{V}}(\lambda^{i}_{l_i+1})^{-1}.\label{eqn_lemma_1_v}
	\end{equation}
	where $C_{\textrm{err}}$ is a uniform constant.
\end{lemma}
The proof of Lemma \ref{lemma_3} is similar to the proof of Lemma \ref{lemma_1}, where the approximation property (\ref{eqn_projection_approximation_velocity_no_oversampling}) will be used, and in fact the uniform constant $C_{\textrm{err}}=2$. In the following, we give the convergence analysis about the above online adaptive enrichment algorithm.
\begin{theorem}\label{theorem_4_1}
Let $\mathbf{u}_h$ be the fine-grid solution, $\mathbf{u}^m_{\textrm{ms}}$ be the multiscale solution at the $m$-th level of the online adaptive enrichment algorithm, $m\ge1$, then we have
\end{theorem}
\begin{equation}\label{eqn_theorem_2}
\|\mathbf{u}_h-\mathbf{u}^{m+1}_{\textrm{ms}}\|^2_{\kappa^{-1}} \le \Big(1-\frac{\Lambda_{\min}}{C_{\textrm{err}}}\frac{\sum_{i\in I_{\textrm{add}}}\|R_i\|^2_{\scriptscriptstyle{V}}(\lambda^i_{l_i+1})^{-1}}{\sum^{N_T}_{i=1}\|R_i\|^2_{\scriptscriptstyle{V}}(\lambda^i_{l_i+1})^{-1}}\Big) \|\mathbf{u}_h-\mathbf{u}^{m}_{\textrm{ms}}\|^2_{\kappa^{-1}},
\end{equation}
{\it where} $\Lambda_{\min} = \min_{i\in I_{\textrm{add}}}\lambda^i_{l_i+1}$,
{\it and} $I_{\textrm{add}}$ {\it is the index set of coarse elements that the corresponding online basis functions are newly added at the enrichment level $m+1$.}

{\bf Proof :} Suppose that we need to add an online basis function $\phi$ on a given coarse element $T_i$. Let $W^{m+1}_{\textrm{ms}}=W^m_{\textrm{ms}}+\textrm{span}\{\phi\}$ be the newly enriched multiscale space, and $(\mathbf{u}^{m+1}_{\textrm{ms}}, p^{m+1}_{\textrm{ms}})$ be the corresponding newly solved multiscale solution. By the relationship in (\ref{eqn_bilinear_form_equal}) and the definition of norm $\|\cdot\|_{\kappa^{-1}}$ following (\ref{eqn_global_quad_rule}), it is easy to observe that
\begin{equation*}
\|\mathbf{u}_h-\mathbf{u}^{m+1}_{\textrm{ms}}\|^2_{\kappa^{-1}} = a(p_h-p^{m+1}_{\textrm{ms}},p_h-p^{m+1}_{\textrm{ms}}),
\end{equation*}
where, the bilinear form $a(\cdot,\cdot)$ defined in (\ref{eqn_definition_apq}) is symmetric, continuous and coercive, then from (\ref{eqn_pressure_form_multiscale}), we know that the multiscale solution $\mathbf{u}^{m+1}_{\textrm{ms}}$ satisfies
\begin{equation*}
\|\mathbf{u}_h-\mathbf{u}^{m+1}_{\textrm{ms}}\|^2_{\kappa^{-1}} = \inf\limits_{q\in\tilde{p}_h+W^{m+1}_{\textrm{ms}}}a(p_h-q,p_h-q).
\end{equation*}
Taking $q=p^m_{\textrm{ms}} + \alpha_0\phi$, with $\alpha_0\in\mathbb{R}$, then we have
\begin{equation}\label{eqn_theorem_2_deduction}
\begin{split}
\|\mathbf{u}_h-\mathbf{u}^{m+1}_{\textrm{ms}}\|^2_{\kappa^{-1}} &\le a(p_h-p^m_{\textrm{ms}}-\alpha_0\phi,p_h-p^m_{\textrm{ms}}-\alpha_0\phi)\\
&=a(p_h-p^m_{\textrm{ms}},p_h-p^m_{\textrm{ms}})-2\alpha_0 a(p_h-p^m_{\textrm{ms}},\phi) + \alpha^2_0a(\phi,\phi)\\
&=\|\mathbf{u}_h-\mathbf{u}^m_{\textrm{ms}}\|^2_{\kappa^{-1}}-2\alpha_0 a(p_h-p^m_{\textrm{ms}},\phi)+\alpha^2_0a(\phi,\phi).
\end{split}
\end{equation}
According to the deduction in (\ref{eqn_theorem_2_deduction}), we select $\alpha_0={a(p_h-p^m_{\textrm{ms}},\phi)}/{a(\phi,\phi)}$ to maximize the quantity
\begin{equation*}
2\alpha_0 a(p_h-p^m_{\textrm{ms}},\phi) - \alpha^2_0a(\phi,\phi),
\end{equation*}
so as to maximize the reduction in error when the online basis function $\phi$ is added into the multiscale space,
and by the above choice of $\alpha_0$, we have
\begin{equation}\label{eqn_theorem_2_before_fine_form}
\|\mathbf{u}_h-\mathbf{u}^{m+1}_{\textrm{ms}}\|^2_{\kappa^{-1}} \le \|\mathbf{u}_h-\mathbf{u}^m_{\textrm{ms}}\|^2_{\kappa^{-1}} - \frac{|a(p_h-p^m_{\textrm{ms}},\phi)|^2}{a(\phi,\phi)}.
\end{equation}
It is evident that $\phi\in W_h(T_i)\subset W_h$, then by the use of (\ref{eqn_pressure_form}), the above inequality (\ref{eqn_theorem_2_before_fine_form}) becomes
\begin{equation*}
\|\mathbf{u}_h-\mathbf{u}^{m+1}_{\textrm{ms}}\|^2_{\kappa^{-1}} \le \|\mathbf{u}_h-\mathbf{u}^{m}_{\textrm{ms}}\|^2_{\kappa^{-1}} - \frac{|(f,\phi)-a(p^m_{\textrm{ms}},\phi)|^2}{a(\phi,\phi)}.
\end{equation*}
On the other hand, the online basis function $\phi$ actually belongs to the local snapshot space $W^{i}_{\textrm{snap}}\subset W_h(T_i)$, since we have employed the correction function $\tilde{\phi}^{i}$ to cope with the souce term $f$ on fine-grid elements in the interior of $T_i$. Let $\mathbf{v}=\psi^+|_{T_i}$, where $\psi^+$ is the velocity field obtained in (\ref{eqn_online_basis_1})-(\ref{eqn_online_basis_2}), then $\mathbf{v}\in V^{i}_{\textrm{snap}}$ is the velocity field with respect to the online basis function $\phi$. By the definition of $a(\cdot,\cdot)$ in (\ref{eqn_definition_apq}), local norm $\|\cdot\|_{\scriptscriptstyle{V_i}}$ in (\ref{eqn_snapshot_norm_2_V}), we know that
\begin{equation*}
a(\phi,\phi) = \|\mathbf{v}\|^2_{\kappa^{-1},i} = \|\phi\|_{\scriptscriptstyle{V_i}},
\end{equation*}
From definitions of the local residual operator $R_i$ in (\ref{eqn_residual_operator_no_oversampling}), the norm $\|\cdot\|_{\scriptscriptstyle{V}}$ of $R_i$ in (\ref{eqn_residual_norm_V}), and the construction of online basis function $\phi$ in (\ref{eqn_online_basis}), we we know that the online basis function $\phi$ maximize the local residual $|(f,\phi)-a(p^m_{\textrm{ms}},\phi)|^2/{a(\phi,\phi)}$ with
\begin{equation*}
|(f,\phi)-a(p^m_{\textrm{ms}},\phi)| = a(\phi,\phi) = \|\phi\|^2_{\scriptscriptstyle{V_i}} = \|R_i\|^2_{\scriptscriptstyle{V}},
\end{equation*}
Thus, we get that
\begin{equation}\label{eqn_theorem_2_inequality_1}
\|\mathbf{u}_h-\mathbf{u}^{m+1}_{\textrm{ms}}\|^2_{\kappa^{-1}} \le \|\mathbf{u}_h-\mathbf{u}^m_{\textrm{ms}}\|^2_{\kappa^{-1}} - \|R_i\|^2_{\scriptscriptstyle{V}}.
\end{equation}
Note that the initial multiscale space use the first $l_i$ offline basis functions on the coarse element $T_i$, and by the results in Lemma \ref{lemma_1}, it is easy to get that
\begin{equation}\label{eqn_theorem_2_inequality_2}
\|\mathbf{u}_h - \mathbf{u}^m_{\textrm{ms}}\|^2_{\kappa^{-1}} \le C_{\textrm{err}}\sum\limits^{N_T}_{i=1}\|R_i\|^2_{\scriptscriptstyle{V}}(\lambda^i_{l_i+1})^{-1}.
\end{equation}
Combining the above two inequalities (\ref{eqn_theorem_2_inequality_1}) and (\ref{eqn_theorem_2_inequality_2}) together, we derive that
\begin{equation*}
\|\mathbf{u}_h-\mathbf{u}^{m+1}_{\textrm{ms}}\|^2_{\kappa^{-1}} \le \Big(1-\frac{\lambda^i_{l_i+1}}{C_{\textrm{err}}}\frac{\|R_i\|^2_{\scriptscriptstyle{V}}(\lambda^i_{l_i+1})^{-1}}{\sum^{N_T}_{i=1}\|R_i\|^2_{\scriptscriptstyle{V}}(\lambda^i_{l_i+1})^{-1}}\Big) \|\mathbf{u}_h-\mathbf{u}^{m}_{\textrm{ms}}\|^2_{\kappa^{-1}}.
\end{equation*}
Let $I_{\textrm{add}}\subset\{1,2,\cdots,N_T\}$ be the index set of some non-overlapping coarse elements that online basis functions need to be added. For each index $i\in I_{\textrm{add}}$, we compute the online basis function $\phi_i$ on the corresponding coarse element $T_i$ by (\ref{eqn_online_basis_1})-(\ref{eqn_online_basis_2}) and using it to enrich the multiscale space. After adding all online basis functions $\phi_i$, $i\in I_{\textrm{add}}$, we obtain the new level multiscale space, represented as $W^{m+1}_{\textrm{ms}} = W^m_{\textrm{ms}} + \textrm{span}\{\phi_i, i\in I_{\textrm{add}}\}$. Therefore, following the above arguments and by the fact that the coarse element $T_i, i\in I_{\textrm{add}}$ are non-overlapping, we get that
\begin{equation*}
\|\mathbf{u}_h-\mathbf{u}^{m+1}_{\textrm{ms}}\|^2_{\kappa^{-1}} \le \Big(1-\sum\limits_{i\in I_{\textrm{add}}}\frac{\lambda^i_{l_i+1}}{C_{\textrm{err}}}\frac{\|R_i\|^2_{\scriptscriptstyle{V}}(\lambda^i_{l_i+1})^{-1}}{\sum^{N_T}_{i=1}\|R_i\|^2_{\scriptscriptstyle{V}}(\lambda^i_{l_i+1})^{-1}}\Big) \|\mathbf{u}_h-\mathbf{u}^{m}_{\textrm{ms}}\|^2_{\kappa^{-1}},
\end{equation*}
and by taking $\Lambda_{\min} = \min_{i\in I_{\textrm{add}}}\lambda^i_{l_i+1}$, we get the inequality (\ref{eqn_theorem_2}) and the proof is completed. $\square$

\noindent\textbf{Remark 4.1.} Form results in Theorem \ref{theorem_4_1}, we know that a faster convergence rate of the proposed online adaptive enrichment algorithm can be derived by adding more online basis functions in each enrichment level, i.e., with more coarse elements selected to add online basis functions. Moreover, the convergence rate can also be improved by choosing more offline basis functions on each coarse element in the formation of the initial multiscale space $W^1_{\textrm{ms}}$. To ensure that the error decays in a rate independent of the permeability contrast for the multiscale space enrichment with online basis functions, we need to take enough initial basis functions so that $\Lambda_{\min}$ is large enough and the initial multiscale space satisfies the following Online Error Reduction Property (ONERP) \cite{chung2015residual}:
\begin{equation*}
\frac{\Lambda_{\min}}{C_{\textrm{err}}}\frac{\sum_{i\in I_{\textrm{add}}}\|R_i\|^2_{\scriptscriptstyle{V}}(\lambda^i_{l_i+1})^{-1}}{\sum^{N_T}_{i=1}\|R_i\|^2_{\scriptscriptstyle{V}}(\lambda^i_{l_i+1})^{-1}} \ge \delta_0,
\end{equation*}
with $0<\delta_0<1$ independent of the contrast of the permeability. By taking the convergence rate as $\varepsilon_{\textrm{on}}=(1-\delta_0)$, we can obtain the following convergence for our online adaptive enrichment algorithm,
\begin{equation*}
\|\mathbf{u}_h-\mathbf{u}^{m+1}_{\textrm{ms}}\|^2_{\kappa^{-1}} \le \varepsilon_{\textrm{on}} \|\mathbf{u}_h-\mathbf{u}^{m}_{\textrm{ms}}\|^2_{\kappa^{-1}}.
\end{equation*}
In addition, we will also show numerically that when the initial multiscale space contains all offline basis functions corresponding to eigenvalues that are sensitive to the contrast of permeability, then the error will decay in a rate independent of the contrast.
\section{Numerical tests}
In this section, we present some numerical examples to demonstrate the performance of the proposed offline and online adaptive enrichment algorithms for solving the single-phase flow in high-contrast and heterogeneous porous media. Denote the fine-grid solution by $(p_h,\mathbf{u}_h)$, multiscale solution by $(p_{\textrm{ms}},\mathbf{u}_{\textrm{ms}})$, then the relative $L^2$ errors for pressure and velocity are defined respectively as
\begin{equation}
\textrm{Erp}(p_{\textrm{ms}})=\|p_{\textrm{ms}} - p_h\|\slash\|p_h\| \quad\textrm{and}\quad \textrm{Eru}(\mathbf{u}_{\textrm{ms}})=\|\mathbf{u}_{\textrm{ms}} - \mathbf{u}_h\|_{\kappa^{-1}}\slash\|\mathbf{u}_h\|_{\kappa^{-1}}.\nonumber
\end{equation}
We will test the offline adaptive enrichment algorithm and online adaptive enrichment algorithm, respectively, in the following two subsections.
\subsection{Offline enrichment tests}
In this subsection, we investigate the performance of our offline adaptive enrichment algorithm with the proposed error indicator in (\ref{eqn_eta_offline}), to make it clear, we denote this indicator by $\eta^{\scriptscriptstyle{W}}_i$, that is,
\begin{equation*}
\eta^{\scriptscriptstyle{W}}_i = \|R_i\|^2(\lambda^{i}_{l_i+1})^{-1}, \quad i=1,2,\cdots,N_T,
\end{equation*}
and for the comparison purpose, we also test the performance of the offline adaptive enrichment algorithm induced by the exact indicator $\eta^{\textrm{ex}}$, defined as
\begin{equation*}
\eta^{\textrm{ex}}_i = \|\mathbf{u}_h-\mathbf{u}_{\textrm{ms}}\|^2_{\kappa^{-1},i} \ \ \ \ \textrm{on}\ \ T_i\in\mathcal{T}_H,\quad i=1,2,\cdots,N_T.
\end{equation*}
We use $(\mathbf{u}_{\textrm{ms},{\scriptscriptstyle{W}}},p_{\textrm{ms},{\scriptscriptstyle{W}}})$, $(\mathbf{u}_{\textrm{ms,ex}},p_{\textrm{ms,ex}})$ to denote the corresponding multiscale solutions of the offline adaptive method with error indicators $\eta^{\scriptscriptstyle{W}}$ and $\eta^{\textrm{ex}}$, respectively, and use  $(\mathbf{u}_{\textrm{ms,um}},p_{\textrm{ms,um}})$ to denote the multiscale solution obtained by the uniform enrichment with offline basis functions. In the following examples of this subsection, we take the parameter $\theta_{\textrm{off}}=0.7$ and add one offline basis function per coarse element selected by the criterion (\ref{eqn_criterion_offline}) at each enrichment iteration. The oversampling size is set to be two fine-grid layers.\\

{\bf Example 1:}
\begin{figure}[h!]
	\centering
	\includegraphics[width=0.49\textwidth]{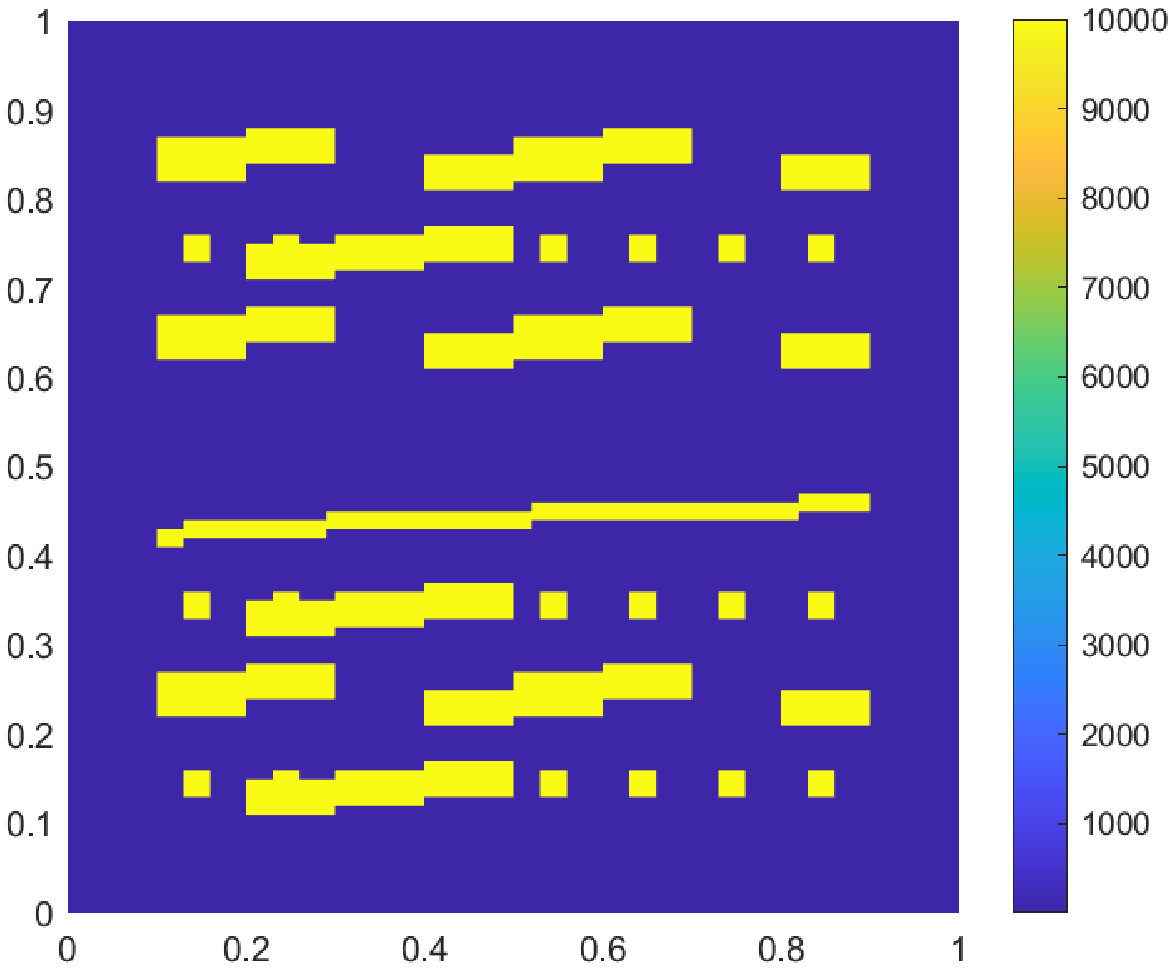}
	\includegraphics[width=0.49\textwidth]{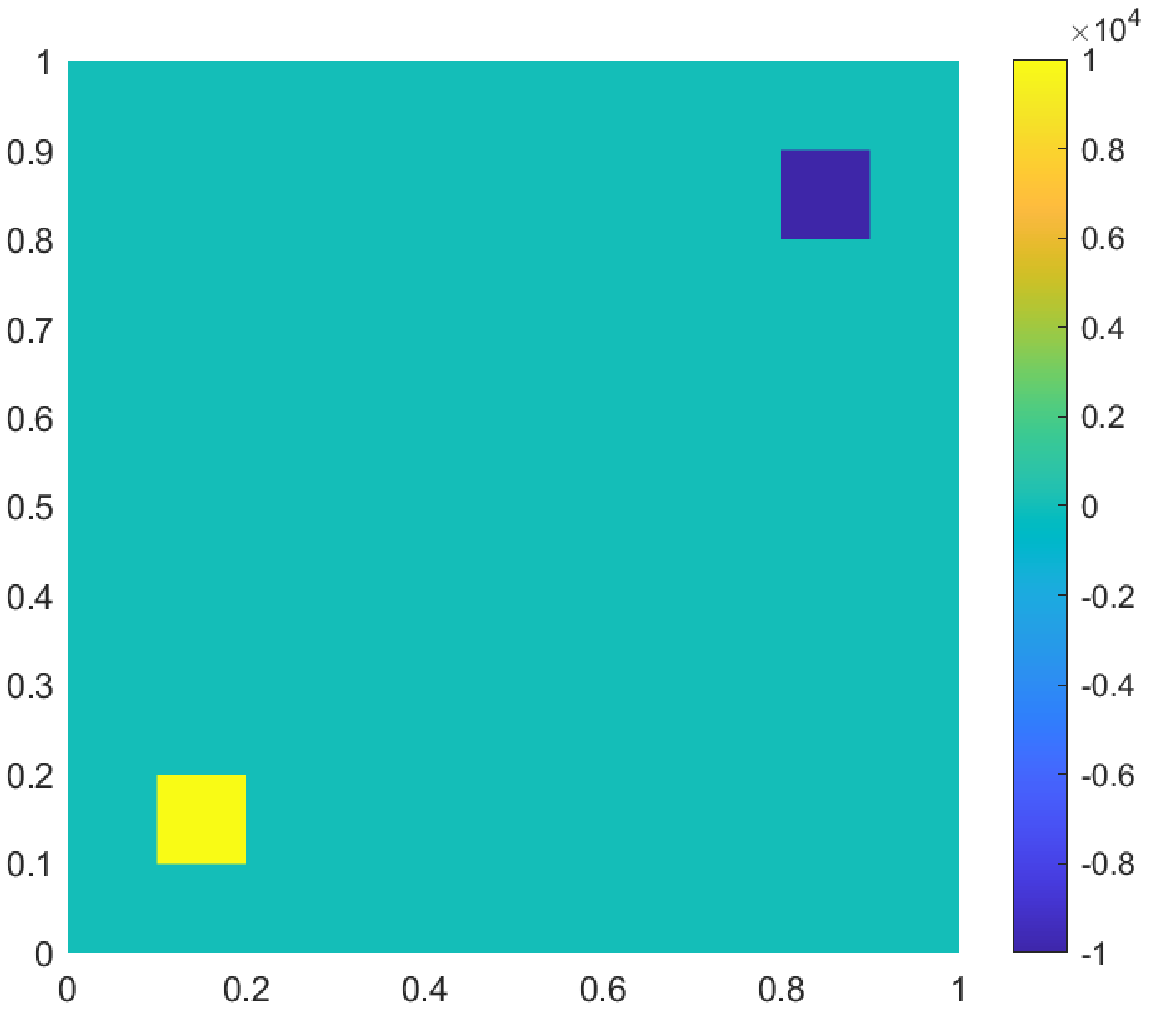}
	\caption{The distribution of permeability field $\kappa$ (left) and source term $f$ (right) in Example 1.}
	\label{fig_perm_1}
\end{figure}
The computational domain is taken to be a square $\Omega=[0, 1]^2$, the boundary condition is set to be $p=0$ on $x=0$, $x=1$, and homogeneous Neumann boundary condition $\mathbf{u}\cdot\mathbf{n}=0$ on $y=0$, $y=1$. The permeability field $\kappa$ and the source term $f$ in (\ref{eqn_model_2}) are shown in the left and right graph of Figure \ref{fig_perm_1}, respectively. The fine grid is taken to be a $100\times100$ uniform mesh, and the coarse grid is a $10\times10$ uniform mesh with $N_T=100$ coarse elements.

In Table \ref{tab_err_ex1_pressure_velocity}, we show the relative errors $\textrm{Erp}(p_{\textrm{ms},{\scriptscriptstyle{W}}})$, $\textrm{Erp}(p_{\textrm{ms,ex}})$, $\textrm{Erp}(p_{\textrm{ms,um}})$ for pressure, and relative errors $\textrm{Eru}(\mathbf{u}_{\textrm{ms},{\scriptscriptstyle{W}}})$, $\textrm{Eru}(\mathbf{u}_{\textrm{ms,ex}})$, $\textrm{Eru}(\mathbf{u}_{\textrm{ms,um}})$ for velocity, with respect to three initial offline basis functions per coarse element, where the symbol '$\#$Dofs' denotes the total number of basis functions that the multiscale space $W_{\textrm{ms}}$ used. It is observed that with the enrichment of the multiscale space, the relative errors of offline adaptive enrichment algorithms induced by indicators $\eta^{\scriptscriptstyle{W}}$ and $\eta^{\textrm{ex}}$ are similar and about half of relative errors of the offline uniform enrichment algorithm with similar total dimensions of the multiscale space, which indicates that the offline adaptive enrichment algorithm with the indicator $\eta^{\scriptscriptstyle{W}}$ is effective, and we only need a smaller number of offline basis functions to obtain the same relative errors compared with the offline uniform enrichment algorithm.
\begin{table}[h!]
	\caption{(Example 1) Relative errors $\textrm{Erp}(p_{\textrm{ms},{\scriptscriptstyle{\textit{W}}}})$, $\textrm{Erp}(p_{\textrm{ms,ex}})$, $\textrm{Erp}(p_{\textrm{ms,um}})$ for pressure and $\textrm{Eru}(\mathbf{u}_{\textrm{ms},{\scriptscriptstyle{\textit{W}}}})$, $\textrm{Eru}(\mathbf{u}_{\textrm{ms,ex}})$, $\textrm{Eru}(\mathbf{u}_{\textrm{ms,um}})$ for velocity of the offline enrichment with three initial offline basis functions per coarse element, $\theta_{\textrm{off}}=0.7$.}\label{tab_err_ex1_pressure_velocity}
	\centering
	\begin{tabular*}
		{1.00\textwidth}{@{\extracolsep{\fill}}|c|c|c||c|c|c||c|c|c|}
		\hline
		$\#$Dofs &$\textrm{Erp}(p_{\textrm{ms},\scriptscriptstyle{W}})$ &$\textrm{Eru}(\mathbf{u}_{\textrm{ms},\scriptscriptstyle{W}})$ &$\#$Dofs &$\textrm{Erp}(p_{\textrm{ms,ex}})$ &$\textrm{Eru}(\mathbf{u}_{\textrm{ms,ex}})$ &$\#$Dofs &$\textrm{Erp}(p_{\textrm{ms,um}})$ &$\textrm{Eru}(\mathbf{u}_{\textrm{ms,um}})$\\
		\hline
		300    &0.3208  &0.3370        &300    &0.3208  &0.3370        &300    &0.3208  &0.3370 \\
		416    &0.0359  &0.0673        &408    &0.0448  &0.0661        &400    &0.0786  &0.1392 \\
		486    &0.0175  &0.0392        &497    &0.0169  &0.0323        &500    &0.0236  &0.0668 \\
		584    &0.0072  &0.0211        &597    &0.0075  &0.0190        &600    &0.0117  &0.0491 \\
		696    &0.0036  &0.0136        &691    &0.0034  &0.0121        &700    &0.0066  &0.0254 \\
		802    &0.0018  &0.0092        &822    &0.0013  &0.0060        &800    &0.0038  &0.0176 \\
		\hline
	\end{tabular*}
\end{table}

In Figure \ref{fig_num_1}, we display distributions of the number of offline basis functions on each coarse element for multiscale spaces whose total dimensions are around $600$ $(N_T\times 6)$ with the proposed error indicators $\eta^{\scriptscriptstyle{W}}$ and $\eta^{\textrm{ex}}$. In the left and right graphs, the total dimensions of multiscale spaces are $584$ and $597$ after $11$ and $8$ iterations of the offline adaptive enrichment, and the associated relative errors for velocity are $0.0211$ and $0.0190$, with respect to indicators $\eta^{\scriptscriptstyle{W}}$ and $\eta^{\textrm{ex}}$, respectively. We observe that the dimension distributions induced by the proposed indicator $\eta^{\scriptscriptstyle{W}}$ and the exact indicator $\eta^{\textrm{ex}}$ follow a slightly different pattern, where the distribution induced by indicator $\eta^{\scriptscriptstyle{W}}$ has more number of basis functions on the boundary coarse elements, and the distribution induced by indicator $\eta^{\textrm{ex}}$ has more number of basis functions on the interior coarse elements. Figure \ref{fig_pressure_velocity_1} presents the corresponding multiscale solutions $(\mathbf{u}_{\textrm{ms},\scriptscriptstyle{W}},p_{\textrm{ms},\scriptscriptstyle{W}})$ and also the fine-grid solution $(\mathbf{u}_h,p_h)$ for comparison, in these graphs, the streamlines of velocity are plotted. We can see that the multiscale solution obtained by the offline adaptive enrichment algorithm using indicator $\eta^{\scriptscriptstyle{W}}$ has a good approximation of the fine-grid solution.
\begin{figure}[h!]
	\centering
	\includegraphics[width=0.40\textwidth]{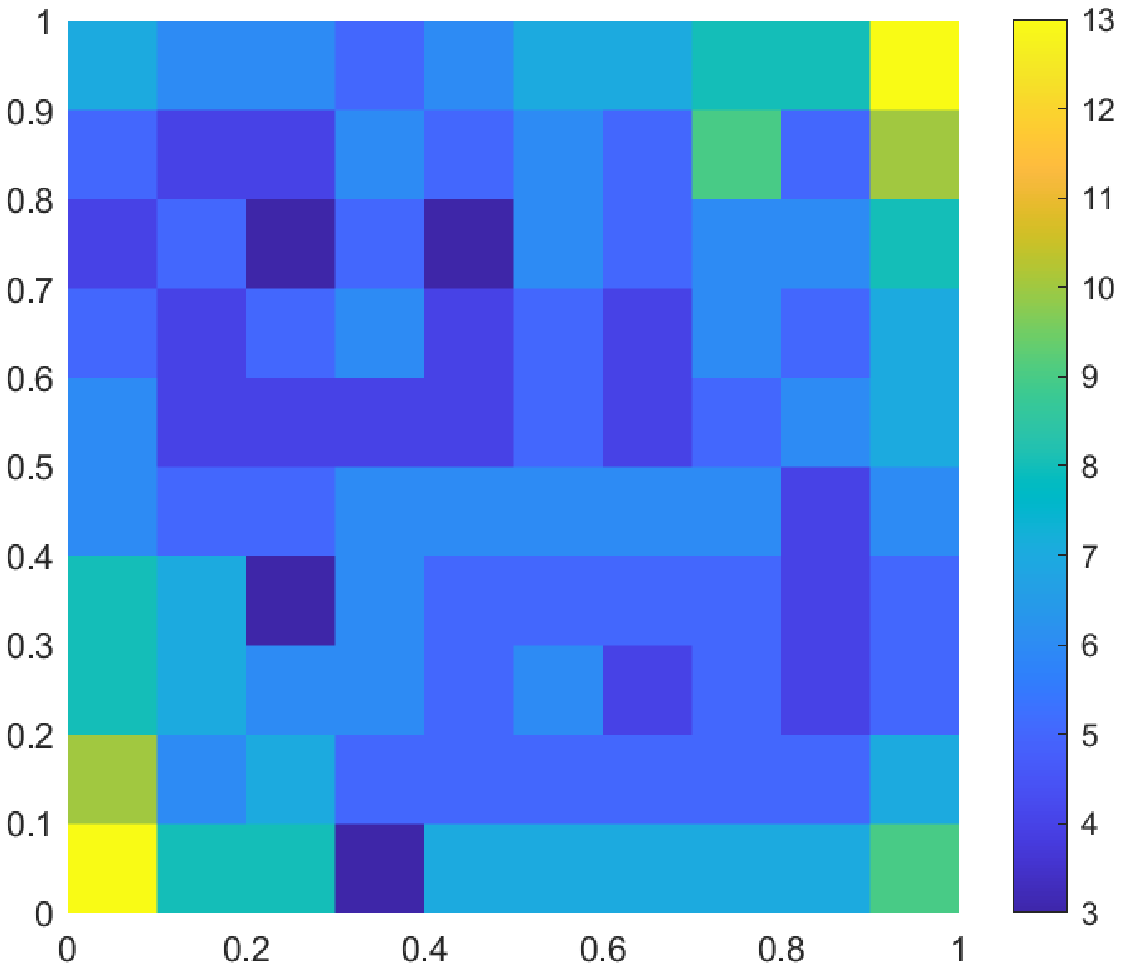}
	\includegraphics[width=0.40\textwidth]{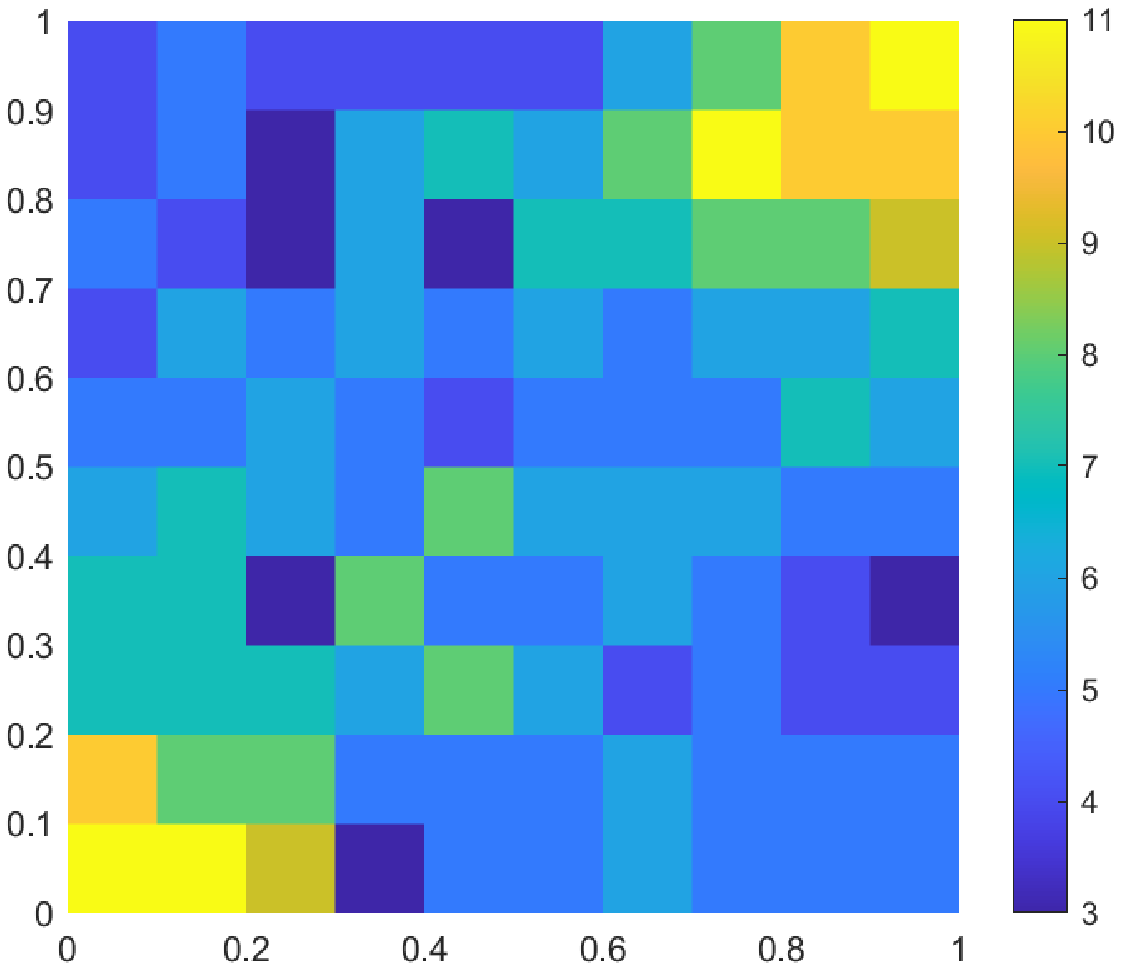}
	\caption{(Example 1) The dimension distributions for multiscale spaces with total dimensions around $N_T\times 6$ of the offline adaptive method. Left: Indicator $\eta^{\scriptscriptstyle{\textit{W}}}$, after $11$ iterations, $\#$Dofs$=584$. Right: Indicator $\eta^{\textrm{ex}}$, after $8$ iterations, $\#$Dofs$=597$.}
	\label{fig_num_1}
\end{figure}
\begin{figure}[h!]
	\centering
	\includegraphics[width=0.40\textwidth]{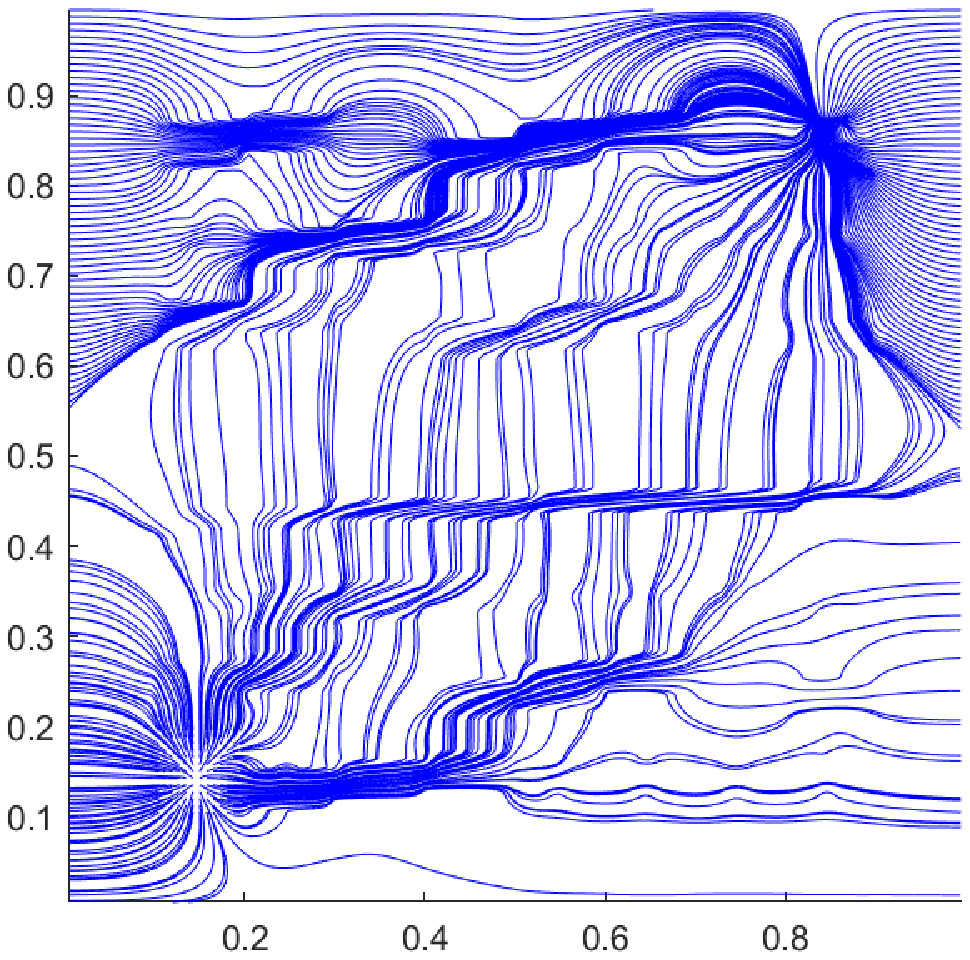}
	\includegraphics[width=0.40\textwidth]{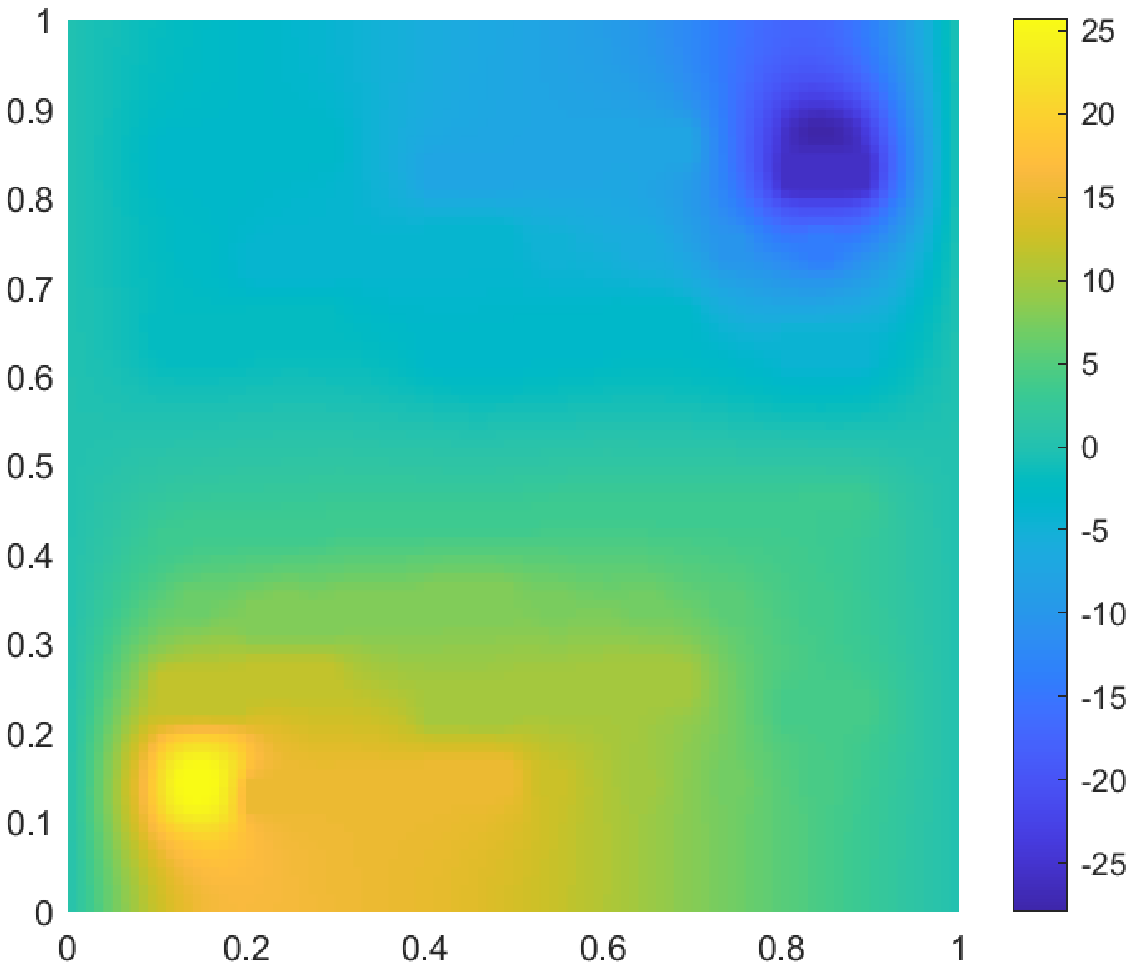}\\
	\includegraphics[width=0.40\textwidth]{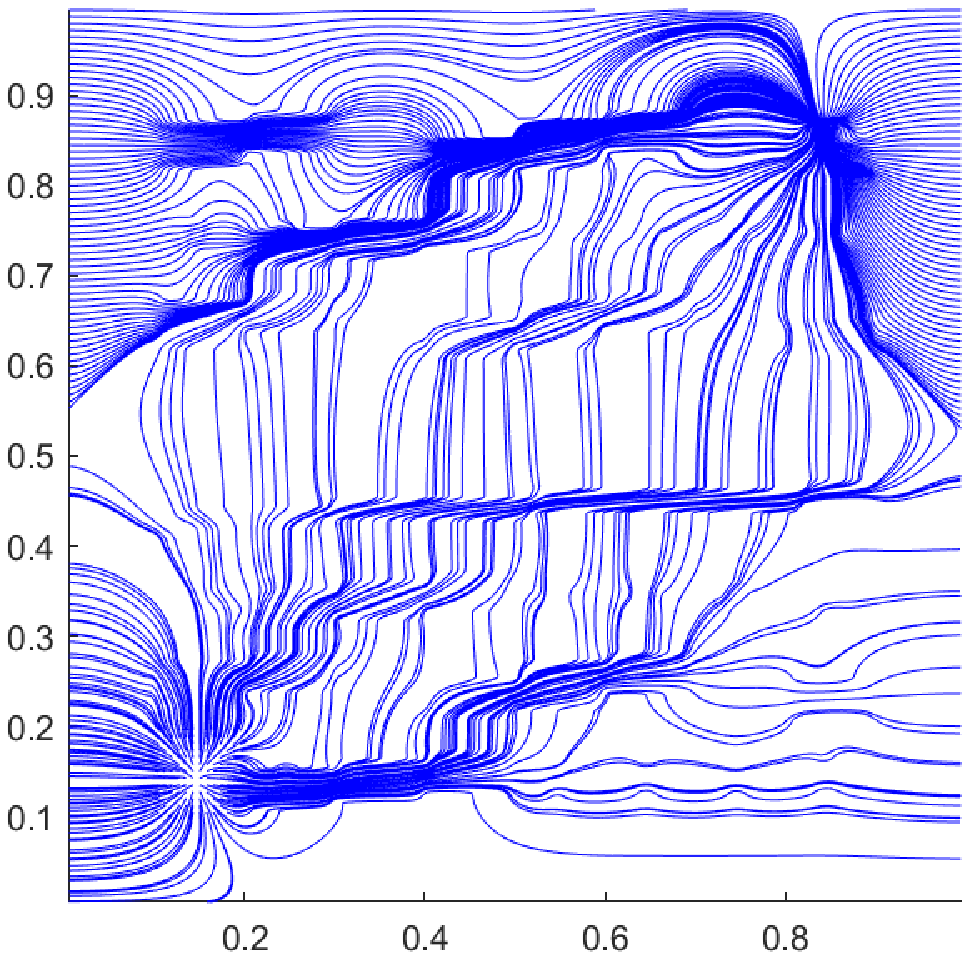}
	\includegraphics[width=0.40\textwidth]{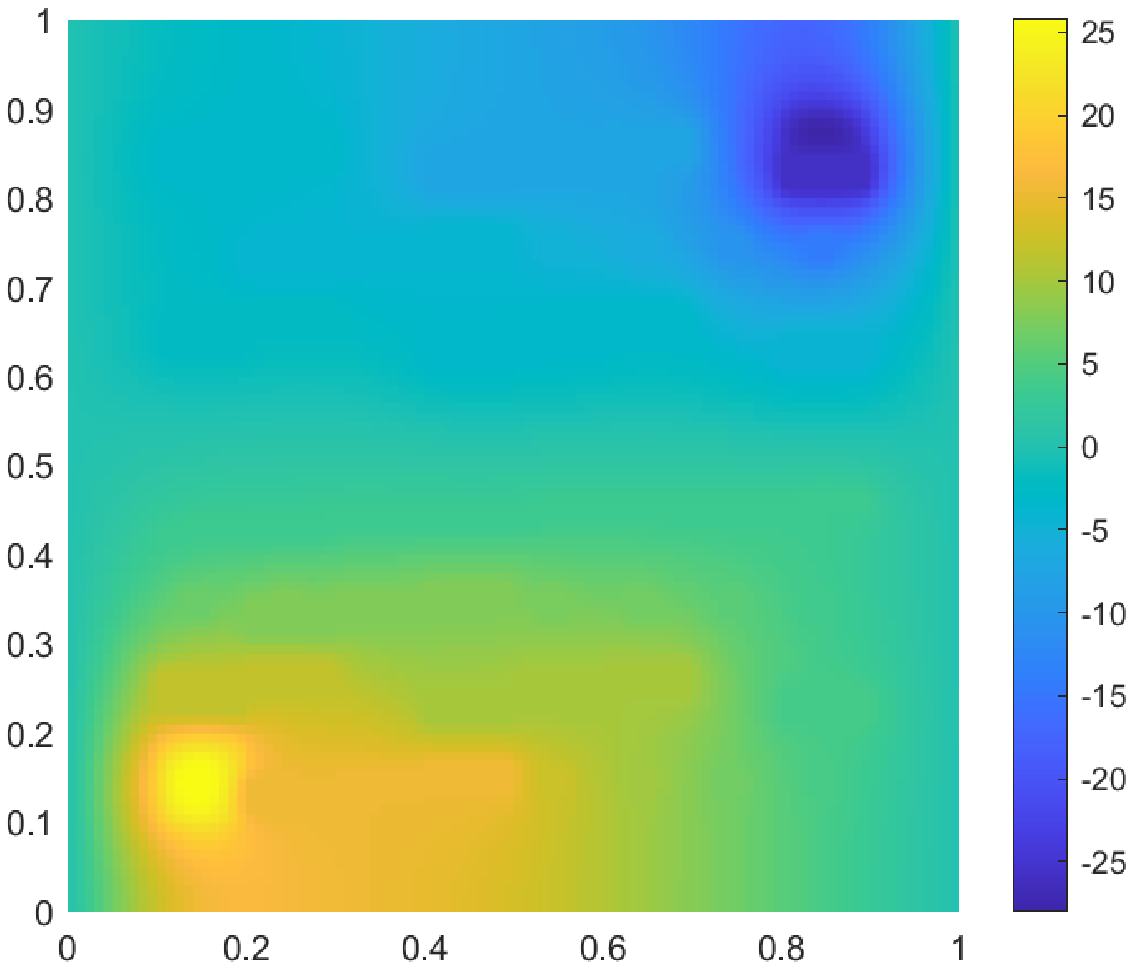}
	\caption{(Example 1) Velocity and pressure. Top: Multiscale solution $(\mathbf{u}_{\textrm{ms},{\scriptscriptstyle{\textit{W}}}},p_{\textrm{ms},{\scriptscriptstyle{\textit{W}}}})$, after $11$ iterations, $\#$Dofs$=584$, $\textrm{Eru}(\mathbf{u}_{\textrm{ms},{\scriptscriptstyle{\textit{W}}}})=0.0211$. Bottom: Fine-grid solution $(\mathbf{u}_h,p_h)$.}
	\label{fig_pressure_velocity_1}
\end{figure}

{\bf Example 2:}
\begin{figure}[h!]
	\centering
	\includegraphics[width=0.96\textwidth]{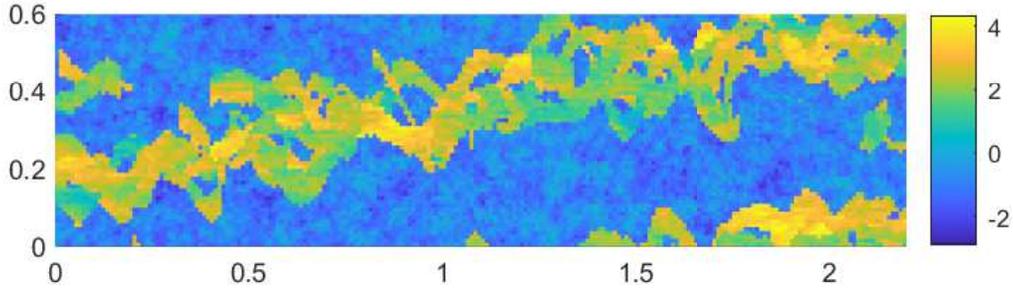}
	\caption{The distribution of permeability field $\kappa$ in logarithmic scale in Example 2.}
	\label{fig_perm_2}
\end{figure}
The computatinal domain is set to be $\Omega=[0,2.2]\times[0,0.6]$, the boundary condition is assumed to be the type of an open-side boundary condition with $p=1$ on $x=0$ and $p=0$ on $x=2.2$, and homogeneous Neumann boundary condition $\mathbf{u}\cdot\mathbf{n}=0$ on the other boundaries. The source term $f$ in (\ref{eqn_model_2}) is taken to be zero. The fine grid is a $220\times60$ uniform mesh, and the coarse grid is a $22\times6$ uniform mesh with $N_T=132$ coarse elements. The permeability field $\kappa$ is a part of the horizontal permeability from the SPE10 data set, as shown in Figure \ref{fig_perm_2}.

Table \ref{tab_err_ex2_pressure_velocity} displays the relative errors $\textrm{Erp}(p_{\textrm{ms},{\scriptscriptstyle{W}}})$, $\textrm{Erp}(p_{\textrm{ms,ex}})$, $\textrm{Erp}(p_{\textrm{ms,um}})$ for pressure, and the relative errors $\textrm{Eru}(\mathbf{u}_{\textrm{ms},{\scriptscriptstyle{W}}})$, $\textrm{Eru}(\mathbf{u}_{\textrm{ms,ex}})$ $\textrm{Eru}(\mathbf{u}_{\textrm{ms,um}})$ for velocity, with respect to three initial offline basis functions per coarse element. We observe semblable results as shown in Example 1. We can see that for multiscale spaces with similar dimensions, the relative errors of the offline adaptive enrichment with indicators $\eta^{\scriptscriptstyle{W}}$ and $\eta^{\textrm{ex}}$ are great smaller than the relative errors of the offline uniform enrichment, and the performances of indicators $\eta^{\scriptscriptstyle{W}}$ and $\eta^{\textrm{ex}}$ are similar. Totally speaking, the offline adaptive method with the proposed error indicator $\eta^{\scriptscriptstyle{W}}$ is effective and reliable that can improve the accuracy of the multiscale solution greatly with smaller number of basis functions than the uniform enrichment with offline basis functions.
\begin{table}[h!]
	\caption{(Example 2) Relative errors $\textrm{Erp}(p_{\textrm{ms},{\scriptscriptstyle{\textit{W}}}})$, $\textrm{Erp}(p_{\textrm{ms,ex}})$, $\textrm{Erp}(p_{\textrm{ms,um}})$ for pressure and $\textrm{Eru}(\mathbf{u}_{\textrm{ms},{\scriptscriptstyle{\textit{W}}}})$, $\textrm{Eru}(\mathbf{u}_{\textrm{ms,ex}})$, $\textrm{Eru}(\mathbf{u}_{\textrm{ms,um}})$ for velocity of the offline enrichment with three initial offline basis functions on each coarse element, $\theta_{\textrm{off}}=0.7$.}\label{tab_err_ex2_pressure_velocity}
	\centering
	\begin{tabular*}
		{1.00\textwidth}{@{\extracolsep{\fill}}|c|c|c||c|c|c||c|c|c|}
		\hline
		$\#$Dofs &$\textrm{Erp}(p_{\textrm{ms},\scriptscriptstyle{W}})$ &$\textrm{Eru}(\mathbf{u}_{\textrm{ms},\scriptscriptstyle{W}})$ &$\#$Dofs &$\textrm{Erp}(p_{\textrm{ms,ex}})$ &$\textrm{Eru}(\mathbf{u}_{\textrm{ms,ex}})$ &$\#$Dofs &$\textrm{Erp}(p_{\textrm{ms,um}})$ &$\textrm{Eru}(\mathbf{u}_{\textrm{ms,um}})$ \\
		\hline
		396    &0.0630     &0.4897        &396    &0.0630     &0.4897        &396    &0.0630  &0.4897 \\
		522    &0.0070     &0.0308        &522    &0.0104     &0.0410        &528    &0.0207  &0.2038 \\
	    654    &0.0020     &0.0148        &644    &0.0036     &0.0145        &660    &0.0079  &0.0746 \\
		794    &0.0011     &0.0091        &812    &0.0026     &0.0077        &792    &0.0066  &0.0618 \\
		931    &7.8298e-4  &0.0061        &908    &0.0017     &0.0056        &924    &0.0047  &0.0436 \\
		1062   &5.4147e-4  &0.0039        &1061   &6.8010e-4  &0.0028        &1056   &0.0022  &0.0254 \\
		\hline
	\end{tabular*}
\end{table}

Likewise, we also display distributions of the number of basis functions in Figure \ref{fig_num_2} for multiscale spaces whose total dimensions are around $792$ $(N_T\times 6)$ induced by indicators $\eta^{\scriptscriptstyle{W}}$ and $\eta^{\textrm{ex}}$. In the left and right graphs, the total dimensions of the multiscale spaces are $794$ and $812$ after $19$ and $14$ iterations of the offline adaptive enrichment, and the corresponding relative errors for velocity are $0.0091$ and $0.0077$, with respect to indicators $\eta^{\scriptscriptstyle{W}}$ and $\eta^{\textrm{ex}}$, respectively. Similar results with Example 1 are also observed, we find that the dimension distribution induced by indicator $\eta^{\scriptscriptstyle{V}}$ has more number of basis functions on the boundary coarse elements. Figure \ref{fig_pressure_velocity_2} compares the corresponding multiscale solution $(\mathbf{u}_{\textrm{ms},\scriptscriptstyle{W}},p_{\textrm{ms},\scriptscriptstyle{W}})$ and the fine-grid solution $(\mathbf{u}_h,p_h)$, we can see that by the use of error indicator $\eta^{\scriptscriptstyle{W}}$, the offline adaptive method has a competitive performance.
\begin{figure}[h!]
	\mbox{\hspace{-0.90cm}}
	\begin{minipage}[b]{0.52\textwidth}
		\centering
		\includegraphics[scale=0.66]{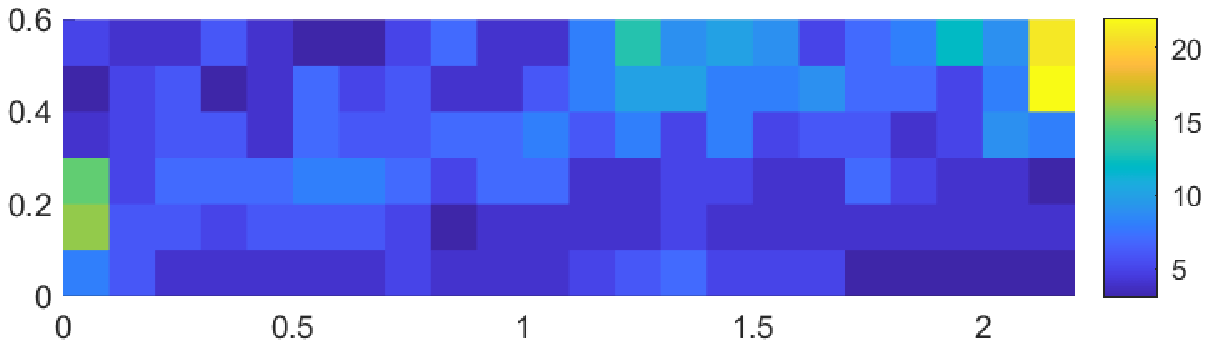}
	\end{minipage}
	\mbox{\hspace{0.00cm}}
	\begin{minipage}[b]{0.52\textwidth}
		\centering
		\includegraphics[scale=0.66]{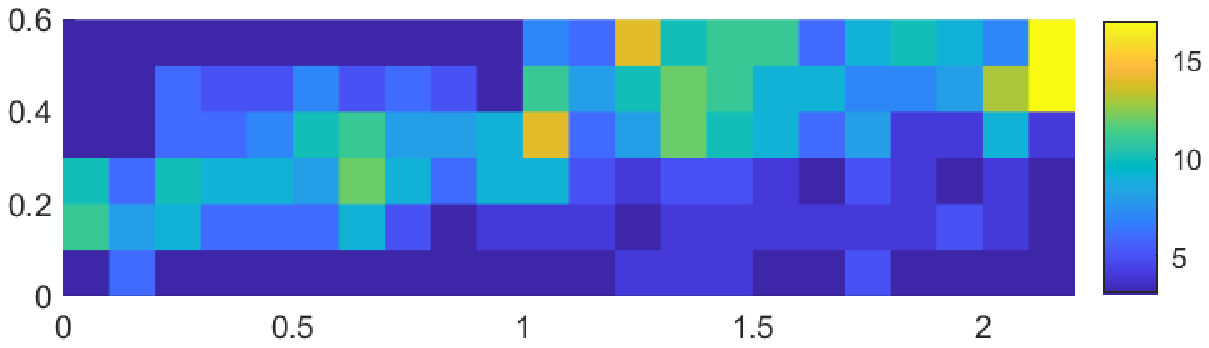}
	\end{minipage}
	\caption{(Example 2) The dimension distributions for multiscale spaces with total dimensions around $N_T\times 6$ of the offline adaptive method. Left: Indicator $\eta^{\scriptscriptstyle{\textit{W}}}$, after $19$ iterations, $\#$Dofs$=794$. Right: Indicator $\eta^{\textrm{ex}}$, after $14$ iterations, $\#$Dofs$=812$.}
	\label{fig_num_2}
\end{figure}
\begin{figure}[h!]
	\mbox{\hspace{-0.90cm}}
	\begin{minipage}[b]{0.52\textwidth}
		\centering
		\includegraphics[scale=0.66]{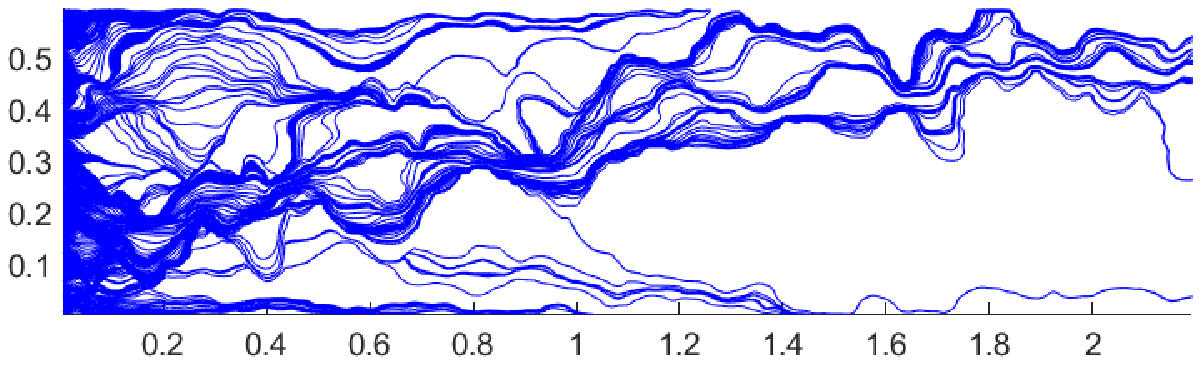}
	\end{minipage}
	\mbox{\hspace{0.00cm}}
	\begin{minipage}[b]{0.52\textwidth}
		\centering
		\includegraphics[scale=0.66]{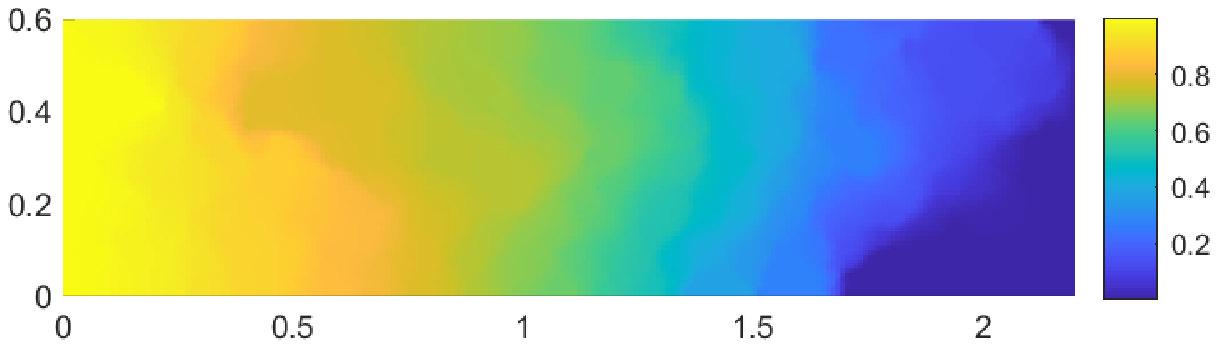}
	\end{minipage}\\
	\mbox{\hspace{-0.90cm}}
	\begin{minipage}[b]{0.52\textwidth}
		\centering
		\includegraphics[scale=0.66]{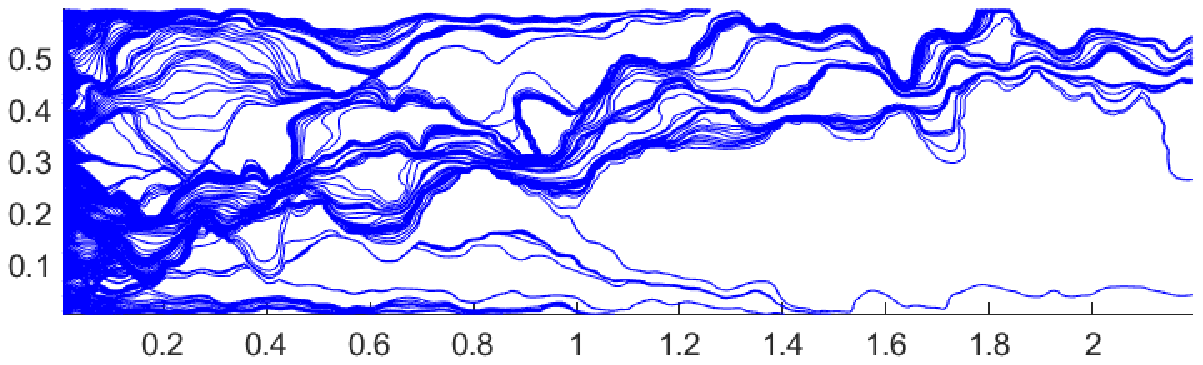}
	\end{minipage}
	\mbox{\hspace{0.00cm}}
	\begin{minipage}[b]{0.52\textwidth}
		\centering
		\includegraphics[scale=0.66]{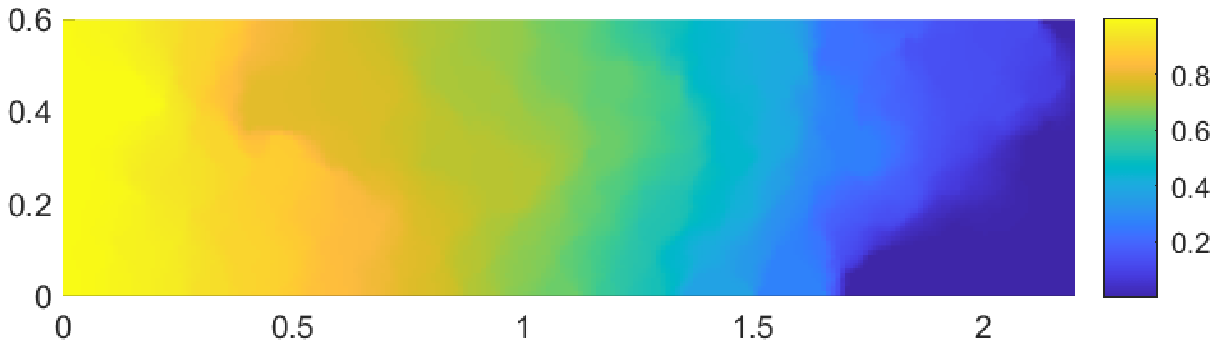}
	\end{minipage}
	\caption{(Example 2) Velocity and pressure. Top: Multiscale solution $(\mathbf{u}_{\textrm{ms},{\scriptscriptstyle{\textit{W}}}},p_{\textrm{ms},{\scriptscriptstyle{\textit{W}}}})$, after $19$ iterations, $\#$Dofs$=794$, $\textrm{Eru}(\mathbf{u}_{\textrm{ms},{\scriptscriptstyle{\textit{W}}}})=0.0091$. Bottom: Fine-grid solution $(\mathbf{u}_h,p_h)$.}
	\label{fig_pressure_velocity_2}
\end{figure}
\subsection{Online enrichment tests}
In this subsection, we investigate the performance of the online enrichment algorithm. We will test two situations respectively that the multiscale space is enriched uniformly for all coarse elements and enriched adaptively for selected coarse elements based on residuals. The online basis functions are calculated and added into the multiscale space only in disjoint regions at a time, for convenience, we use a two-index notation to enumerate all coarse elements, i.e., the coarse elements are indexed by $T_{ij}$, with $i=1,2,\cdots,N_x$ and $j=1,2,\cdots,N_y$, where $N_x$ and $N_y$ are the number of partitions of the coarse grid $\mathcal{T}_H$ along the $x$ and $y$ directions, respectively. Let $I_x=\{1,2,\cdots,N_x\}$ and $I_y=\{1,2,\cdots,N_y\}$. we denote $I_{x,1}$, $I_{x,2}$ be the subsets composed of the odd, even indices of $I_x$ respectively, and $I_{y,1}$, $I_{y,2}$ be the subsets composed of the odd, even indices of $I_y$ respectively. We can separate all coarse elements into four disjoint subsets $I_1$, $I_2$, $I_3$ and $I_4$, respectively, $I_1=I_{x,1}\times I_{y,1}$, $I_2=I_{x,1}\times I_{y,2}$, $I_3=I_{x,2}\times I_{y,1}$ and $I_4=I_{x,2}\times I_{y,2}$. Each iteration of the online multiscale space enrichment contains four subiterations, in particular, these four subiterations are defined by adding online basis functions with respect to coarse elements $T_{ij}\in I_1$, $T_{ij}\in I_2$, $T_{ij}\in I_3$ and $T_{ij}\in I_4$, respectively.

Firstly, we test the situation that the multiscale space is uniformly enriched by adding one online basis function per coarse element at each enrichment iteration.

{\bf Example 3:} In Table \ref{tab_err_ex1_online_uniform} and Table \ref{tab_err_ex2_online_uniform}, we present the relative errors for pressure and velocity of the online uniform enrichment using data in Example 1 and Example 2, respectively, where different number of initial basis functions (the first one, three and five offline basis functions) per coarse element are tested. It is evident that the accuracy of the multiscale solution is improved a lot by several iterations of the uniform enrichment with online basis functions. Comparing results of the online uniform enrichment using three initial basis functions in the middle columns of Table \ref{tab_err_ex1_online_uniform} with results of the offline uniform enrichment in the last columns of Tables \ref{tab_err_ex1_pressure_velocity}, we observe that the relative errors of online uniform enrichment decay more quickly than the offline uniform enrichment. Meanwhile, we obtain the same observation from the comparison of results in the middle columns of Table \ref{tab_err_ex2_online_uniform} by the online uniform enrichment and results in the last columns of Tables \ref{tab_err_ex2_pressure_velocity} by the offline uniform enrichment for the use of data in Example 2. We come to the conclusion that the online basis functions are able to behave better than the offline basis functions.

To compare the performance of the online enrichment with different number of initial basis functions evidently, the convergence histories are also plotted in the top of Figure \ref{fig_uniform_online_velocityerr_12} for the use of data in Example 1 and in the bottom of Figure \ref{fig_uniform_online_velocityerr_12} for the use of data in Example 2, where the relative error and the logarithm of the relative error for velocity against the dimensions of the multiscale space $W_{\textrm{ms}}$ at each enrichment iteration are depicted. It can be observed that, for both data cases, the convergence rates become faster by utilizing more initial basis functions per coarse element, which conform with the convergence analysis in Theorem \ref{theorem_4_1}.
\begin{table}[h!]
	\caption{(Example 3) Relative errors $\textrm{Erp}(p_{\textrm{ms}})$, $\textrm{Eru}(\mathbf{u}_{\textrm{ms}})$ of the online uniform enrichment with different number of initial basis functions per coarse element using data in Example 1. Left: One initial basis function, $\Lambda_{\min}=0.0146$. Middle: Three initial basis functions, $\Lambda_{\min}=0.0337$. Right: Five initial basis functions, $\Lambda_{\min}=0.0526$.}\label{tab_err_ex1_online_uniform}
	\centering
	\begin{tabular*}
		{0.93\textwidth}{@{\extracolsep{\fill}}|c|c|c||c|c|c||c|c|c|}
		\hline
		\multicolumn{3}{|c||}{Initial Dofs per $T=1$}  &\multicolumn{3}{c||}{Initial Dofs per $T=3$} &\multicolumn{3}{c|}{Initial Dofs per $T=5$}\\
		\hline
		$\#$Dofs &$\textrm{Erp}(p_{\textrm{ms}})$ &$\textrm{Eru}(\mathbf{u}_{\textrm{ms}})$ &$\#$Dofs &$\textrm{Erp}(p_{\textrm{ms}})$ &$\textrm{Eru}(\mathbf{u}_{\textrm{ms}})$ &$\#$Dofs &$\textrm{Erp}(p_{\textrm{ms}})$ &$\textrm{Eru}(\mathbf{u}_{\textrm{ms}})$\\
		\hline
		100   &0.8803    &0.7932               &-     &-         &-                    &-     &-         &-      \\
		300   &0.0762    &0.0953               &300   &0.2940    &0.3389               &-     &-         &-      \\
		500   &0.0033    &0.0105               &500   &0.0068    &0.0205               &500   &0.1006    &0.1555 \\
		600   &9.2967e-4 &0.0029               &600   &0.0012    &0.0044               &600   &0.0133    &0.0253 \\
		700   &2.3171e-4 &7.4302e-4            &700   &2.6200e-4 &9.4020e-4            &700   &7.1004e-4 &0.0037 \\
		800   &4.0392e-5 &1.9977e-4            &800   &7.2793e-5 &2.5491e-4            &800   &8.8408e-5 &5.1306e-4 \\
		900   &1.0080e-5 &7.5486e-5            &900   &8.5836e-6 &3.7774e-5            &900   &1.1875e-5 &6.9438e-5 \\
		1000  &3.4558e-6 &2.2781e-5            &1000  &1.4129e-6 &6.8889e-6            &1000  &1.4402e-6 &9.2594e-6 \\
		\hline
	\end{tabular*}
\end{table}
\begin{table}[h!]
	\caption{(Example 3) Relative errors $\textrm{Erp}(p_{\textrm{ms}})$, $\textrm{Eru}(\mathbf{u}_{\textrm{ms}})$ of online uniform enrichment with different number of initial basis functions per coarse element using data in Example 2. Left: One initial basis function, $\Lambda_{\min}=1.2779e-5$. Middle: Three initial basis functions, $\Lambda_{\min}=0.0051$. Right: Five initial basis functions, $\Lambda_{\min}=0.0280$.}\label{tab_err_ex2_online_uniform}
	\centering
	\begin{tabular*}
		{0.93\textwidth}{@{\extracolsep{\fill}}|c|c|c||c|c|c||c|c|c|}
		\hline
		\multicolumn{3}{|c||}{Initial Dofs per $T=1$}  &\multicolumn{3}{c||}{Initial Dofs per $T=3$} &\multicolumn{3}{c|}{Initial Dofs per $T=5$}\\
		\hline
		$\#$Dofs &$\textrm{Erp}(p_{\textrm{ms}})$ &$\textrm{Eru}(\mathbf{u}_{\textrm{ms}})$ &$\#$Dofs &$\textrm{Erp}(p_{\textrm{ms}})$ &$\textrm{Eru}(\mathbf{u}_{\textrm{ms}})$ &$\#$Dofs &$\textrm{Erp}(p_{\textrm{ms}})$ &$\textrm{Eru}(\mathbf{u}_{\textrm{ms}})$\\
		\hline
		132   &0.1194    &1.1716               &-     &-         &-                    &-     &-         &-      \\
		396   &0.0411    &0.3000               &396   &0.0487    &0.4524               &-     &-         &-      \\
		660   &0.0143    &0.0970               &660   &0.0026    &0.0412               &660   &0.0074    &0.1743 \\
		792   &0.0106    &0.0748               &792   &7.4940e-4 &0.0119               &792   &0.0032    &0.0549 \\
		924   &0.0047    &0.0372               &924   &2.7816e-4 &0.0028               &924   &5.9840e-4 &0.0100 \\
		1056  &0.0011    &0.0085               &1056  &8.6825e-5 &7.4109e-4            &1056  &9.4250e-5 &0.0016 \\
		1188  &5.4700e-4 &0.0020               &1188  &1.9785e-5 &2.0897e-4            &1188  &1.2473e-5 &2.2622e-4 \\
		1320  &1.7744e-4 &5.2569e-4            &1320  &6.5514e-6 &4.9012e-5            &1320  &3.2814e-6 &4.2655e-5 \\
		\hline
	\end{tabular*}
\end{table}
\begin{figure}[h!]
	\centering
	\mbox{\hspace{0.00cm}}
	\begin{minipage}[b]{0.44\textwidth}
		\centering
		\includegraphics[scale=0.555]{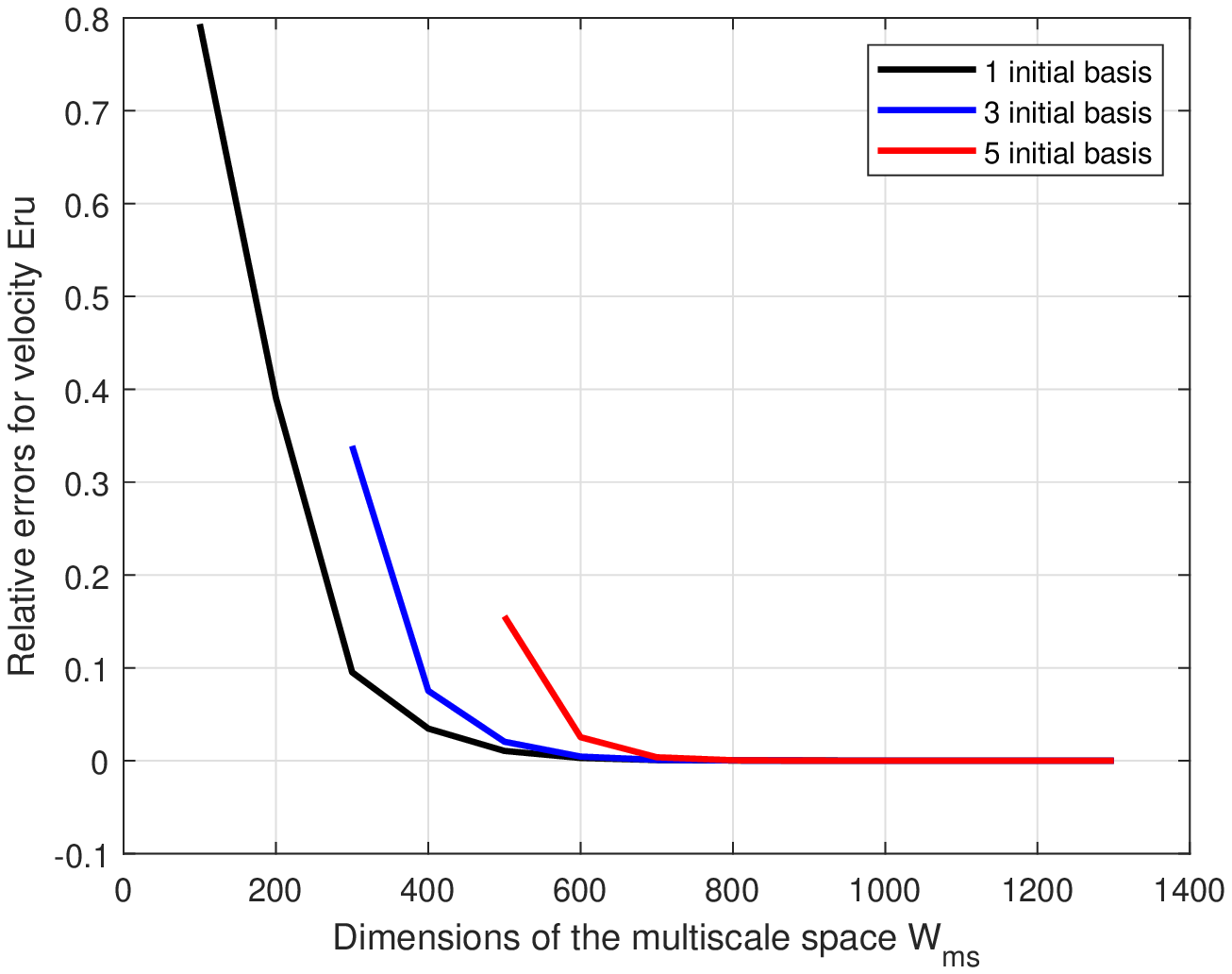}
	\end{minipage}
	\mbox{\hspace{0.01cm}}
	\begin{minipage}[b]{0.44\textwidth}
		\centering
		\includegraphics[scale=0.555]{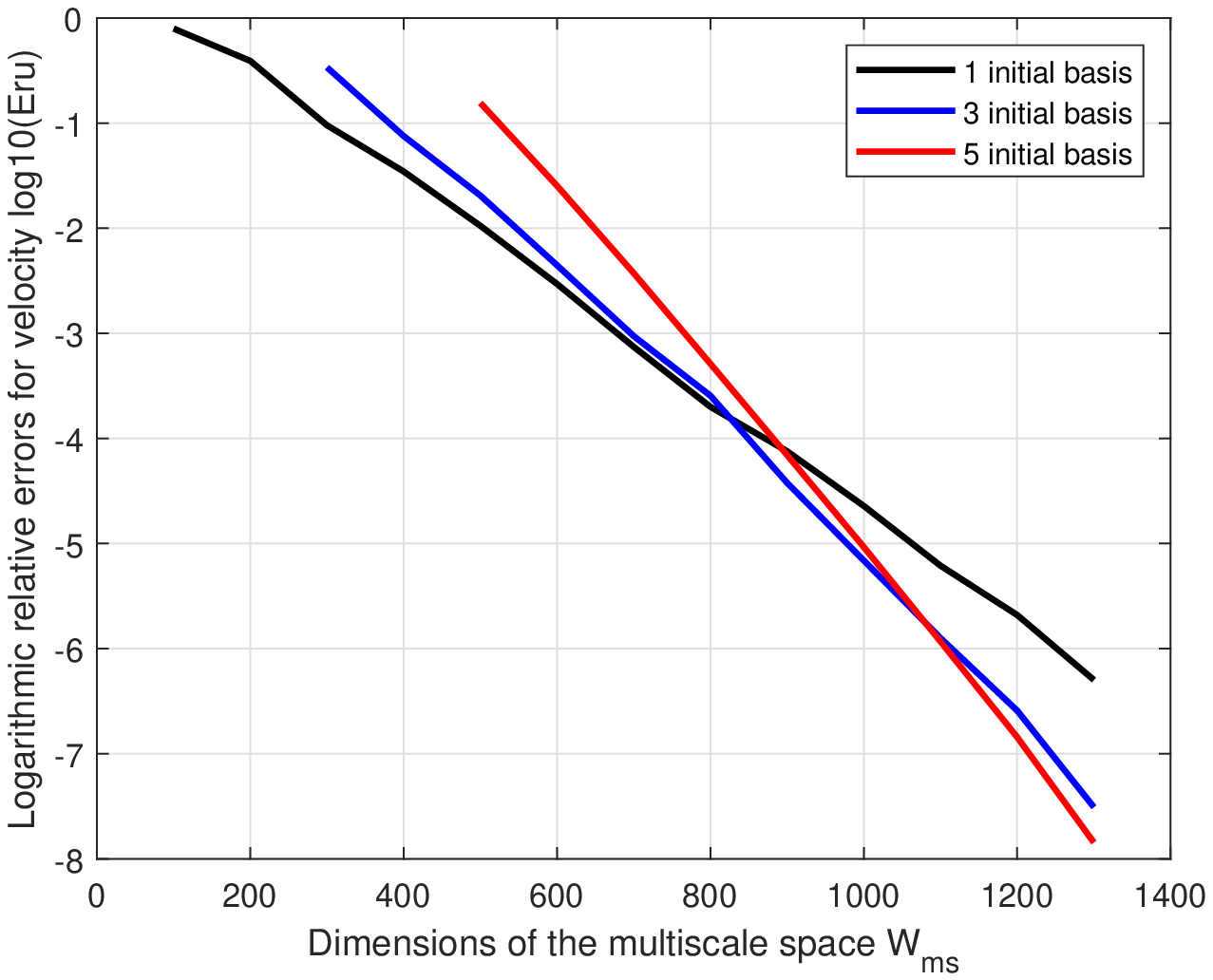}
	\end{minipage}\\
	\mbox{\hspace{0.00cm}}
	\begin{minipage}[b]{0.44\textwidth}
		\centering
		\includegraphics[scale=0.555]{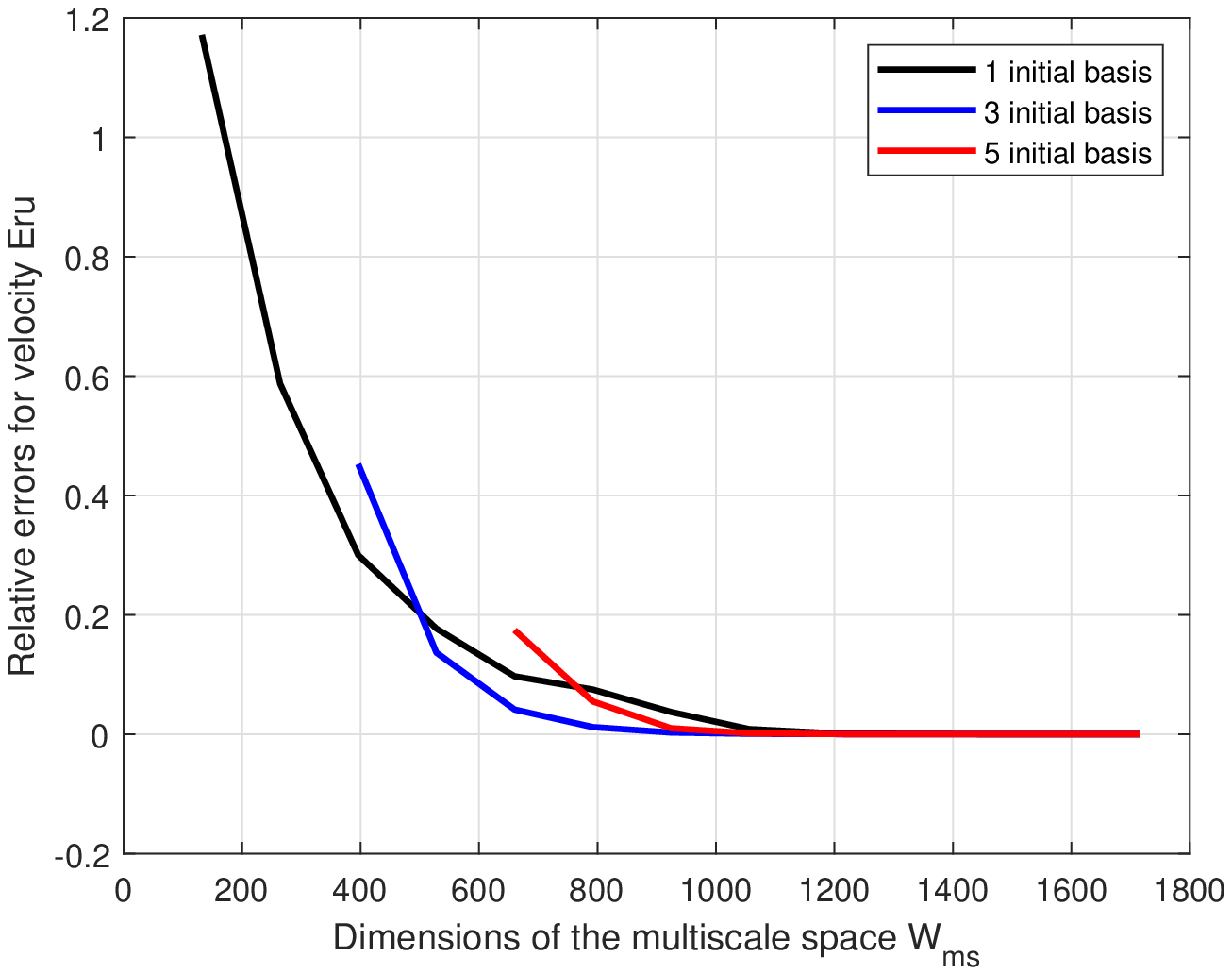}
	\end{minipage}
	\mbox{\hspace{0.01cm}}
	\begin{minipage}[b]{0.44\textwidth}
		\centering
		\includegraphics[scale=0.555]{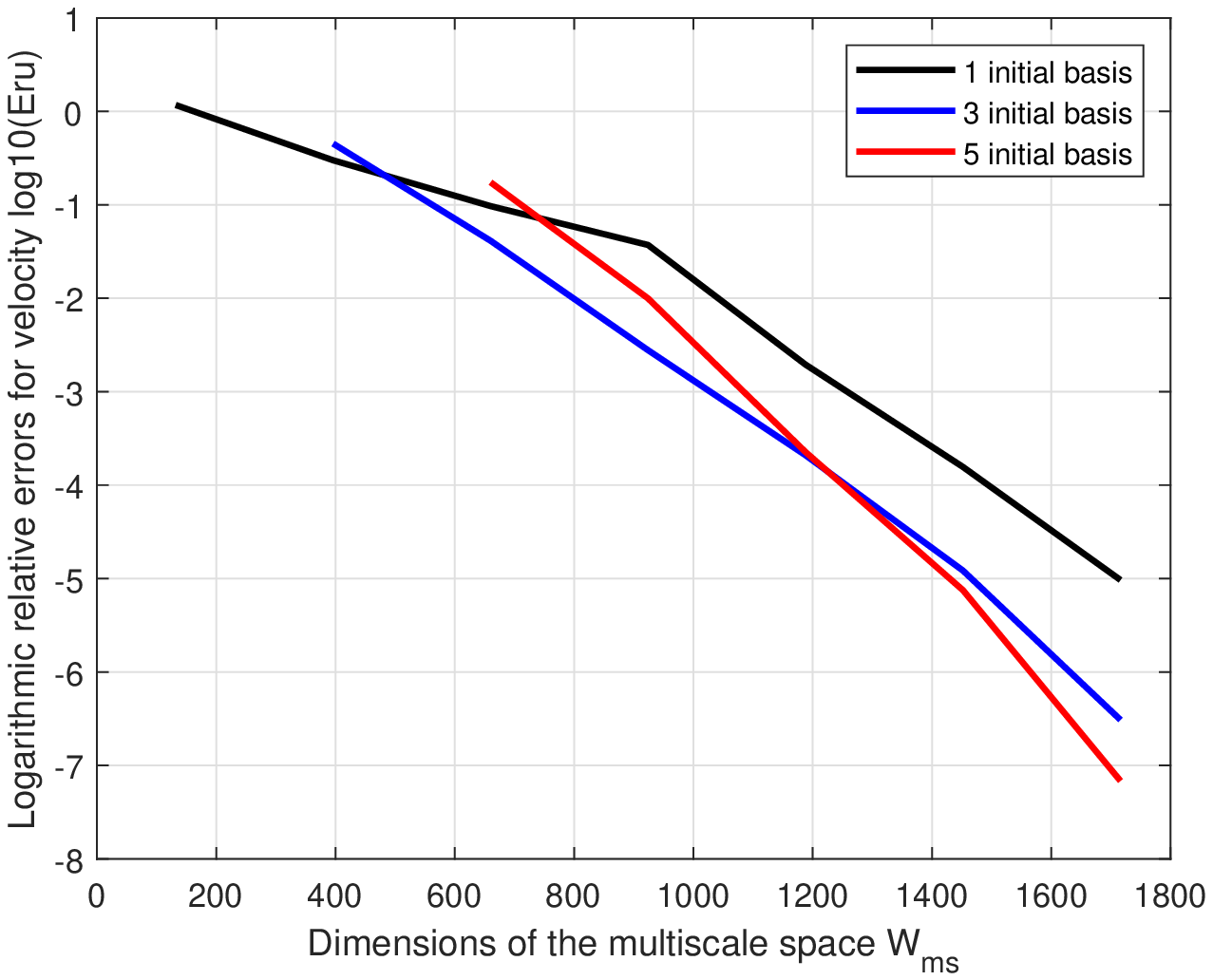}
	\end{minipage}
	\caption{(Example 3) Convergence comparisons between the usage of one, three and five initial basis functions. Top left: Relative errors for velocity using data in Example 1. Top right: Logarithmic of relative errors for velocity using data in Example 1. Bottom left: Relative errors for velocity using data in Example 2. Bottom right: Logarithmic of relative errors for velocity using data in Example 2.}\label{fig_uniform_online_velocityerr_12}
\end{figure}

{\bf Example 4:}
\begin{figure}[h!]
	\centering
	\includegraphics[width=0.49\textwidth]{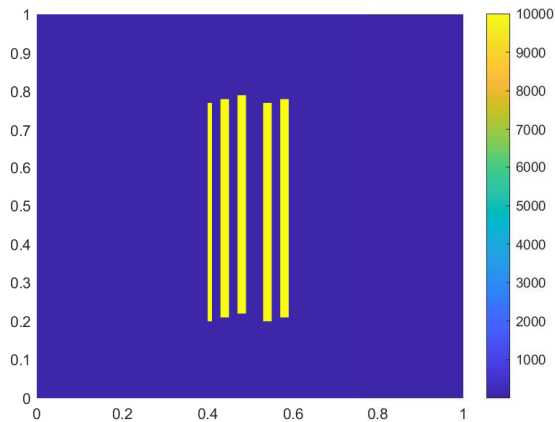}
	\caption{The distribution of permeability field $\kappa$ in Example 4.}
	\label{fig_perm_3}
\end{figure}
To further study the significance of the number of initial basis functions, we present another permeability field $\kappa$ shown in Figure \ref{fig_perm_3}. The computational domain $\Omega=[0, 1]^2$ is divided into a $100\times100$ uniform mesh for the fine grid, and a $10\times10$ uniform mesh for the coarse grid. We consider three cases with different permeability contrasts $10^2$, $10^4$ and $10^6$. In Table \ref{tab_err_ex3_online_uniform_1}, Table \ref{tab_err_ex3_online_uniform_2}, and Table \ref{tab_err_ex3_online_uniform_3}, we present the convergence histories for the use of one, two and three initial basis functions per coarse element, respectively. In these tables, from left to right, contrasts of the permeability $\kappa$ are $10^2$, $10^4$ and $10^6$, respectively, that is, conductivities of inclusions (yellow color) in Figure \ref{fig_perm_3} are $10^2$, $10^4$ and $10^6$, respectively. In Table \ref{tab_err_ex3_online_uniform_1} and Table \ref{tab_err_ex3_online_uniform_2}, we observe that the error decay for lower contrast case is faster than the higher contrast case with one and two initial basis functions. For the use of three initial basis functions shown in Table \ref{tab_err_ex3_online_uniform_3}, we observe that the relative errors decay rapidly for all contrast cases, and the convergence rates are similar, i.e., independent of the permeability contrasts, which further verified the theoretical analysis in Theorem \ref{theorem_4_1}.
\begin{table}[h!]
	\caption{(Example 4) Relative errors $\textrm{Erp}(p_{\textrm{ms}})$, $\textrm{Eru}(\mathbf{u}_{\textrm{ms}})$ of the uniform online enrichment using one initial basis function per coarse element with different contrasts. Left: Contrast $10^2$, $\Lambda_{\min}=9.4283e-4$. Middle: Contrast $10^4$, $\Lambda_{\min}=9.5935e-6$. Right: Contrast $10^6$, $\Lambda_{\min}=9.5951e-8$.}\label{tab_err_ex3_online_uniform_1}
	\centering
	\begin{tabular*}
		{0.93\textwidth}{@{\extracolsep{\fill}}|c|c|c||c|c|c||c|c|c|}
		\hline
		$\#$Dofs &$\textrm{Erp}(p_{\textrm{ms}})$ &$\textrm{Eru}(\mathbf{u}_{\textrm{ms}})$ &$\#$Dofs &$\textrm{Erp}(p_{\textrm{ms}})$ &$\textrm{Eru}(\mathbf{u}_{\textrm{ms}})$ &$\#$Dofs &$\textrm{Erp}(p_{\textrm{ms}})$ &$\textrm{Eru}(\mathbf{u}_{\textrm{ms}})$\\
		\hline
		100   &0.0344    &0.7705               &100   &0.0345    &0.7699               &100   &0.0345    &0.7697 \\
		200   &0.0220    &0.2828               &200   &0.0259    &0.3074               &200   &0.0260    &0.3077 \\
		300   &0.0196    &0.2026               &300   &0.0250    &0.2422               &300   &0.0251    &0.2426 \\
		400   &0.0105    &0.1163               &400   &0.0250    &0.2325               &400   &0.0251    &0.2332 \\
		500   &0.0039    &0.0503               &500   &0.0247    &0.2302               &500   &0.0250    &0.2323 \\
		600   &0.0013    &0.0202               &600   &0.0235    &0.2213               &600   &0.0250    &0.2322 \\
		700   &3.5209e-4 &0.0086               &700   &0.0178    &0.1876               &700   &0.0250    &0.2321 \\
		800   &1.7974e-4 &0.0059               &800   &0.0119    &0.1407               &800   &0.0249    &0.2317 \\
		\hline
	\end{tabular*}
\end{table}

\begin{table}[h!]
	\caption{(Example 4) Relative errors $\textrm{Erp}(p_{\textrm{ms}})$, $\textrm{Eru}(\mathbf{u}_{\textrm{ms}})$ of the uniform online enrichment using two initial basis functions per coarse element with different contrasts. Left: Contrast $10^2$, $\Lambda_{\min}=0.0026$. Middle: Contrast $10^4$, $\Lambda_{\min}=2.5751e-5$. Right: Contrast $10^6$, $\Lambda_{\min}=2.5747e-7$.}\label{tab_err_ex3_online_uniform_2}
	\centering
	\begin{tabular*}
		{0.93\textwidth}{@{\extracolsep{\fill}}|c|c|c||c|c|c||c|c|c|}
		\hline
		$\#$Dofs &$\textrm{Erp}(p_{\textrm{ms}})$ &$\textrm{Eru}(\mathbf{u}_{\textrm{ms}})$ &$\#$Dofs &$\textrm{Erp}(p_{\textrm{ms}})$ &$\textrm{Eru}(\mathbf{u}_{\textrm{ms}})$ &$\#$Dofs &$\textrm{Erp}(p_{\textrm{ms}})$ &$\textrm{Eru}(\mathbf{u}_{\textrm{ms}})$\\
		\hline
		200   &0.0928    &0.7461               &200   &0.0860    &0.7103               &200   &0.0757    &0.7255 \\
		300   &0.0258    &0.2936               &300   &0.0261    &0.2933               &300   &0.0285    &0.3378 \\
		400   &0.0071    &0.0985               &400   &0.0070    &0.0995               &400   &0.0108    &0.1552 \\
		500   &0.0015    &0.0277               &500   &0.0015    &0.0323               &500   &0.0040    &0.0605 \\
		600   &3.6771e-4 &0.0088               &600   &7.1039e-4 &0.0211               &600   &0.0010    &0.0232 \\
		700   &1.8341e-4 &0.0036               &700   &6.3902e-4 &0.0197               &700   &7.0677e-4 &0.0208 \\
		800   &4.2469e-5 &9.4641e-4            &800   &5.2819e-4 &0.0163               &800   &6.8788e-4 &0.0207 \\
		900   &8.9878e-6 &3.2320e-4            &900   &1.3385e-4 &0.0071               &900   &6.8080e-4 &0.0206 \\
		\hline
	\end{tabular*}
\end{table}
\begin{table}[h!]
	\caption{(Example 4) Relative errors $\textrm{Erp}(p_{\textrm{ms}})$, $\textrm{Eru}(\mathbf{u}_{\textrm{ms}})$ of uniform online enrichment using three initial basis functions per coarse element with different contrasts. Left: Contrast $10^2$, $\Lambda_{\min}=0.0232$. Middle: Contrast $10^4$, $\Lambda_{\min}=0.0223$. Right: Contrast $10^6$, $\Lambda_{\min}=0.0223$.}\label{tab_err_ex3_online_uniform_3}
	\centering
	\begin{tabular*}
		{0.93\textwidth}{@{\extracolsep{\fill}}|c|c|c||c|c|c||c|c|c|}
		\hline
		$\#$Dofs &$\textrm{Erp}(p_{\textrm{ms}})$ &$\textrm{Eru}(\mathbf{u}_{\textrm{ms}})$ &$\#$Dofs &$\textrm{Erp}(p_{\textrm{ms}})$ &$\textrm{Eru}(\mathbf{u}_{\textrm{ms}})$ &$\#$Dofs &$\textrm{Erp}(p_{\textrm{ms}})$ &$\textrm{Eru}(\mathbf{u}_{\textrm{ms}})$\\
		\hline
		300   &0.0159    &0.1674               &300   &0.0168    &0.1726               &300   &0.0168    &0.1727 \\
		400   &0.0030    &0.0414               &400   &0.0034    &0.0443               &400   &0.0034    &0.0443 \\
		500   &5.1293e-4 &0.0095               &500   &4.8449e-4 &0.0087               &500   &4.8437e-4 &0.0087 \\
		600   &9.5984e-5 &0.0024               &600   &8.4381e-5 &0.0019               &600   &8.4450e-5 &0.0019 \\
		700   &3.0699e-5 &7.7508e-4            &700   &2.5300e-5 &5.4401e-4            &700   &2.5525e-5 &4.6136e-4 \\
		800   &7.9786e-6 &2.4366e-4            &800   &4.1344e-6 &1.4041e-4            &800   &4.1418e-6 &8.5442e-5 \\
		900   &3.4821e-6 &8.8916e-5            &900   &8.7758e-7 &4.4792e-5            &900   &8.2877e-7 &2.1322e-5 \\
		1000  &9.4343e-7 &2.4935e-5            &1000  &2.0419e-7 &6.8467e-6            &1000  &2.0197e-7 &5.6335e-6 \\
		\hline
	\end{tabular*}
\end{table}
\begin{table}[h!]
	\caption{(Example 5) Relative errors at the last six iterations of the online adaptive enrichment with three initial basis functions per coarse element, $\theta_{\textrm{on}}=0.7$ and $tol=10^{-3}$. Left: Using data in Example 1. Right: Using data in Example 2.}\label{tab_err_ex12_online_adaptive_velocity}
	\centering
	\mbox{\hspace{1.30cm}}
	\begin{minipage}[b]{0.46\textwidth}
		\begin{tabular*}
			{0.82\textwidth}{@{\extracolsep{\fill}}|c|c|c|}
			\hline
			\multicolumn{3}{|c|}{Initial Dofs per $T=3$, data in Example 1} \\
			\hline
			$\#$Dofs &$\qquad\textrm{Erp}(p_{\textrm{ms}})$ &$\textrm{Eru}(\mathbf{u}_{\textrm{ms}})$\\
			\hline
			709   &9.0882e-5 &4.6765e-4  \\
			756   &4.4725e-5 &2.1259e-4  \\
			802   &1.9867e-5 &1.0231e-4  \\
			848   &8.5440e-6 &4.4013e-5  \\
			892   &3.6393e-6 &2.0677e-5  \\
			933   &1.8553e-6 &1.0259e-5  \\
			\hline
		\end{tabular*}
	\end{minipage}
    \mbox{\hspace{-0.30cm}}
    \begin{minipage}[b]{0.46\textwidth}
    	\begin{tabular*}
    		{0.82\textwidth}{@{\extracolsep{\fill}}|c|c|c|}
    		\hline
    		\multicolumn{3}{|c|}{Initial Dofs per $T=3$, data in Example 2} \\
    		\hline
    		$\#$Dofs &$\qquad\textrm{Erp}(p_{\textrm{ms}})$ &$\textrm{Eru}(\mathbf{u}_{\textrm{ms}})$\\
    		\hline
    		550   &7.7395e-3 &2.3308e-2  \\
    		612   &3.5899e-3 &1.1867e-2  \\
    		670   &2.1758e-3 &6.0973e-3  \\
    		747   &8.0973e-4 &2.9364e-3  \\
    		808   &4.0425e-4 &1.4809e-3  \\
    		874   &2.1741e-4 &7.3230e-4  \\
    		\hline
    	\end{tabular*}
    \end{minipage}
\end{table}
\begin{figure}[h!]
	\mbox{\hspace{-0.30cm}}
	\begin{minipage}[b]{0.52\textwidth}
		\centering
		\includegraphics[scale=0.36]{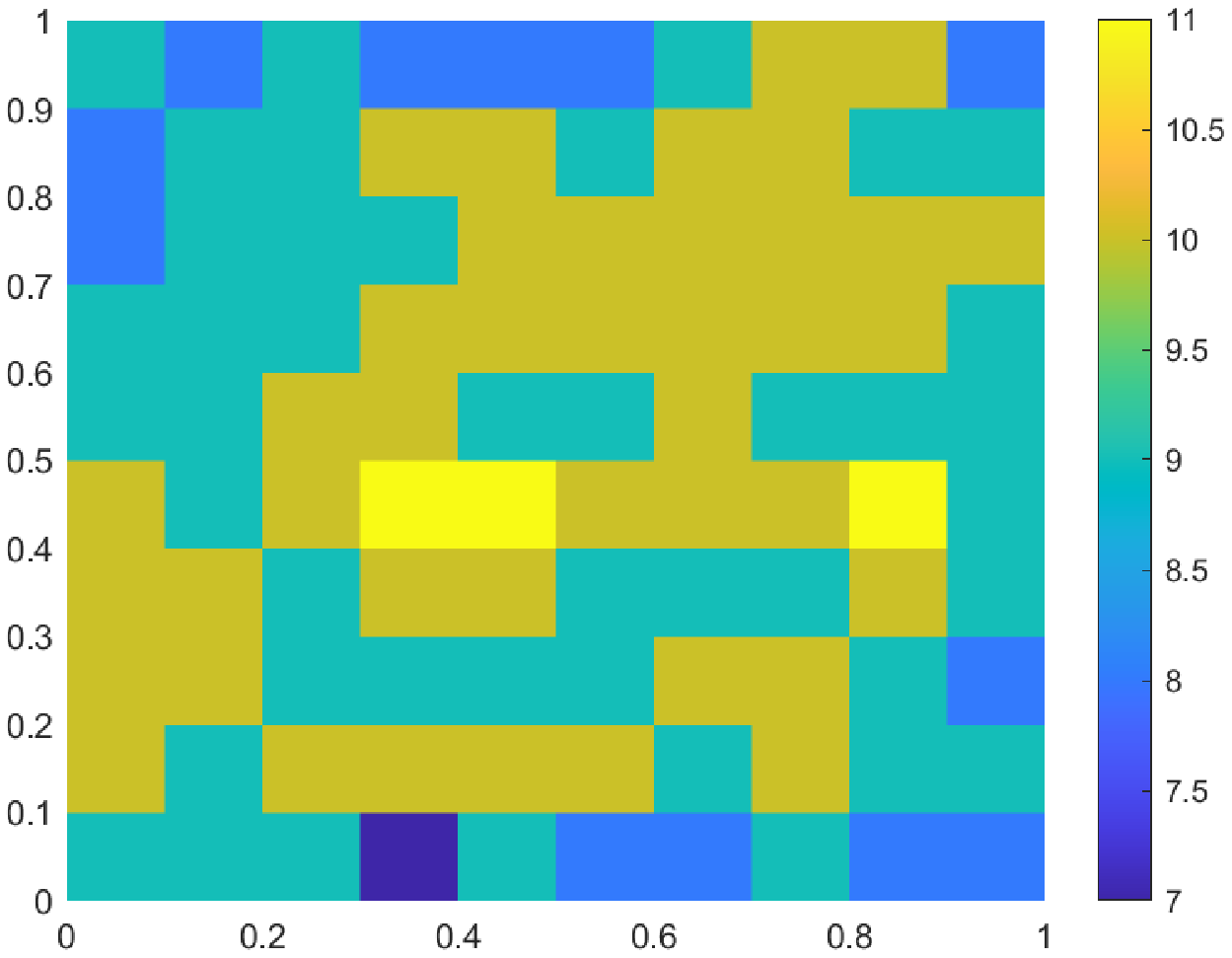}
	\end{minipage}
	\mbox{\hspace{-0.70cm}}
	\begin{minipage}[b]{0.52\textwidth}
		\centering
		\includegraphics[scale=0.36]{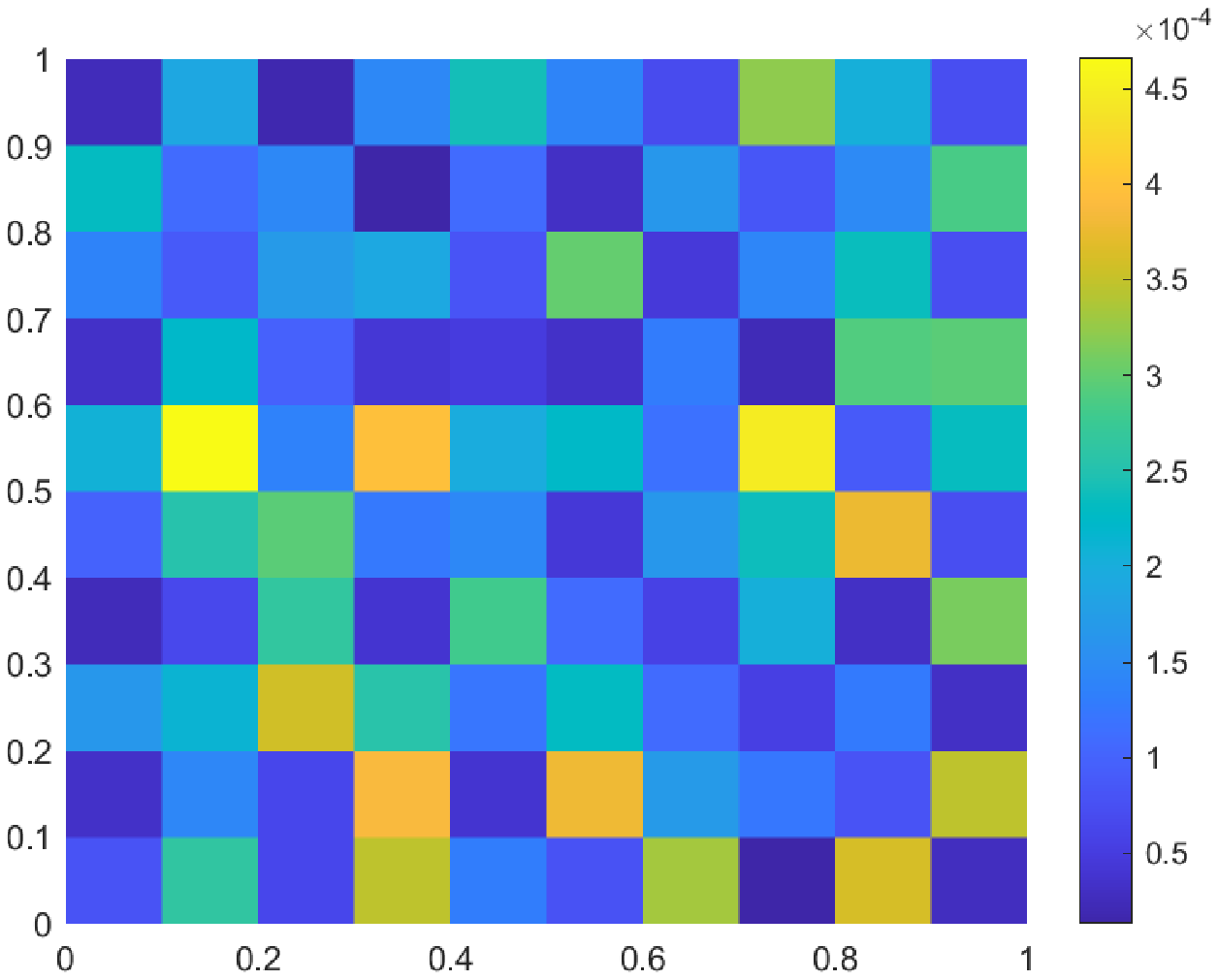}
	\end{minipage}\\
	\mbox{\hspace{-0.90cm}}
	\begin{minipage}[b]{0.52\textwidth}
		\centering
		\includegraphics[scale=0.66]{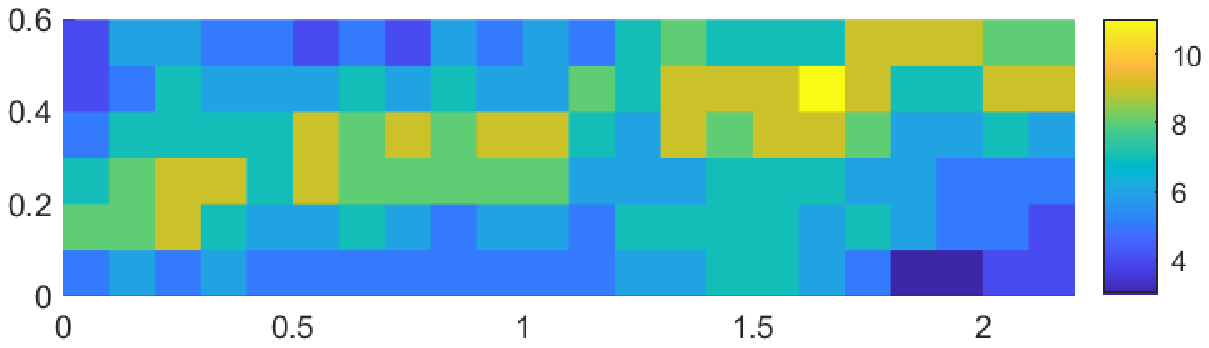}
	\end{minipage}
	\mbox{\hspace{-0.20cm}}
	\begin{minipage}[b]{0.52\textwidth}
		\centering
		\includegraphics[scale=0.66]{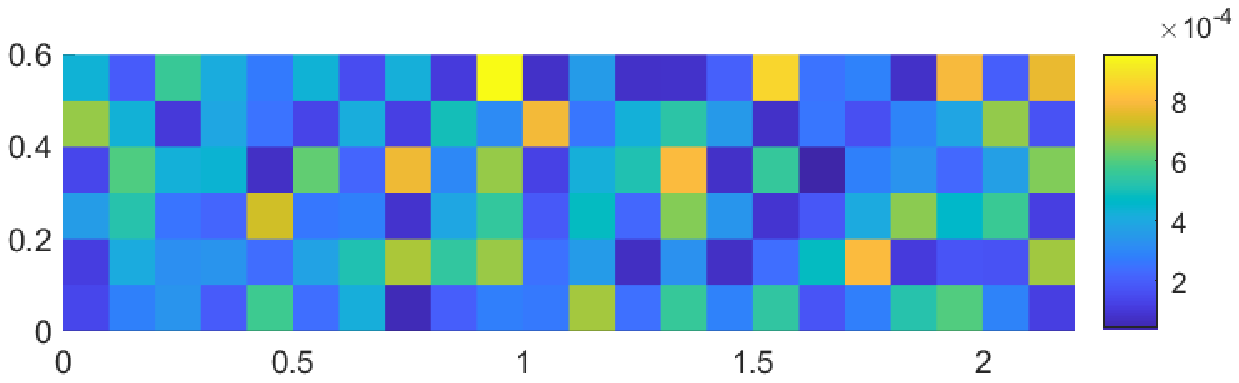}
	\end{minipage}
	\caption{(Example 5) The dimension distributions of the multiscale space and the error estimators $\eta$ after the last iteration of the online adaptive enrichment with $\theta_{\textrm{on}}=0.7$ and $tol=10^{-3}$. Top: Using data in Example 1, $\#$Dofs$=933$. Bottom: Using data in Example 2, $\#$Dofs$=874$.}
	\label{fig_num_residual_12}
\end{figure}
\begin{figure}[h!]
	\centering
	\mbox{\hspace{0.00cm}}
	\begin{minipage}[b]{0.44\textwidth}
		\centering
		\includegraphics[scale=0.56]{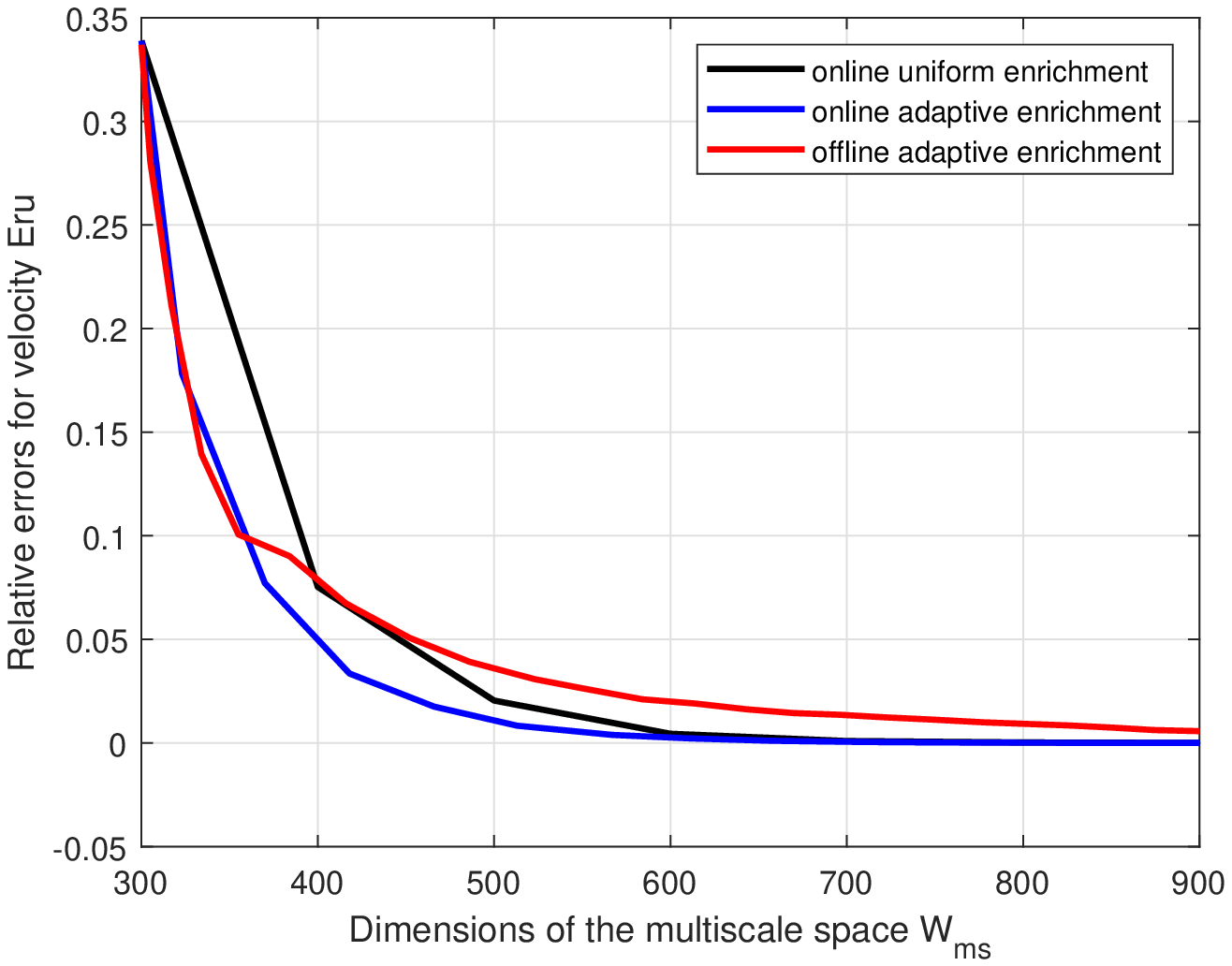}
	\end{minipage}
	\mbox{\hspace{0.02cm}}
	\begin{minipage}[b]{0.44\textwidth}
		\centering
		\includegraphics[scale=0.56]{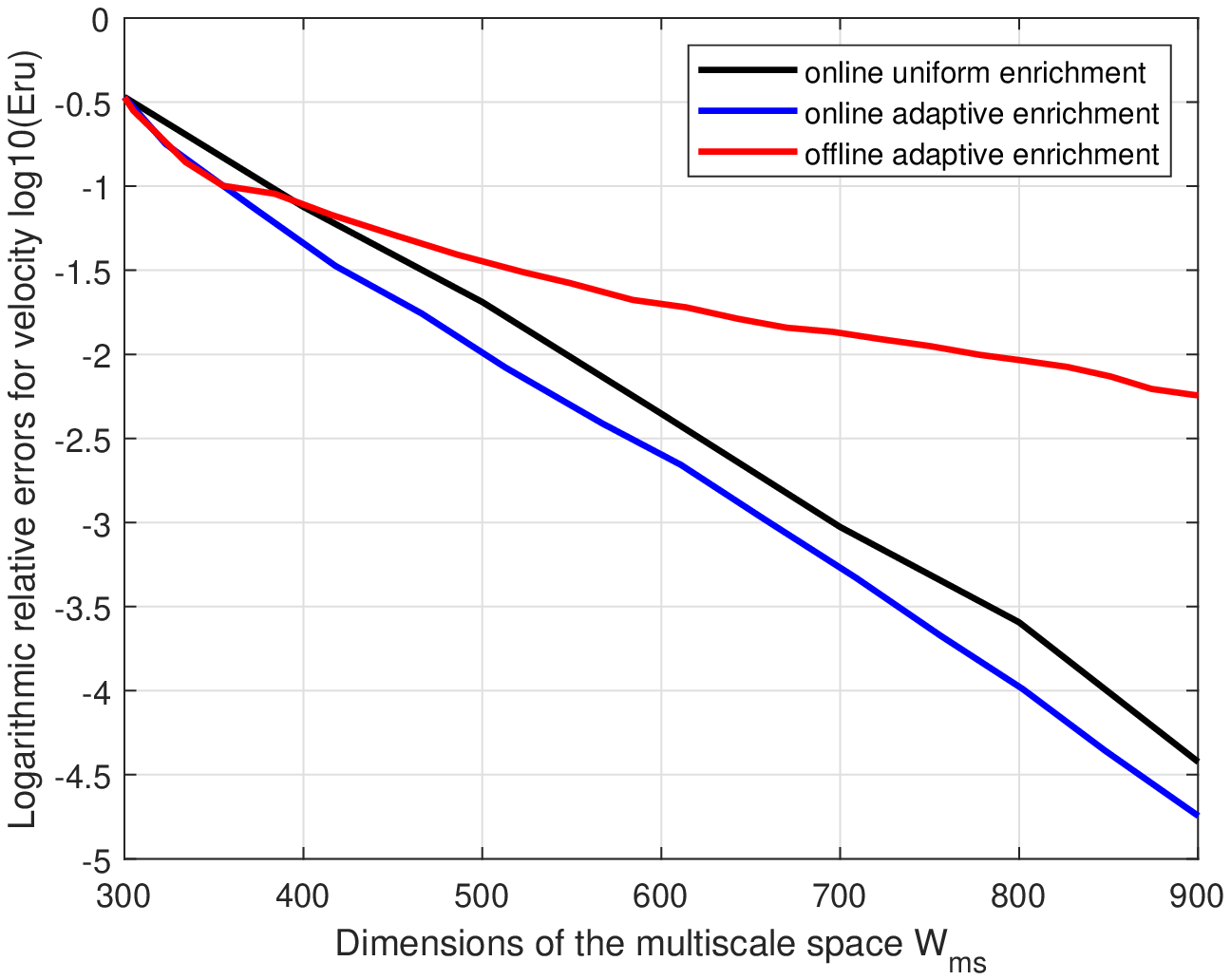}
	\end{minipage}\\
	\mbox{\hspace{0.00cm}}
	\begin{minipage}[b]{0.44\textwidth}
		\centering
		\includegraphics[scale=0.56]{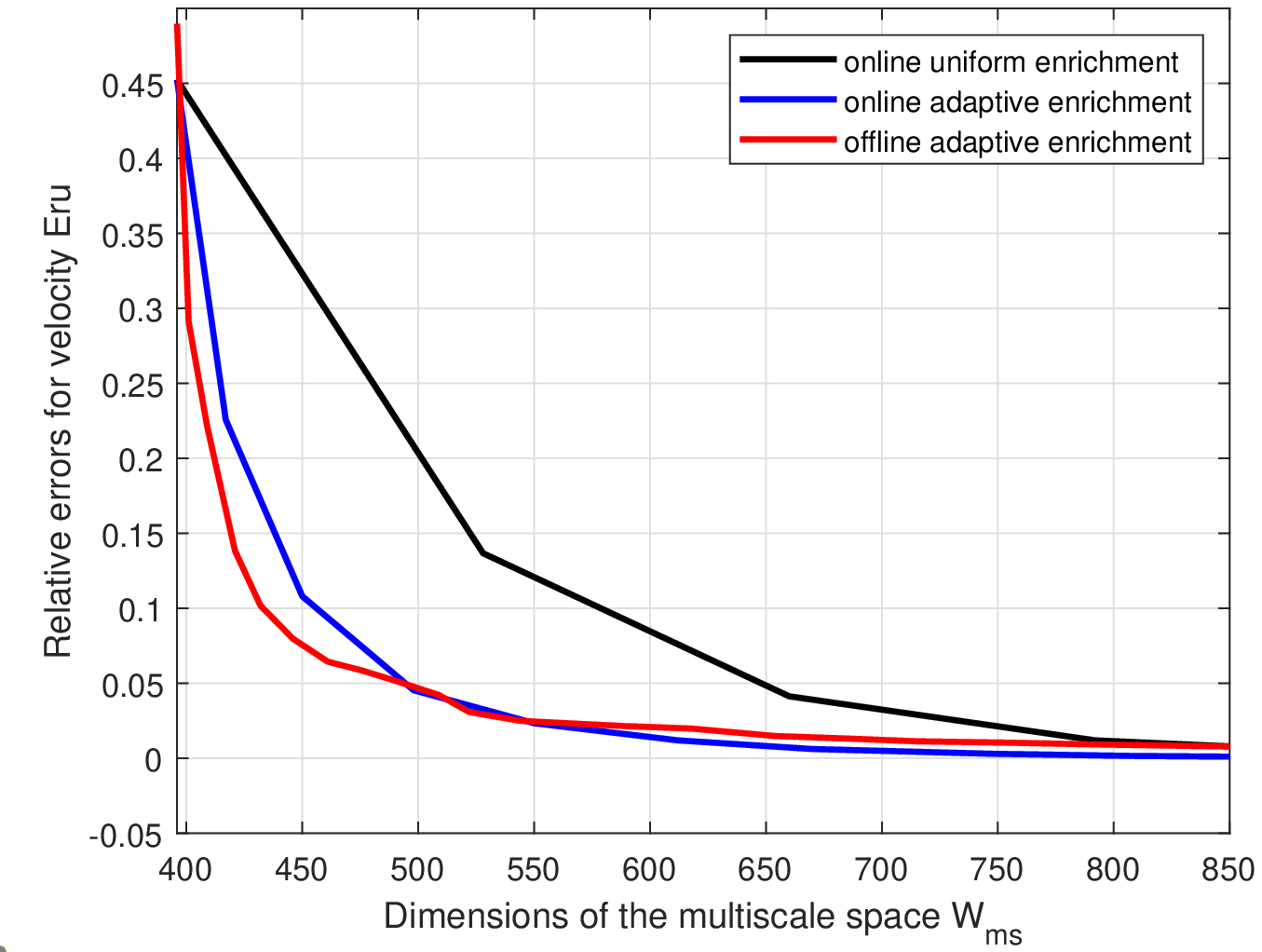}
	\end{minipage}
	\mbox{\hspace{0.02cm}}
	\begin{minipage}[b]{0.44\textwidth}
		\centering
		\includegraphics[scale=0.56]{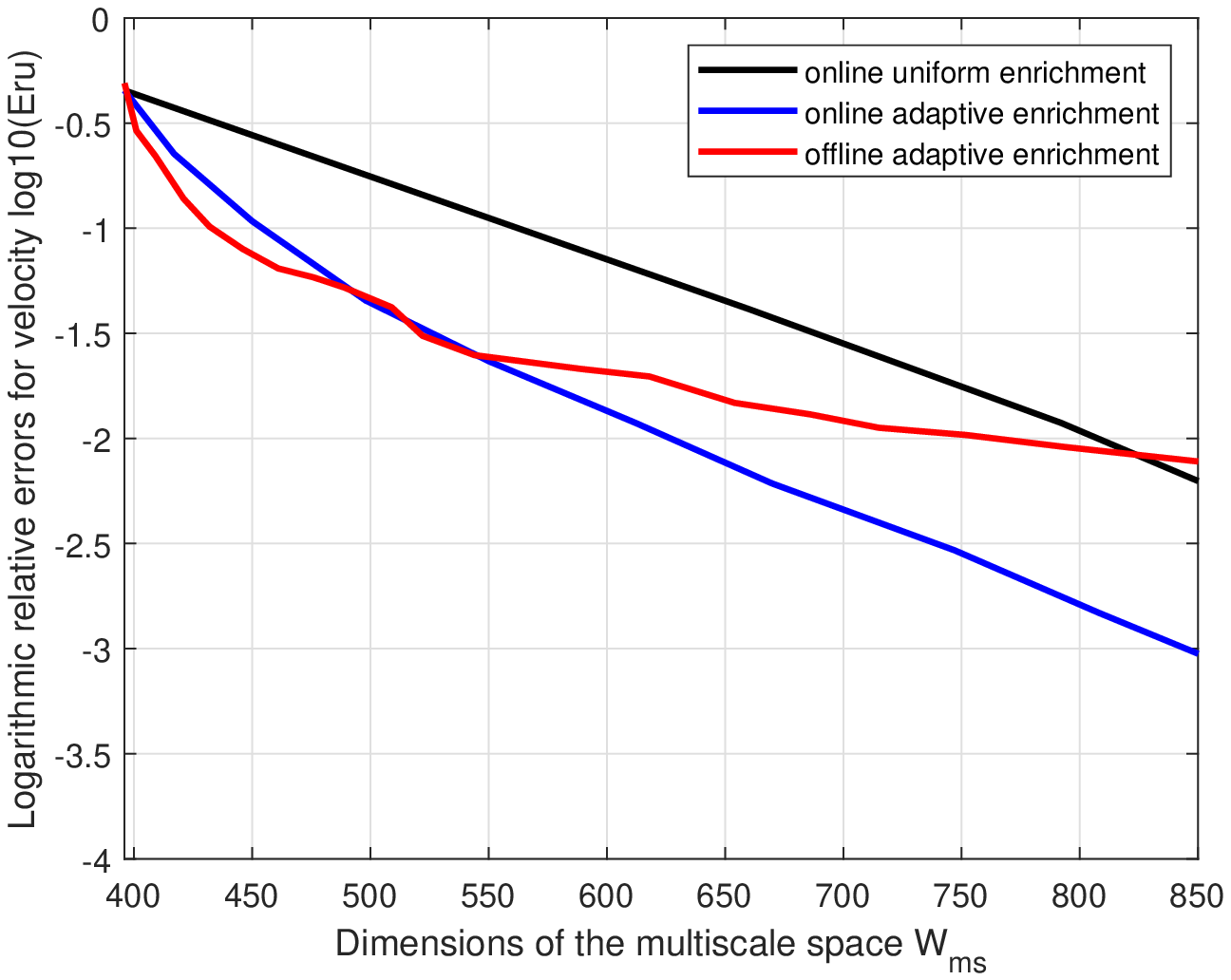}
	\end{minipage}
	\caption{(Example 5) Convergence comparisons between the online uniform enrichment, online adaptive enrichment and offline adaptive enrichment with three initial basis functions, $\theta_{\textrm{off}}=0.7$ and $\theta_{\textrm{on}}=0.7$. Top left: Relative errors for velocity using data in Example 1. Top right: Logarithmic of relative errors for velocity using data in Example 1. Bottom left: Relative errors for velocity using data in Example 2. Bottom right: Logarithmic of relative errors for velocity using data in Example 2.}\label{fig_uniform_adaptive_velocityerr_12}
\end{figure}

Next, we perform the online adaptive enrichment by computing and adding online basis functions only on coarse elements where the the corresponding residuals are large enough according to (\ref{eqn_criterion_online}). We will stop the enrichment until the local estimator (\ref{eqn_online_adaptive_estimator}) reaches a certain threshold $\eta_i\le tol$, $i=1,2,\cdots,N_T$, where $tol$ is the specified tolerance.

{\bf Example 5:} In this example, we take parameters $\theta_{\textrm{on}}=0.7$ and $tol=10^{-3}$. As Example 3, we use the data in Example 1 and Example 2 for the numerical tests. In Table \ref{tab_err_ex12_online_adaptive_velocity}, we present numerical results for the last six iterations of the online adaptive enrichment using three initial basis functions per coarse element. In Figure \ref{fig_num_residual_12}, we present distributions of the number of basis functions and local error estimators after the last iteration of the online adaptive enrichment. We find that the local error estimator is smaller than $tol=10^{-3}$ on each coarse element, the distribution of the number of basis functions are different form the resulting distributions of the offline adaptive enrichment with offline basis functions shown in Figure \ref{fig_num_1} and Figure \ref{fig_num_2}. It seems that the number distributions of basis functions associated with online adaptive enrichment are more average than the offline adaptive enrichment. The convergence histories of the online adaptive enrichment and online uniform enrichment are plotted in Figure \ref{fig_uniform_adaptive_velocityerr_12}, where we plots the relative errors for the velocity against the dimensions of the multiscale space $W_{\textrm{ms}}$, we find that the performance of online adaptive enrichment is better than the online uniform enrichment with the same dimensions of the multiscale space. In addition, the convergence history of offline adaptive enrichment also plotted in the red line with parameter $\theta_{\textrm{off}}=0.7$ and three initial basis functions per coarse element, we observe that at beginning iterations of the multiscale space enrichment, the performances of adaptive enrichment with online basis functions and offline basis functions are similar, however, after several iterations of the enrichment, when the multiscale space has sufficient basis functions, the adaptive enrichment with online basis functions becomes more effective than the adaptive enrichment with offline basis functions.
\section{Conclusions}
In this paper, we develop offline and online adaptive methods, respectively, to enrich the multiscale space for the generalized multiscale approximation of a mixed finite element method with velocity elimination. We derive an a-posteriori error indicator depending on the pressure-related weighted $L^2$-norm of the local residual operator, where the eigenvalue structures of spectral decompositions in the offline stage are also coupled into the error indicator. Based on the proposed error indicator, we present the offline adaptive method to enrich the multiscale space by adding offline basis functions iteratively on coarse elements with large local residuals, where offline multiscale basis functions are computed in the offline stage before the enrichment. We also propose the online adaptive method that makes use of online basis functions for the multiscale space enrichment on selected coarse elements relying on the velocity-based weighted $L^2$-norm of the local residual. Online basis functions are calculated in the actual simulation based on the solution of the previous iteration and some optimal minimum energy principles. We give the theoretical analysis for the convergence of both these two adaptive methods, the analysis shows that sufficient initial basis functions leads to a faster convergence rate. Moreover, we have conducted plenty of numerical examples to demonstrate the performance of these two adaptive methods and also confirm the theoretical analysis. We find that both the offline and online adaptive methods have competitive performances that can achieve higher accuracies compared with the offline uniform enrichment using the same number of offline basis functions. Besides, when the multiscale space has adequate basis functions, the online adaptive method generally performs better than the offline adaptive method as online basis functions contain important global information such as distant effects that offline basis functions cannot capture. In particular, when the initial multiscale space contains all offline basis functions corresponding to the relative smaller eigenvalues (which are contrast sensitive) of the local spectral decompositions in the offline stage, the convergence rate of the online enrichment is independent of the contrast of the permeability.
\bibliographystyle{elsarticle-num}
\bibliography{Darcy_GMsFEM_Adaptive_REF}
\end{document}